\documentclass{amsart}
\usepackage{amsmath} 
\usepackage{amssymb}
\usepackage{mathrsfs}

\newtheorem{theorem}{Theorem}[section] 
\newtheorem{claim}[theorem]{Claim}
\newtheorem{lemma}[theorem]{Lemma} 
 
\newtheorem{conclusion}[theorem]{Conclusion}
\newtheorem{observation}[theorem]{Observation}

\theoremstyle{definition}
\newtheorem{definition}[theorem]{Definition}
\newtheorem{convention}[theorem]{Convention}

\newtheorem{fact}[theorem]{Fact}
\newtheorem{discussion}[theorem]{Discussion}

\theoremstyle{remark}
\newtheorem{remark}[theorem]{Remark}

\newtheorem{notation}[theorem]{Notation}
\newtheorem{context}[theorem]{Context}
\newtheorem{content}[theorem]{Content}

\newcommand{\Dp}{{\rm Dp}}

\newcommand{\bas}{{\rm bas}}
\newcommand{\can}{{\rm can}}
\newcommand{\lex}{{\rm lex}}
\newcommand{\wed}{{\rm wd}}
\newcommand{\vhm}{{\rm vhm}}
\newcommand{\glb}{{\rm glb}}
\newcommand{\nic}{{\rm nice}}

\newcommand{\Th}{{\rm Th}}
\newcommand{\Rk}{{\rm Rk}}
\newcommand{\rk}{{\rm rk}}
\newcommand{\Sk}{{\rm Sk}}
\newcommand{\tp}{{\rm tp}}
\newcommand{\tr}{{\rm tr}}
\newcommand{\cb}{{\rm cb}}

\newcommand{\pp}{{\rm pp}}
\newcommand{\cl}{{\rm cl}}
\newcommand{\all}{{\rm all}}

\newcommand{\ged}{{\rm gd}}
\newcommand{\Sep}{{\rm Sp}}
\newcommand{\Dep}{{\rm Dp}}
\newcommand{\res}{{\rm res}}
\newcommand{\otp}{{\rm otp}}
\newcommand{\inv}{{\rm inv}}
\newcommand{\pcf}{{\rm pcf}}
\newcommand{\Reg}{{\rm Reg}}
\newcommand{\Max}{{\rm Max}}
\newcommand{\GCH}{{\rm GCH}}

\newcommand{\complete}{{\rm complete}}

\newcommand{\Succ}{{\rm Succ}}
\newcommand{\Suc}{{\rm Suc}}
\newcommand{\ZFC}{{\rm ZFC}}
\newcommand{\rang}{{\rm rang}}
\newcommand{\Range}{{\rm Range}}
\newcommand{\splitt}{{\rm split}}
\newcommand{\IRCar}{{\rm IRCar}}

\newcommand{\Ord}{{\rm Ord}}

\newcommand{\bs}{{\rm bs}}
\newcommand{\bd}{{\rm bd}}

\newcommand{\cov}{{\rm cov}}
\newcommand{\tcf}{{\rm tcf}}

\newcommand{\Min}{{\rm Min}}
\newcommand{\Dom}{{\rm Dom}}
\newcommand{\Rang}{{\rm Rang}}

\newcommand{\rest}{{\restriction}}
\newcommand{\dom}{{\rm dom}}

\newcommand{\wilog}{{\rm without loss of generality}}
\newcommand{\Wilog}{{\rm Without loss of generality}}

\newcommand{\then}{{\underline{then}}}
\newcommand{\when}{{\underline{when}}}
\newcommand{\Then}{{\underline{Then}}}

\newcommand{\If}{{\underline{if}}}
\newcommand{\Iff}{{\underline{iff}}}
\newcommand{\mn}{{\medskip\noindent}}
\newcommand{\sn}{{\smallskip\noindent}}

\newcommand{\vare}{{\varepsilon}}

\newcommand{\bbB}{{\mathbb B}}
\newcommand{\bbI}{{\mathbb I}}
\newcommand{\gb}{{\mathfrak b}}
\newcommand{\gB}{{\mathfrak B}}
\newcommand{\gc}{{\mathfrak c}}
\newcommand{\ga}{{\mathfrak a}}

\newcommand{\cE}{{\mathscr E}}
\newcommand{\cD}{{\mathscr D}}

\newcommand{\cH}{{\mathscr H}}

\newcommand{\cF}{{\mathscr F}}
\newcommand{\cG}{{\mathscr G}}
\newcommand{\cI}{{\mathscr I}}

\newcommand{\bbL}{{\mathbb L}}

\newcommand{\bbP}{{\mathbb P}}
\newcommand{\cP}{{\mathscr P}}

\newcommand{\cS}{{\mathscr S}}
\newcommand{\cT}{{\mathscr T}}
 
\newcommand{\cU}{{\mathscr U}}

\newcommand{\cf}{{\rm cf}}
\newcommand{\pr}{{\rm pr}}

\newcount\skewfactor
\def\mathunderaccent#1#2 {\let\theaccent#1\skewfactor#2
\mathpalette\putaccentunder}
\def\putaccentunder#1#2{\oalign{$#1#2$\crcr\hidewidth
\vbox to.2ex{\hbox{$#1\skew\skewfactor\theaccent{}$}\vss}\hidewidth}}
\def\name{\mathunderaccent\tilde-3 }

\newenvironment{PROOF}[2][\proofname.]
   {\begin{proof}[#1]\renewcommand{\qedsymbol}{\textsquare$_{\rm #2}$}}
   {\end{proof}}

\begin{document}

\title {Combinatorial Background for Non-structure}
\author {Saharon Shelah}
\address{Einstein Institute of Mathematics\\
Edmond J. Safra Campus, Givat Ram\\
The Hebrew University of Jerusalem\\
Jerusalem, 91904, Israel\\
 and \\
 Department of Mathematics\\
 Hill Center - Busch Campus \\ 
 Rutgers, The State University of New Jersey \\
 110 Frelinghuysen Road \\
 Piscataway, NJ 08854-8019 USA}
\email{shelah@math.huji.ac.il}
\urladdr{http://shelah.logic.at}
\thanks{In reference like \cite[1.16=L7.7]{Sh:331}, the 1.16 is
  the number of claim (or definition) and L7.7 is its label; so
  intended just to help the author to correct it if the number will be
  changed. The author thanks Alice Leonhardt for the beautiful typing. Number 
E62}



% This is E62
% was: e-appendix, appendix  of nonstructure.  (mg, Jan 2010)

%\input t3totex.tex 
% \input macros.sty 
%\input moredefs.sty % \input xatexdef.sty
%\documentclass{amsart} 
%\ifx\shlhetal\undefinedcontrolsequence\let\shlhetal\relax\fi
%\ifx \shlhetal\relax
%      \def\publicprivatemode{PUBLIC}\else
% \usepackage{amssymb}
%  \usepackage{amsthm}
%\let\frak\mathfrak
%\input  setup.sty
%\input atalya.sty
% \input markers.sty
% \renewcommand{\paperno}{app}

% Converted to Alice style 1/4/2013 - someone else typed originally

% previous version:  Dec/23/2014

\date{June 24, 2015}

\begin{abstract}
This was supposed to be an appendix to the book {\it Non-structure}, 
and probably will be if it materializes.

It presents relevant material sometimes new, which used in works
which were supposed to be part of that book.

In \S1 we deal with partition theorems on trees with $\omega$ levels;
it is self contained.
In \S2 we deal with linear orders which are countable union
of scattered ones with  unary predicated, it is self contained.
In \S3 we deal mainly with pcf theory but just quote.  
In \S4, on normal ideals, we repeat \cite{Sh:247}.  This 
is used in \cite{Sh:331}.
\end{abstract}

\maketitle
\numberwithin{equation}{section}
\newpage

%\centerline{Correction AT 9/02, 12/02,8/03, 1/04, 3/07}

\section{Partitions on trees}

See \cite{RuSh:117}, \cite[2.4,2.5]{Sh:136}, 
\cite{Sh:b}, \cite[XI3.5,XI3.5A,XI3.7,XI5.3,XV2.6,XV2.6A,
XV2.6B,XV2.6C]{Sh:f} on those theorems.

See Rubin, Shelah \cite{RuSh:117} pp 47-48 on the history of such
theorems, more in \cite{Sh:136}.

\begin{definition}
\label{a4}
1) $\bold I$ is an ideal on $S$ when it is a family of subsets of $S$
including the singleton, closed under union of two and $S \notin S$.

\noindent
2) An ideal $\bold I$ is $\lambda$-complete \If \, any union of less
than $\lambda$ members of $\bold I$ is still a member of $\bold I$.
\end{definition}

\noindent
In \cite[1.1=L1.1]{Sh:511}, \cite[1.2=L1.2]{Sh:511} we use
\begin{definition}
\label{p2.1}
1)  A tagged tree is a pair $({\cT},\bar{\bold I})$ such that:
\mn
\begin{enumerate}
\item[$(a)$]   ${\cT}$ is a $\omega$-tree, which in this section 
means a non-empty set of finite sequences of ordinals such that 
if $\eta \in {\cT}$ then any initial segment of $\eta$ belongs
to ${\cT}$.  We understand that ${\cT}$ is ordered by initial segments,
i.e., $\eta \le_{\cT} \nu$ means $\eta$ is an initial segment of $\nu$ that is
$\eta \trianglelefteq \nu$
\sn
\item[$(b)$]  $\bar{\bold I}$ is a function but only $\bar{\bold I} \rest
(\Dom(\bold I)\cap {\cT})$ matters, such that for every $\eta \in
  {\cT}$: if $\bar{\bold I}(\eta) = \bold I_\eta$ is defined
\then \, $\bar{\bold I}(\eta)$ is an ideal of subsets of some set called
the domain of $\bold I_\eta,\Dom(\bold I_\eta)$ and $\Dom(\bold I_\eta) 
\notin \bold I_\eta$, and
\[
\Succ_{\cT}(\eta) := \{\nu:\nu \text{ is an immediate successor of }
\eta \text{ in }{\cT}\} \subseteq \Dom(\bold I_\eta).
\]
The interesting case is when $\Succ_{\cT}(\eta) \notin \bold I_\eta$ and 
 usually $\bold I_\eta$ is $\aleph_2$-complete
\sn
\item[$(c)$]  For every $\eta \in {\cT}$ we have $\Succ_{\cT}(\eta) \ne 
\emptyset$.
\end{enumerate}
\mn
2)  We call $(\cT,\bar{\bold I})$ normal \when \, for every $\eta \in \Dom
(\bold I_\eta)$ we have: $\Dom(\bold I_\eta) = \Succ_{\cT}(\eta)$.
\end{definition}

\begin{convention}
\label{p2.1A}
1) For any tagged tree $({\cT}, \bar{\bold I})$ we can define
the function $\bar{\bold I}^\dagger,$ by:

\[
\Dom(\bar{\bold I}^\dagger) = \{\eta:\eta\in \Dom(\bar{\bold I}) 
\text{ and } \Succ_{{\cT}}(\eta) \subseteq \Dom(\bold I_\eta),\text{
  and } \Succ_{{\cT}}(\eta) \notin \bold I_\eta\}
\]

\[
\bold I^\dagger_\eta = \{\{\alpha:\eta \char 94 \langle \alpha
\rangle \in A\}: A \in \bold I_\eta\}.
\]

\mn
2) We sometimes, in an abuse of notation, do not distinguish between 
$\bar{\bold I}$ and $\bar{\bold I}^\dagger$.  
Also if $\bold I^\dagger_\eta$ is constantly $\bold I^*$, we may 
write $\bold I^*$ instead of $\bar{\bold I}$.

\noindent
3) We use ${\cT}$ only to denote $\omega$--trees.
\end{convention}

\begin{definition}
\label{p2.2}
1)  We say that $\eta$ is a splitting point of $({\cT}, \bar{\bold
  I})$ \when \, $\eta \in \cT,\bold I_\eta$ is defined and 
$\Succ_{\cT}(\eta) \notin \bold I_\eta$.
Let $\splitt({\cT},\bar{\bold I})$ be the set of splitting points of 
$({\cT},\bar{\bold I})$.  Usually, we will be interested only in trees 
where each branch meets $\splitt({\cT},\bar{\bold I})$ infinitely often. 

\noindent
2) For $\eta \in {\cT}$, let ${\cT}^{[\eta]} := \{\nu \in {\cT}:
\nu = \eta$ or $\nu \triangleleft \eta$ or $\eta \triangleleft \nu\}$. 
\end{definition}

\begin{definition}
\label{p2.3}
We now define several orders between tagged trees:

\noindent
1) $({\cT}_1, \bar{\bold I}_1) \le ({\cT}_2, \bar{\bold I}_2)$ \If \, 
${\cT}_2 \subseteq {\cT}_1$, and $\splitt({\cT}_2,\bar{\bold I}_2) 
\subseteq \splitt({\cT}_1,\bar{\bold I}_1)$, and for every 
$\eta \in \splitt({\cT}_2,\bar{\bold I}_2)$ we have $\bar{\bold
  I}_2(\eta) \rest \Succ_{{\cT}_2}(\eta) = \bar{\bold I}_1(\eta) \rest 
\Succ_{\cT_2}(\eta)$ (where $\bold I \rest A = \{B:B \subseteq A$ and 
$B \in \bold I\}$).
(So every splitting point of ${\cT}_2$ is a splitting point of 
$({\cT}_1,\bar{\bold I}_1)$, and $\bar{\bold I}_2 \rest
\text{ split}(\cT_2,\bar{\bold I}_2)$ is completely determined by 
$\bar{\bold I}_1$ and $\splitt({\cT}_2,\bar{\bold I}_2)$ provided 
that $\bar{\bold I}_2$ is normal, see \ref{p2.1}(2).)

\noindent
2) $({\cT}_1, \bar{\bold I}_1 ) \le^* ({\cT}_2,\bar{\bold I}_2)$
\when \, $({\cT}_1, \bar{\bold I}_1 ) \le ({\cT}_2,
\bar{\bold I}_2)$ and $\splitt({\cT}_2,\bar{\bold I}_2) =
\splitt({\cT}_1,\bar{\bold I}_1) \cap {\cT}_2$.

\noindent
3) $({\cT}_1,\bar{\bold I}_1) \le^{\otimes} ({\cT}_2,\bar{\bold
  I}_2)$ \If \, $({\cT}_1,\bar{\bold I}_1) \le^* ({\cT}_2,\bar{\bold
  I}_2)$ and $\eta \in {\cT}_2 \setminus \splitt({\cT}_1,\bar{\bold
  I}_1) \Rightarrow \Succ_{{\cT}_2}(\eta) = \Succ_{{\cT}_1}(\eta)$.

\noindent
4) $({\cT}_1,\bar{\bold I}_1) \le^\otimes_\mu ({\cT}_2,\bar{\bold I}_2)$ 
\If \, $({\cT}_1,\bar{\bold I}_1) \le^* ({\cT}_2,\bar{\bold I}_2)$ 
and $\eta \in {\cT}_2$ and $|\Succ_{{\cT}_1}(\eta)| < \mu 
\Rightarrow \Succ_{{\cT}_2}(\eta) = \Succ_{{\cT}_1}(\eta)$.
\end{definition}

\begin{definition}
\label{p2.4}
1)  For a set $\bbI$ of ideals, a tagged tree $({\cT},\bar{\bold I})$ is an
$\bbI$-tree \If \, for every splitting point $\eta \in {\cT}$ we have
$\bold I_\eta \in \bbI$ (up to an isomorphism) or just $\bold I_\eta$
is isomorphic to $\bold I \rest A$ for some $\bold I \in \bbI,A
\subseteq \dom(\bold I),A \notin \bold I$; but we usually use
restriction-closed $\bbI$, see Definition \ref{p2.5}(2).

\noindent
2) For a set $\bold S$ of regular cardinals, an $\bold S$-tree ${\cT}$
is a tree such that for any point $\eta \in {\cT}$ we have: 
$|\Succ_{\cT}(\eta)| \in \bold S$ or $|\Succ_{\cT}(\eta)| = 1$. 

\noindent
3) We may omit $\bar{\bold I}$ and denote a tagged tree $({\cT}, 
\bar{\bold I})$ by ${\cT}$ whenever ${\cT} \subseteq \Dom(\bar{\bold
  I})$ and $\bold I_\eta = \{A \subseteq \Succ_T(\eta):|A| < 
|\Succ_{\cT}(\eta)|\}$ and $|\Succ_{\cT}(\eta)| \in \IRCar \cup \{1\}$
for every $\eta \in {\cT}$, recalling $\IRCar$ is the class of infinite
regular cardinals.

\noindent
4) For a tree ${\cT},\lim({\cT})$ is the set of branches
of ${\cT}$, i.e. all $\omega$-sequences of ordinals, such that every
finite initial segment of them is a member of ${\cT}$, that is
$\lim({\cT}) = \{\eta \in {}^\omega\Ord:(\forall n)\,\eta \rest n\in \cT\}$.

\noindent
5)  A subset $J$ of a tree ${\cT}$ is a front if: $\eta \ne \nu
\in J$ implies none of them is an initial segment of the other, and every
$\eta \in \lim({\cT})$ has an initial segment which is a member of $J$.

\noindent
6) $({\cT},\bar{\bold I})$ is standard \If \, for 
every non-splitting point $\eta \in {\cT}$ we have $|\Succ_{\cT}(\eta)| = 1$.

\noindent
7) $({\cT},\bar{\bold I})$ is full \If \, every $\eta \in {\cT}$ is 
a splitting point.

\noindent
8) The natural topology on $\lim(I)$ for an $\omega$-tree ${\cT}$ is 
defined by ${\cU} \subseteq \lim({\cT})$ is open \when \,
for every $\eta \in {\cU}$ for some $n<\omega$ we have 
$\lim({\cT}^{[\eta \rest n]}) \subseteq {\cU}$.
\end{definition}

\noindent
Recall
\begin{observation}
\label{a18}
1) The set $\lim({\cT})$ is not absolute, i.e., if $\bold V_1
\subseteq \bold V_2$ are two universes of set theory then in 
general $(\lim({\cT}))^{\bold V_1}$ will be a proper subset of 
$(\lim({\cT}))^{\bold V_2}$. 

\noindent
2) However, the notion of being a front is absolute: 
if $\bold V_1 \models ``A$ is a front in ${\cT}$", 
then there is a depth function $f:{\cT} \rightarrow \Ord$ satisfying 

\[
\eta \triangleleft \nu \text{ and } \forall k \le \ell g(\eta)
[\eta \rest k \notin A] \rightarrow f(\eta) > f(\nu).
\]

\mn
This function will also witness in $\bold V_2$ that $A$ is a front.

\noindent
3) $A \subseteq {\cT}$ contains a front \underline{if and only if} 
$A$ meets every branch of ${\cT}$.  So if $A \subseteq {\cT}$ contains 
a front of ${\cT}$ and ${\cT}' \subseteq {\cT}$ is a subtree, then 
$A \cap {\cT}'$ contains a front of ${\cT}'$.  Also this notion is absolute.
\end{observation}

\begin{notation}
\label{a21}
In several places in this section we will have an occasion to use the
following notation: Assume that $({\cT},\bar{\bold I})$ is a tagged 
tree, and for each $\eta \in {\cT}$ we are given a family ${\cP}_\eta$
of subsets of ${\cT}^{[\eta]}$ such that 

\[
\eta \triangleleft \nu \Rightarrow (\forall A \in {\cP}_\eta)(\exists
B \in {\cP}_\nu)[B \subseteq A].
\]

\mn
1) We inductively define for all 
$\alpha \in \Ord \cup\{\infty\}$ the property $\Dep_\alpha(\eta)$ by:
$\Dep_\alpha(\eta)$ \underline{if and only if} $(\forall \beta < \alpha)
(\forall A \in {\cP}_\eta)(\exists \nu \in A \cap \splitt({\cT}))
[\eta\triangleleft \nu$ and $\Dep_\beta(\eta)$ and 
$\{\rho:\rho \in \Succ_{\cT}(\nu)$ and $\Dep_\beta(\rho)\} \notin 
\bold I_\nu]$.

\noindent
2) Then it is easy to see that
$\Dep(\eta) := \max\{\alpha\in \Ord \cup\{\infty\}:\Dep_\alpha(\eta)\}$
is well defined, and $\Dep_\alpha(\eta) \Leftrightarrow \Dep(\eta) \ge
\alpha$.  We call $\Dep(\eta)$ the ``depth'' of $\eta$ (with respect to
the family $\bar{\cP} = \langle {\cP}_\eta:\eta \in {\cT}\rangle$ 
and the tagged tree $({\cT},\bar{\bold I})$).  
It is easy to check that $\eta \triangleleft \nu \Rightarrow 
\Dep(\eta) \ge \Dep(\nu)$.

\noindent
3) Similarly we can define $\Dp'_\alpha(\eta),\Dp'(\eta)$, when in the
definition of $\Dp_\alpha(\eta)$ we replace $\eta \triangleleft
\nu$ by $\eta = \nu$ in this case.
\end{notation}

\begin{definition}
\label{p2.5}
1) A tagged tree $({\cT},\bar{\bold I})$ is $\lambda$-complete \If \, for each
$\eta \in {\cT} \cap \Dom(\bar{\bold I})$ the ideal $\bold I_\eta$ 
is $\lambda$-complete.

\noindent
2) A family $\bbI$ of ideals is $\lambda$-complete \If \, each $\bold I \in
\bbI$ is $\lambda$-complete. We will only consider $\aleph_2$-complete
families $\bbI$.

\noindent
3) A family $\bbI$ is restriction-closed \If \, $\bold I \in \bbI,
A \subseteq \Dom (\bold I),A \notin \bold I$ implies $\bold I \rest A = 
\{B \in \bold I:B \subseteq A\}$ belongs to $\bbI$. 

\noindent
4) The restriction closure of $\bbI$ is 

\[
\res-\cl(\bbI) = \{\bold I \restriction A:\bold I \in \bbI, A
\subseteq \Dom (\bold I), A \notin \bold I\}.
\]

\mn
5) $\bold I$ is $\lambda$-indecomposable \If \, for every $A \subseteq \Dom
(\bold I),A \notin \bold I$, and $h:A \rightarrow \lambda$ there is 
$Y \subseteq \lambda,|Y| < \lambda$ such that $h^{-1}(Y) \notin \bold
I$. We say $\bar{\bold I}$ or ${\bbI}$, is $\lambda$-indecomposable if 
each $\bold I_\eta$ (or $\bold I\in {\bbI})$ is
$\lambda$-indecomposable; similarly in part (7). 

\noindent
6) $\bold I$ is strongly $\lambda$-indecomposable \If \,for $A_i \in
\bold I$ ($i < \lambda$) and $A \subseteq \Dom(\bold I),A \notin \bold
I$ we can find $B \subseteq A$ of cardinality $< \lambda$ 
such that for no $i < \lambda$ does $A_i$ include $B$. 
\end{definition}

\begin{observation}
\label{1.6n}
1) If an ideal $\bold I$ is $\lambda^+$-complete \then \, it is 
$\lambda$-indecomposable.

\noindent
2) If $\bold I$ is an ideal and $|\Dom(\bold I)|<\lambda$ \then \,
$\bold I$ is $\lambda$-indecomposable.

\noindent
3) If $\bold I$ is a strongly $\theta$-indecomposable ideal \then \, 
$\bold I$ is a $\theta$-decomposable ideal.
\end{observation}
  
\begin{lemma}
\label{1.7}
1) If $(\cT,\bar{\bold I})$ is a $\lambda^{+}$-complete tree and
  $\bold H$ is a function from $\lim({\cT})$ to $\lambda$ 
such that for every $ \alpha < \lambda$ the set $\bold
  H^{-1}(\{\alpha\})$ is a Borel subset of $\lim(\cT)$ 
(in the topology that was defined in Definition \ref{p2.4}(8)) \then \, there
is a tagged subtree $ ({\cT}^{\dagger},\bar{\bold I})$ satisfying 
$(\cT,\bar{\bold I}) 
\le^{*} ({\cT}^{\dagger},\bar{\bold I})$ (see Definition \ref{p2.3}(2)) 
such that $\bold H$ is constant on $\lim(\cT^{\dagger})$.

\noindent
2) In part (1) we can let $\bold H$ be multivalued, 
i.e. assume $\lim({\cT})$ is $\bigcup\limits_{\alpha<\lambda} 
\dot{\bbB}_\alpha$, each $\dot{\bbB}_\alpha$ is a Borel
subset of $\lim({\cT})$. If $({\cT},\bar{\bold I})$ is $\lambda^+$-complete
\then\ there is $({\cT}^\dagger,\bar{\bold I})$ such that $({\cT},
\bar{\bold I}) \le^* ({\cT}^\dagger,\bar{\bold I})$ and for 
some $\alpha<\lambda$ we have $\lim({\cT}^\dagger) \subseteq 
\dot{\bbB}_\alpha$.

\noindent
3) We can allow in (1) the function ${\bold H}$ to have values outside 
$\lambda$ as long as $|\Rang({\bold H})| \le \lambda$. Similarly (2).
\end{lemma}

\begin{PROOF}{\ref{1.7}}
1) First note that \If \, ${\cT}_{1} \subseteq {\cT}$ satisfies $(*)$ below
then $({\cT},\bar{\bold I}) \le^{*} ({\cT}_{1},\bar{\bold I} \rest {\cT}_{1})$
where:
\mn 
\begin{enumerate}
\item[$(*)$]  $\langle \rangle \in {\cT}_{1};\eta \triangleleft \nu \in 
{\cT}_1 \Rightarrow \eta \in {\cT}_1$; for every $\eta \in {\cT}_{1}$ 
if $\eta$ is a splitting point of $({\cT},\bar{\bold I})$ 
then $\Succ_{{\cT}_1}(\eta) = \Succ_{{\cT}}(\eta)$; and if $\eta$ 
is not a splitting point of $T$ then $|\Succ_{{\cT}_1}(\eta)|=1$.
\end{enumerate}
\mn
 So \wilog \, we can assume that in ${\cT}$ every point is either a
splitting  point or it has only one immediate extension i.e. 
$({\cT},\bar{\bold I})$ is standard. 

For each $\alpha<\lambda$ let us define a game $\Game_{\alpha}$: in the
first move the first player chooses the node $\eta_0$ in the tree such that
$\lg(\eta_0)=0$, the second player responds by choosing a proper subset
$A_0$ of $\Succ_{{\cT}}(\eta_0)$ such that $A_0 \in \bold I_{\eta_0}$. 
For $n>0$, in the $n$-th move, the first player chooses an immediate
extension $\eta_{n}$ of $\eta_{n-1}$,  such that $\eta_{n} \notin A_{n-1}$
or $\eta_{n-1}$ is not a splitting point of $({\cT},\bar{\bold I})$, and the
second player responds by choosing $A_{n} \in \bold I_{\eta_n}$.  

The first player wins if for the infinite branch $\eta$ defined by $\eta_0, 
\eta_1,\eta_2,\ldots$ we have $\bold H(\eta)=\alpha$. 
By the assumption of the lemma this is a Borel game so by 
Martin's Theorem, \cite{Mar75}  one of the
players has a winning strategy. We claim that for some $\alpha<\lambda$,
the first player has a winning strategy in the game $\Game_{\alpha}$. Assume 
otherwise, i.e., for every $ \alpha < \lambda $ the second player has a
winning strategy ${\bold f}_{\alpha} $. We construct an infinite branch
inductively: let $ \eta_{0} = \langle \rangle$ recalling 
$\eta_{0} \in {\cT}$. At stage $n$
let $A_n$ be $\bigcup\limits_{\alpha<\lambda} \bold f_\alpha(\eta_0,
\eta_{1},\ldots,\eta_{n-1})$; now if $\eta_{n-1}$ is a 
splitting point (of $({\cT},\bar{\bold I}))$ then 
$\bold I_{\eta_{n-1}}$ is $\lambda^{+} $-complete and each
$\bold f_\alpha(\eta_0,\ldots,\eta_{n-1})$ is a member of it, because
$\eta_0,F_\alpha(\eta_0),\eta_1,F_\alpha(\eta_0,\eta_1),\dotsc,\eta_{n-1}$
is an initial segment of a play of the game $\Game_\alpha$ in which
the second player uses the winning strategy $\bold f_\alpha$,
hence $A_{n} \in \bold I_{\eta_{n-1}}$, so 
clearly $\Succ_T(\eta_{n-1}) \nsubseteq A_n$. 

If $\eta_{n-1}$ is not a splitting point, it has only one immediate
successor and let it be $\eta_n$, otherwise since $\Succ(\eta_{n-1})
\notin \bar{\bold I}_{\eta_{n-1}},A_n \in \bar{\bold I}_{\eta_{n-1}}$, 
we have $(\Succ(\eta_{n-1}) \setminus A_n) \ne \emptyset$ so we 
can choose $\eta_n\in (\Succ_{{\cT}}(\eta_{n-1}) \setminus A_{n})$. 
Let $\eta = \bigcup\limits_{n<\omega} \eta_n$ be the infinite branch 
that we define by our construction and let
$\alpha(*) = \bold H(\eta)$. Now, in the game $\Game_{\alpha(*)}$, if the first
player chooses $\eta_n$ at stage $n$ (for all $n$) and the second player
plays by his strategy ${\bold f}_{\alpha(*)}$, the first player 
will win although the second player has used his winning strategy 
${\bold f}_{\alpha(*)}$, a contradiction.  

So there must be $\alpha(*)$ such that the first player has a winning
strategy ${\bold f}_{\alpha(*)}$ for $\Game_{\alpha(*)}$, and let 
${\cT}^{\dagger}$ be the subtree of ${\cT}$ defined by 
$\{\eta\in {\cT}:\eta=\langle\rangle$, or letting $n = \ell g(\eta)+1$ 
we have that $\langle \eta \rest 0,\ldots,\eta \restriction n\rangle$ 
are the first $n+1$ moves of the first player in a play in which he
plays according to ${\bold f}_{\alpha(*)}\}$. Now, for $\eta \in {\cT}^\dagger
\cap \splitt({\cT},\bar{\bold I})$, let $A = \Succ_{{\cT}^\dagger}(\eta)$. 
Then $A \notin \bold I_\eta$, otherwise the second player 
could have played it as $A_n$.  So $({\cT},\bar{\bold I}) \le^* 
({\cT}^\dagger,\bar{\bold I})$, and ${\cT}^\dagger$ is as required. 

\noindent
2) Same proof replacing $\bold H^{-1}(\{\alpha\})$ by 
$\dot{\bbB}_\alpha$, so $\bold H(\eta) = \alpha(*)$ by $\eta \in 
\dot{\bbB}_{\alpha(*)}$.              

\noindent
3) Trivial. 
\end{PROOF}

\begin{proof}
E.g.

\noindent
3) So let $A \subseteq \Dom(\bold I),A \in \bold I$ and $h:A
\rightarrow \lambda$ be given and we should find $Y \subseteq \lambda$
of cardinality $< \lambda$ such that $h^{-1}(Y) \notin \bold I$.  For
$i < \lambda$ let $A_i := h^{-1}\{i\}$, so as $\bold I$ is strongly
$\lambda$-indecomposable there $B \subseteq A$ of cardinality $<
\lambda$.  Let $Y = \{h(t):t \in B\}$ so clearly $Y$ is a subset of
$\lambda$ of cardinality $\le |B| < \lambda$, so it suffices to prove
that $h^{-1}(Y) \notin \bold I$.
\end{proof}

\begin{conclusion}
\label{1.9}
If $({\cT},\bar{\bold I})$ is a $\lambda^+$-complete tree, and 
$g$ is a function from ${\cT}$ into $\lambda$, and
${\lambda}^{\aleph_0} = \lambda$, \then \, there is a tagged subtree 
$({\cT}^\dagger,\bar{\bold I})$ satisfying $({\cT},\bar{\bold I})
\le^* ({\cT}^\dagger,\bar{\bold I})$ and such that $g\restriction 
{\cT}^\dagger$ depends only on the length of its argument, i.e.,  for some
function $g^\dagger:\omega \rightarrow \lambda$, for all $\eta \in 
{\cT}^\dagger$ we  have $g(\eta) = g^\dagger(\ell g(\eta))$. 
\end{conclusion}

\begin{PROOF}{\ref{1.9}}
Follows by \ref{1.7} for the function ${\bold H}, {\bold H}(\eta)
= \langle g(\eta\rest n):n < \omega\rangle$.
\end{PROOF}

\begin{lemma}
\label{1.10}
1) Assume that $\lambda$ is a regular uncountable cardinal,  and $({\cT},
\bar{\bold I})$ is a tagged tree such that for every $ \eta\in {\cT}$ 
$\bold I_{\eta} $ is $\lambda^+$-complete or $|\Succ_{\cT}(\eta)| <
\lambda$. If $\bold H:\lim(\cT) \rightarrow \lambda$ satisfies
``$\dot{\bbB}_\alpha := \{\eta\in \lim({\cT}):
\bold H(\eta)< \alpha\}$ is a Borel subset of $\lim({\cT})$ for any successor
$\alpha < \lambda$", \then \, there are $\alpha < \lambda$ and
$(\cT',\bar{\bold I})$ satisfying $({\cT},\bar{\bold I}) \le^{*} 
({\cT}',\bar{\bold I})$ and such that for all $\eta \in {\cT}'$ 
we have $\bold H(\eta) < \alpha$, and for all $\eta$ in ${\cT}'$, 
if $|\Succ_{\cT}(\eta)| < \lambda$, then $\Succ_{{\cT}'}(\eta) =
\Succ_{\cT}(\eta)$.

\noindent
2) Like part (1) but we omit the function $\bold H$ and just assume
$\dot{\bbB}_\alpha$ is a Borel subset of $\lim(\cT)$ for 
$\alpha < \lambda$ but demand $\bigcup\limits_{\alpha < \lambda}
\dot{\bbB}_\alpha = \lim(\cT)$; moreover every
   $X \subseteq \lim(\cT)$ of cardinality $< \lambda$ is included in
   some $\dot{\bbB}_\alpha,\alpha < \lambda$.

\noindent
3) Let $\lambda,\mu$ be uncountable cardinals satisfying
$\lambda^{<\mu} = \lambda$ and let 
$(\cT,\bar{\bold I})$ be a tree in which for each
$\eta \in \cT$ either $|\Succ_{\cT}(\eta)| < \mu$ or $\bar{\bold I}(\eta)$ is
$\lambda^+$-complete.  For $A \subseteq \cT$ and $\eta \in \cT$ we
define $\rest_{\cT}(\eta,A)$ as the sequence $\langle x_\ell:\ell <
\ell g(\eta)\rangle$ when $x_\ell$ is $\eta(\ell)$ if $\eta \rest \ell
\in A$ and zero if $\eta \rest \ell \in A$.
\Then \, for every function $\bold H:\cT \rightarrow
   \lambda$ there exists $\cT',(\cT,\bar{\bold I}) \le^* 
(\cT',\bar{\bold I})$ such that (letting $A = \{\eta \in
\cT:|\Suc_{\cT}(\eta)| 
< \mu\}$ hence $\rest_{\cT}(\eta,A) \in {}^{\omega >}\mu$ for $\eta \in \cT$):
\mn
\begin{enumerate}
\item[$\bullet$]  for $\eta,\eta' \in \cT'':\rest_{\cT}(\eta,A) = 
\rest_{\cT}(\eta',A)$ implies: $\bold H(\eta) =
\bold H(\eta')$ and $\eta \in A$ iff $\eta' \in A$, and if $\eta \in \cT'
   \cap A$, then $\Suc_{cT}(\eta) = \Suc_{\cT'}(\eta)$.  
\end{enumerate}
\end{lemma}

\begin{PROOF}{\ref{1.10}}
1) We define for each successor $\alpha<\lambda$ a game $\Game_{\alpha}$
very much like the way we did it for proving Lemma \ref{1.7}, the only
difference being that if $ |\Succ_{\cT} (\eta_n)| < \lambda$, the
second player chooses $A_n$ such that $|\Succ_{\cT}(\eta_n)
\setminus A_n|=1$, otherwise the second player chooses $A_n \in 
\bold I_{\eta_n}$ just like in \ref{1.7}.  
The first player wins if $\bold H(\eta_n) <\alpha$ for every
$n<\omega$.  Here again the game $\Game_{\alpha} $ is determined for every
$\alpha$ (here simply because if the second player wins a play he does so at
some finite stage). Again we claim that there should be at least one
successor $\alpha<\lambda$ for which the first player has a winning
strategy. Assume the contrary, and for each $\alpha<\lambda$ let
${\bold f}_{\alpha}$ be a winning strategy of the second player in the game
$\Game_{\alpha+1}$. We construct a subtree ${\cT}^{*} $ deciding by
induction on the length of the members of ${\cT}$ which of them are 
members of ${\cT}^*$. For $\eta $ that is already in ${\cT}^{*} $, if
$|\Succ_{\cT}(\eta)|<\lambda$ we include all the members of 
$\Succ_{\cT}(\eta)$ in ${\cT}^{*} $; otherwise $\bold I_{\eta}$ 
is $\lambda^{+}$-complete so $\Succ_{\cT}(\eta) \setminus
\bigcup\limits_{\alpha<\lambda} {\bold f}_{\alpha}(\eta \restriction
0,\eta \restriction 1,\ldots,\eta)$ is not
empty; pedantically you use $\Succ_{\cT} (\eta) \setminus \cup 
\{\bold f_\alpha(\eta\rest 0,\eta\rest 1,\ldots,\eta):{\bold f}_\alpha 
(\eta \rest 0,\eta \rest 1,\ldots,\eta)$ is well defined$\}$, 
so we pick one extension of $\eta$ from this set and the rest of 
$\Succ_{\cT}(\eta)$ will not be in ${\cT}^{*} $. Now ${\cT}^{*}$ is 
a tree of height $\omega$ such that each member has less than
$\lambda$ immediate successors. So, as $\lambda$ is regular uncountable, we
get $|{\cT}^{*}| < \lambda$ and hence there is some successor ordinal
$\alpha^{*}<\lambda$
such that $\eta \in {\cT}^{*}$ implies 
$\bold H(\eta)<\alpha^{*}$. Regarding the
game $\Game_{\alpha^*} $, there is a play of it in which the first player
chooses all along the way members of ${\cT}^{*} $ and the second player
plays according to ${\bold f}_{\alpha^*} $; of 
course the first player wins this
game contradicting the assumption that ${\bold f}_{\alpha^*}$ is a winning
strategy for the second player.

Hence, for some successor $\alpha^*$, the second player has a winning
strategy in the game $\Game_{\alpha^*}$. We define ${\cT}'$ just
like we did in the proof of Lemma \ref{1.7},  collecting all the initial
segments of plays of the first player in the game $\Game_{\alpha^{*}}$ when he
plays according to his winning strategy $\bold H_{\alpha^{*}}$. 

\noindent
2) Same proof, (pedantically, \wilog \, $\bbB_\alpha = \emptyset$ for
   $\alpha$ limit).

\noindent
3) Similarly.
\end{PROOF}

\noindent
The following (really part (2)) will be used in the proof of \ref{a45}.
\begin{lemma}
\label{a39}
1) Assume
\mn
\begin{enumerate}
\item[$(a)$]  $(\cT,\bar{\bold I})$ is an ${\bbI}$--tree, ${\bbI}$ a family of
  ideals
\sn
\item[$(b)$]   $\lim({\cT}) = \bigcup\limits_{i<\theta} \, 
\bigcup\limits_{\epsilon < \theta_i} \dot{\bbB}_{i,\epsilon}$, 
each $\dot{\bbB}_{i,\epsilon}$ is a Borel set, increasing with
$\epsilon$
\sn
\item[$(c)$]  $(\alpha) \quad {\bbI}$ is $\partial$-complete, and
\sn
\item[${{}}$]  $(\beta) \quad$ each $\bold I \in {\bbI}$ is strongly 
$\theta$-indecomposable
\sn
\item[$(d)$]  $E_i$ is a $\partial$-complete filter on $\theta_i$
\sn
\item[$(e)$]  if $i < \theta,A_\varepsilon \in \bold I_\eta$ for
  $\varepsilon < \theta_i$ then for some $A \in I_\varepsilon$ we have
  $\sup\{\varepsilon < \theta_i:A_\varepsilon \subseteq A\} \in E_i$
\sn
\item[$(f)$]   $\partial = \cf(\partial)$ and $\partial + \aleph_1 \le
  \theta = \cf(\theta)$
\sn
\item[$(g)$]  $(\forall \alpha < \theta)(|\alpha|^{\aleph_0} <
  \theta)$
\sn
\item[$(h)$]  $(\forall \alpha<\partial)
(|\alpha|^{\aleph_0}<\partial)$ or each 
$\dot{\bbB}_{\zeta,\epsilon}$ is closed 
\sn
\item[$(i)$]   $\dot{\bbB}_i := \bigcup\limits_{\epsilon < \theta_i} 
\dot{\bbB}_{i,\epsilon}$ is increasing with $i$ 
\sn
\item[$(j)$]  $\eta\in {\cT} \setminus \splitt({\cT},\bar{\bold I}) \Rightarrow
|\Succ_{\cT}(\eta)|<\partial$.
\end{enumerate}
\mn
\Then \, for some $i<\theta$ and $\epsilon < \theta_i$ and ${\cT}'$ 
we have $({\cT},\bar{\bold I}) \le^\otimes ({\cT}',\bar{\bold I})$, 
and $\lim({\cT}') \subseteq \dot{\bbB}_{i,\epsilon}$; see Definition \ref{p2.3}(3).

\noindent
2) Assume $({\cT},\bold I)$ be an ${\bbI}$--tree, ${\bbI}$ a family 
of ideals, $\lim({\cT}) = \bigcup\limits_{i<\theta}
\bigcup\limits_{\varepsilon < \varepsilon_i} \dot{\bbB}_{i,\varepsilon}$, 
each $\dot{\bbB}_{i,\varepsilon}$ is a Borel set, 
$i < \theta \Rightarrow \varepsilon_i < \theta]$, 
${\bbI}$ is $\theta$-complete, $\theta$ is regular uncountable 
and each $\bold I \in {\bbI}$ is strongly 
$\theta$-indecomposable, and $\dot{\bbB}_i := 
\bigcup\limits_{\varepsilon < \varepsilon_i}
\dot{\bbB}_{i,\varepsilon}$ is increasing with $i$ and
 
\[
\eta \in {\cT} \setminus \splitt({\cT},\bar{\bold I}) \Rightarrow 
|\Succ_{\cT}(\eta)| < \theta.
\]

\mn
\Then \, for some $i<\theta$ and $\varepsilon<\varepsilon_i$ and 
$\cT'$ we have $(\cT,\bar{\bold I}) \le^\otimes (\cT',\bar{\bold I}$), and 
$\lim(\cT) \subseteq \dot{\bbB}_{i,\varepsilon}$.
\end{lemma}

\begin{PROOF}{\ref{a39}}
1) We first prove part (2).
\smallskip

\noindent
\underline{Proof of part (2)}:  We define, for 
$i < \theta$ and $\epsilon < \varepsilon_i$ a game
$\Game_{i,\epsilon}$ as in the proof of \ref{1.7}, \ref{1.9} 
for the set $\dot{\bbB}_{i,\epsilon}$.
If for some $i <\theta,\epsilon < \varepsilon_i$ the first player wins, then
we get the desired conclusion as in the earlier proofs. 
Otherwise, as each such game is determined
(as $\bbB_{i,\epsilon}$ is a Borel set) there is a winning strategy
${\bold f}_{i,\epsilon}$ for the second player in the game 
$\Game_{i,\epsilon}$. 
Let $\eta \in \splitt({\cT},\bar{\bold I})$. For each $i <\theta$ 
we define a set $A^i_\eta \subseteq \Succ_{\cT}(\eta)$ by 
$A^i_\eta = \cup \{A \subseteq \Succ_{\cT} (\eta)$: 
for some $\varepsilon < \varepsilon_i$ in some play of
the game $\Game_{i,\varepsilon}$ in the $n$-th move 
the first player chooses $\eta$ and the second player 
chooses $A$ by the strategy ${\bold f}_{i,\epsilon}\}$.  Recalling $i <
\theta \Rightarrow \varepsilon_i < \theta$, as
$\bold I_\eta$ is $\theta$-complete clearly $A^i_\eta \in 
\bold I_\eta$.  As $\bold I_\eta$ is strongly 
$\theta$-indecomposable applying the definitions to 
$\langle A^i_\eta:i < \theta\rangle$ we 
can find $B_\eta \subseteq \Suc_{\cT}
(\eta)$ of cardinality $<\theta$ such that $i <\theta
\Rightarrow B_\eta \nsubseteq A^i_\eta$. 
(If we add $\Dom(\bold I_\eta) = \Succ_{\cT}(\eta)$ 
we can in Definition \ref{p2.5}(5) use $A = \Dom(\bold I)$).
Now as in the proof of \ref{1.7} we choose ${\cT}'_n\subseteq 
\{\eta\in {\cT}:\ell g(\eta)=n\}$ by induction on $n$ as follows:
${\cT}'_0=\{\langle\rangle\}, {\cT}'_{n+1}=\cup \{\nu$: for some 
$\eta\in {\cT}_n,\nu\in \Succ_{\cT}(\eta)$ and 
$[\eta\in \splitt({\cT},\bar{\bold I}) \Rightarrow \nu \in 
B_\eta]\}$. 

Let ${\cT}' = \cup \{{\cT}'_n:n<\omega\}$, clearly 
${\cT}' \subseteq {\cT}$ is non-empty, closed under initial 
segments. As $\theta$ is regular and $\eta\in {\cT}'
\setminus \splitt({\cT},\bar{\bold I}) \Rightarrow |\Succ_{\cT} 
(\eta)| < \theta$ and $\eta\in \splitt({\cT},\bar{\bold I}) \Rightarrow 
|B_\eta|<\theta$ clearly $n < \omega \Rightarrow |{\cT}'_n| <\theta$ 
and as $\theta$ is uncountable also $|{\cT}'|<\theta$ hence 
$\lim({\cT}')$ has cardinality $< \theta$. 
As $\langle \dot{\bbB}_i:i < \theta\rangle$ is $\subseteq$-increasing 
with union 
$\lim({\cT})$, clearly for some $i(*)<\theta$ we have 
$\lim({\cT}')\subseteq \dot{\bbB}_{i(*)}$. 

Clearly there is $\eta \in \lim(\cT')$, hence for some $\varepsilon <
\varepsilon_i$ we have $\eta \in \dot{\bbB}_{i(*),\varepsilon}$, but
there is a play of the game $\Game_{i,\varepsilon}$ in which the moves
of the first player are $\langle \eta \rest n:n < \omega\rangle$.
Easy contradiction.
\smallskip

\noindent
\underline{Proof of part (1)}:  We begin as in the 
proof of part (2) until. ``For each $i < \theta$
we define a set $A^i_\eta \ldots$".  Now for each $i < \theta$ and
   $\varepsilon < \theta_i$ we define a set $A^{i,\varepsilon}_\eta
\subseteq \Succ_{\cT}(\eta)$ by: if there is a play of the game
 $\Game_{i,\varepsilon}$ in which the second player uses the strategy
$\bold f_{i,\varepsilon}$ and the first player chooses $\eta$ in
   the $n$-th move, then the second player chooses
   $A^{i,\varepsilon}_\eta$ (note there is at most one such play); if
   there is no such play then let $A^{i,\varepsilon}_\eta =
   \emptyset$.  As $\bbI$ satisfies clause $(e)$ of the assumption
   there is a set $A^i_\eta \subseteq \Succ_{\cT}(\eta)$ satisfying
   $A^i_\eta \in \bold I_\eta$ such that $\{\varepsilon <
   \theta_i:A^{i,\varepsilon}_\eta \subseteq A^i_\eta\} \in E_i$.

Now we continue as in the rest of the proof of part (2) after the
choice of $A^i_\eta$.  In particular, we choose $B_\eta$ (for every
$\eta \in \cT$) and $\cT'_n$ for $n < \omega$ and $\cT'$ and $i(*)$
such that $\lim(\cT') \subseteq \dot{\bbB}_{i(*)}$.

Now for every $\eta\in {\cT}'_n \cap \splitt({\cT},\bar{\bold I})$ 
we know that $B_\eta = \Succ_{{\cT}'} (\eta) \subseteq \Succ_{\cT}
(\eta)$ so there is $\rho_\eta\in B_\eta \setminus A_\eta^{i(*)}$. 
We now choose ${\cT}''_n \subseteq {\cT}'_n$ 
by induction on $n$ as follows: ${\cT}''_n = h \langle\rangle\}$,
${\cT}''_{n+1} = \{\nu$: for some $\eta\in {\cT}''_n, \nu\in 
\Suc_{{\cT}'} (\eta)={\cT}'_{n+1} \cap \Suc_{\cT}
(\eta)$, and $[\eta \in \splitt({\cT},\bar{\bold I}) \Rightarrow 
\nu = \rho_\eta]\}$. 
So ${\cT}'' = \cup\{{\cT}''_n:n<\omega\}$ is a non-empty subset of 
${\cT}'$, closed under initial segments and 
$|{\cT}''_n| < \partial$ and $\lim({\cT}')\subseteq 
\dot{\bbB}_{i(*)}= \bigcup\{\dot{\bbB}_{i(*),\epsilon}:
\epsilon < \varepsilon_{i(*)}\},\dot{\bbB}_{i(*),\epsilon}$ 
increasing with $\epsilon$.  
As $(\forall \alpha<\partial) (|\alpha|^{\aleph_0}<\partial)$ or each 
$\dot{\bbB}_{i(*),\epsilon}$ is closed for some $\epsilon < \theta_{i(*)}$ 
we have $\lim({\cT}'')\subseteq \dot{\bbB}_{i(*),\epsilon}$.  As $E_i$
is $\partial$-complete increasing $\varepsilon$ we have: $\eta \in
\cT' \Rightarrow A^{i,\varepsilon}_\eta \in A^i_\eta$. 
But easily we can find a play of the game $\Game_{i(*),\epsilon}$ 
in which the second player uses the strategy ${\bold f}_{i(*),\epsilon}$ 
and the first player choose $\eta_n$ from ${\cT}''$. In such 
a play the first player wins, contradicting the choice of 
${\bold f}_{i(*),\epsilon}$ .
\end{PROOF}

\noindent
The following uses pcf in its phrasing (hence in its proof)
\begin{lemma}
\label{a42} 
Suppose $({\cT},\bar{\bold I})$ is an ${\bbI}$-tree, $\theta$ regular
uncountable, $\langle A_\eta:\eta\in {\cT}\rangle$ is such that:
$A_\eta$ is a set of ordinals, $[\eta\triangleleft \nu \Rightarrow A_\eta
\subseteq A_\nu ]$ and
\mn
\begin{enumerate}
\item[$(*)$]   $(a) \quad {\bold S}$ is a set of uncountable regular cardinals
\sn
\item[${{}}$]  $(b) \quad \bbI' := \bbI \setminus \{\bold I \in \bbI:
|\Dom(\bold I)| < \mu\}$ is $\mu^+$-complete or at least strongly

\hskip25pt $\mu$-indecomposable for every $\mu$ such that $\mu \in \bold
S$ or $\mu \in \pcf({\bold S} \cap A_\eta)$ 

\hskip25pt for some $\eta\in {\cT}$
\sn
\item[${{}}$]  $(c) \quad \bbI$ is $\theta$-complete and $|\pcf({\bold S}\cap
A_\eta)| < \theta$ for $\eta\in {\cT}$ and $\theta \le \min({\bold
  S})$,
\sn
\item[${{}}$]  $(d) \quad |A_\eta| < \min({\bold S})$ for $\eta \in {\cT}$
\end{enumerate}
\mn
\Then \, there is ${\cT}^\dagger$  satisfying $({\cT},\bar{\bold I}) \le^*
({\cT}^\dagger,\bar{\bold I})$ and such that:
\mn
\begin{enumerate}
\item[$(**)$]   if $\lambda \in A_\nu \cap {\bold S}$ and $\nu \in 
{\cT}^\dagger$ \then \, for some $\alpha_\nu(\lambda) < \lambda$ 
for every $\rho$ such that $\nu \triangleleft \rho \in
\lim(\cT^\dagger)$ we have
$\alpha_\nu(\lambda) \ge \sup(\lambda \cap \bigcup\limits_{n <\omega} 
A_{\rho \rest n})$.
\end{enumerate}
\end{lemma}

\begin{PROOF}{\ref{a42}}
It is enough to prove the existence of a ${\cT}^\dagger$ as required just
for $\nu = \langle \rangle$, (as we can repeat the 
proof going up in the tree). This will be proved by induction on 
$\max(\pcf({\bold S} \cap A_{\langle \rangle})$)
(exists,  see \cite[Ch.I,1.9]{Sh:g}).  Let $\alpha_\lambda(\eta) 
= \sup(A_\eta \cap \lambda)$. 

We assume knowledge of \cite{Sh:g} and use its notation. 

Let $\ga := {\bold S} \cap A_{\langle\rangle}$ (if $\ga$ is empty we have
nothing to do), let $\mu = \max \pcf(\ga)$, and let $\langle f_\zeta:\zeta<\mu
\rangle$ be $<_{\bold J_{<\mu}[\ga]}$-increasing and cofinal in 
$\Pi \ga$, recalling that the later means that $(\forall f \in \Pi \ga)(\exists
\zeta<\mu)(f<_{\bold J_{<\mu}[\ga]} f_\zeta)$. 
Let $\{\gb_\varepsilon:\varepsilon < \varepsilon(*)\}$ be cofinal 
in $J_{<\mu}[\ga]$, e.g., this set is 
$\{\bigcup\limits_{\theta \in \gc}{\gb}_\theta[\ga]:\gc \subseteq \pcf
(\ga) \setminus \{\mu\}$ is finite$\}$, so by clause (c) of the assumption 
$(*)$ we can have $\epsilon(*)<\theta$ and hence by assumption (c)
${\bbI}'$ is $|\varepsilon(*)|^+$-complete.

For $\varepsilon < \varepsilon(*)$ and $\zeta<\mu$ we consider the statement:
\mn
\begin{enumerate}
\item[$(*)^\varepsilon_\zeta$]  there is a subtree $\cT'$ of $\cT$
  satisfying $({\cT},\bar{\bold I}) \le^* ({\cT}',\bar{\bold I})$ 
such that for every $\eta\in \lim({\cT}')$ and $\lambda \in \ga \setminus 
\gb_\varepsilon$ and $n$ we have $\alpha_\lambda(\eta \rest n) \le 
f_\zeta(\lambda)$.
\end{enumerate}
\mn
It suffices to find such ${\cT}'$ (for some $\varepsilon,\zeta)$ because:
we can apply the induction
hypothesis on $({\gb}_\epsilon,{\cT}')$, this is justified as 
$\max \pcf({\gb}_\epsilon) < \max \pcf({\ga})$.

In ${\bold V}$ define for $\zeta<\mu$ and $\varepsilon <
\varepsilon(*)$ the following set:

\[
\dot{\bbB}_{\zeta,\varepsilon} := \{\eta\in \lim({\cT}): \text{ for
every } \lambda \in \ga \setminus \gb_\varepsilon, n<\omega\Rightarrow
(\lambda \cap A_{\eta \rest n})\subseteq f_\zeta(\lambda)\}.
\]

\mn
Clearly $\dot{\bbB}_{\zeta,\varepsilon}$ is closed and
$\dot{\bbB}_\zeta = \bigcup\limits_{\varepsilon < \vare(*)}
\dot{\bbB}_{\zeta,\vare}$.  Now, $\zeta<\xi<\mu \Rightarrow \dot{\bbB}_\zeta
\subseteq \dot{\bbB}_\xi$ (as $f_\zeta<_{\bold J_{<\mu}[\ga]}f_\xi$) 
and $\lim({\cT})= \bigcup\limits_{\zeta<\mu} \dot{\bbB}_\zeta$ 
(as $\langle f_\zeta:\zeta < \mu \rangle$ is cofinal in
$(\prod,<_{J_{< \mu}[\ga]}$), hence using \ref{a39}(2) above 
(with $\mu$, $\vare(*)$ here 
standing for $\theta,\vare_i$ there) for some $\zeta(*)<\mu$ and $\vare
< \vare(*)$ and ${\cT}'$ we have $({\cT},\bar{\bold I}) \le^* ({\cT}',
\bold I)$ and $\lim({\cT}')\subseteq \dot{\bbB}_{\zeta,\vare}$. 
So $(*)^\vare_\zeta$ holds, but as said above this suffices.
\end{PROOF}

\noindent
The following is used in \cite[1.11,1.13]{Sh:511}
\begin{lemma} 
\label{a45}
Let $\theta$ be an uncountable  regular cardinal (the main case here is
$\theta = \aleph_1)$.  Let ${\bbI}$ be a family of $\theta^+$-complete
ideals, $({\cT}_0,\bar{\bold I})$ a tagged tree, $A=\{\eta\in \cT_0: 
0 < |\Succ_{{\cT}_0} (\eta)| \le \theta\},[\eta\in {\cT}_0
\setminus A \Rightarrow \bold I_\eta \in {\bbI}$ and $\Succ_{{\cT}_0} 
(\eta) \notin \bold I_\eta]$, and $[\eta\in A \Rightarrow 
\Succ_{{\cT}_0}(\eta) \subseteq \{\eta \char 94 \langle i\rangle: 
i < \theta\}]$, and $\bold H:{\cT}_0 \rightarrow \theta$ and 
$\bar{\bold c} = \langle \bar{\bold c}_\eta: \eta \in
A\rangle$, is such that for all $\eta \in A,\bold c_\eta$ is a 
club of $\theta$. \Then \, 
there is a club $C$ of $\theta$ such that: for each 
$\delta \in C$ there is ${\cT}_\delta \subseteq {\cT}_0$ satisfying:
\mn
\begin{enumerate}
\item[$(a)$]  ${\cT}_\delta$ a tree
\sn
\item[$(b)$]  if $\eta\in {\cT}_\delta$ and $|\Succ_{{\cT}_0}(\eta)|
<\theta$, then $\delta \in \bold c_\eta$ and
$\Succ_{{\cT}_\delta}(\eta) = \Succ_{{\cT}_0}(\eta)$,
and if in addition $|\Succ(\eta)| = \theta$, then $\Succ_{{\cT}_\delta}(\eta)
=\{\eta \char 94 \langle i \rangle:i < \delta\}\cap
\Succ_{{\cT}_0}(\eta)$ 
\sn
\item[$(c)$]  $\eta \in {\cT}_\delta \setminus A$ implies 
$\Succ_{{\cT}_\delta}(\eta) \notin \bold I_\eta$
\sn
\item[$(d)$]  for every $\eta \in {\cT}_\delta$ we have 
$\bold H(\eta)<\delta$.
\end{enumerate}
\end{lemma}

\begin{PROOF}{\ref{a45}}
For each $\zeta < \theta$ we define a game $\Game_\zeta$.  The game lasts
$\omega$ moves, in the $n$th move $\eta_n \in {\cT}_0$ of 
length $n$ is chosen. 

\noindent
For $n = 0$:  necessarily $\eta_0 = \langle \rangle$.

\noindent
For $n=m+1$:  If $|\Succ_{{\cT}_0}(\eta_m)| = \theta$,
then the {\it second} player chooses $\eta_{m+1} \in 
\Succ_{{\cT}_0}(\eta_m)$ satisfying $\eta_{m+1}(m)<\zeta$.

\noindent
If $|\Succ_{{\cT}_0}(\eta_m)| < \theta$, then the
{\it second} player chooses any $\eta_{m+1} \in
\Succ_{{\cT}_0}(\eta_m)$. 

\noindent
If $\eta_m \notin A$, then the second player chooses $A_m\in
\bold I_{\eta_m}$, and then the {\it first} player chooses
$\eta_{m+1} \in \Succ_{{\cT}_0}(\eta_m)\setminus A_m$.

At the end, the first player wins if for all $n,\bold H(\eta_n)<\zeta$ and
$|\Succ_{{\cT}_0}(\eta_n)| = \theta \Rightarrow \zeta\in \bold
c_{\eta_n}$.  

Now clearly
\mn
\begin{enumerate}
\item[$(*)$]   if for a club of $\zeta < \theta$ the first
player has a winning strategy for the game $\Game_\zeta$, {\it then} 
there are trees ${\cT}_\delta$ as required.
\end{enumerate}
\mn
Let $S = \{\delta < \theta$: first player does not have a winning
strategy for the game $\Game_\delta\} $;  we assume that the set $S$ is
stationary, and get a contradiction, this suffice. 

For $\delta \in S$ let ${\bold f}_\delta$ be a winning strategy 
for the second player in $\Game_\delta$ (he has a winning strategy 
as the game is determined being closed for the first player). 
So ${\bold f}_\delta$ gives for the first
$(n-1)$-moves of the first player, the $n$-th move of the second player.

Let $\chi$ be a large enough regular cardinal, and let $N_0 \prec
({\cH}(\chi),\in)$ be such that $\theta+1\subseteq N_0$, $\|N_0\|
= \theta$, $({\cT}_0,\bar{\bold I})\in N_0,\bar{\bold c} \in N_0$, and 
$\bar{\bold f} = \langle {\bold f}_\delta: \delta \in S \rangle \in N_0$. 
We can find $N_1 \prec N_0$ such that 
$\|N_1\| < \theta,N_1 \cap \theta$ is an ordinal and
$({\cT}_0,\bar{\bold I}) \in N_1,\langle {\bold f}_\delta: 
\delta \in S \rangle\in N_1$ and $\bar{\bold c} \in N_1$.  
Let $\delta := N_1 \cap \theta$. Since $S$ was
assumed to be stationary, we may assume that $\delta\in S$. 

Now we shall choose by induction on $n, \eta_n \in T_0 \cap N_1$ of length
$n$, such that $\langle \eta_\ell: \ell \le n \rangle$ is an initial segment
of a play of the game $\Game_\delta$ in which the second player uses his
winning strategy ${\bold f}_\delta$. (The $A_\ell \in 
\bold I_{\eta_\ell}$ are not
mentioned as they are not arguments of ${\bold f}_\delta)$.
\medskip

\noindent
\underline{Case 1}.   $n=0$:

We let $\eta_0=\langle\rangle$.
\medskip

\noindent
\underline{Case 2}.   $n = m+1$, $\eta_m \in A$:

Recall that as $\delta\in S$, the second player has the 
winning strategy ${\bold f}_\delta$ for the
game $\Game_\delta$ \underline{but} in general ${\bold f}_\delta\notin N_1$.  
So ${\bold f}_\delta$ gives us $\eta_n$. Now if $|\Succ_{{\cT}_0}
(\eta_m)| < \theta$ then $\Succ_{{\cT}_0}(\eta_m) \subseteq 
N_1$ (because ${\cT}_0$, $\eta_m$ belong to $N_1$ 
and $N_1\cap \theta$ is an ordinal), and hence $\eta_n\in N_1$ as
required. If $|\Succ_{{\cT}_0}(\eta_m)|  = \theta$ then
necessarily $\Succ_{{\cT}_0}(\eta_m) \subseteq \{\eta_m \char 94 \langle
i \rangle:i<\theta\},\eta_n = \eta_m \char 94 \langle i \rangle, i<\delta$ (as
the play is of the game $\Game_\delta$), but $N_1 \cap \theta =
\delta$ so necessarily $i\in N_1$ hence (as
$\eta_m \in N_1)$ also $\eta_n \in N_1$.   

Lastly,
\medskip

\noindent
\underline{Case 3}.  $n = m+1$, $\eta_m \notin A$:

So ${\bold f}_\delta$ gives us $A_m^\delta \in \bold I_{\eta_m}$ which is not
necessarily in $N_1$,  however we let $A^* = \bigcup\{A^\zeta_m: \zeta\in
S$, and there is a play of $\Game_\zeta$ in which $\langle \eta_\ell: \ell \le
m \rangle$ were played (by the first player) 
and the second player plays according to ${\bold f}_\zeta$
(this play is unique) and the strategy ${\bold f}_\zeta$ dictates to the second
player to choose $A^\zeta_m\}$. 

\noindent
Now, $A^*$ belongs to $N_1$ (as $\bar {\bold f} \in N_1$) and 
being the union of $\le \theta$ members of $\bold I_{\eta_m}$ it 
belongs to $\bold I_{\eta_m}$, and hence
$A^*\cap \Succ_{{\cT}_0}(\eta_m)$ is a proper subset of 
$\Succ_{{\cT}_0}(\eta_m)$.  Consequently, there is $\eta_m \char 94
\langle  i \rangle \in \Succ_{{\cT}_0} (\eta_m) \setminus A^*$, and 
thus there is such $i\in N_1$.
Let the first player choose $\eta_n = \eta_m \char 94 \langle i \rangle$. 

So we have played a sequence $\langle\eta_n:n < \omega \rangle$ of elements of
$N_1$, always obeying ${\bold f}_\delta$ so this sequence was 
produced by a play of $\Game_\delta$ in which the second player 
plays according to the strategy ${\bold f}_\delta$. 
But then, for all $n,\eta_n\in N_1\Rightarrow \bold H(\eta_n) \in
N_1$, so $\bold H(\eta_n)<\delta$, and

\[
\eta_n \in N_1\ \Rightarrow \bold c_{\eta_n}\in N_1 \Rightarrow \delta = \sup
(\bold c_{\eta_n} \cap \delta) \Rightarrow \delta \in \bold c_{\eta_n};
\]

\mn
hence the first player wins in this play. So ${\bold f}_\delta$ 
cannot be a winning strategy for the second player in 
$\Game_\delta$. A contradiction, so $S$ is not stationary and we are done.
\end{PROOF}

\begin{claim}
\label{a48} 
Assume $\kappa < \lambda$ and $\cf([\lambda]^{< \kappa^+},\subseteq) =
\lambda$ and $\lambda = \lambda^{\aleph_0}$.

\noindent
1) If $\chi > \lambda^+$ and $x\in {\cH}(\chi)$ \then \, we can find 
$\bar N = \langle N_\eta: \eta\in \cT\rangle$ such that: 
\mn
\begin{enumerate}
\item[$(a)$]  $\cT$ is a subtree of ${}^{\omega>}(\lambda^+)$, each $\eta\in
\cT$ is (strictly) increasing, 
\sn
\item[$(b)$]  $N_\eta \prec ({\cH}(\chi), \in, <^*_\chi)$,
\sn
\item[$(c)$]  $x \in N_\eta$ and $\kappa+1 \subseteq N_\eta$ and $\|N_\eta\| =
\kappa$, 
\sn
\item[$(d)$]  $\nu \triangleleft \eta\in {\cT} \Rightarrow N_\nu\prec N_\eta$,
\sn
\item[$(e)$]  $N_\eta\cap N_\nu = N_{\nu\cap \eta}$ for $\eta$, $\nu\in\cT$, 
\sn
\item[$(f)$]  $\eta\in N_\eta$,
\sn
\item[$(g)$]  if $\eta_\ell \char 94 \langle \alpha_\ell \rangle \in
\cT$ for $\ell=1,2$ and $\alpha_1<\alpha_2$ then $\sup(N_{\eta_1 \char 94
\langle\alpha_1 \rangle} \cap \lambda^+) < 
\min(N_{\eta_2 \char 94 \langle\alpha_2\rangle}\cap \lambda^+
\backslash \alpha_1)$.
\end{enumerate}
\mn
Recall that $\cf([\kappa^{+n}]^{\le \kappa},{\subseteq})=\kappa^{+n}$.

\noindent
2) If in addition $\lambda= \lambda^\kappa$ (equivalently, $2^\kappa
\le \lambda$) \then \, we can add:
\mn
\begin{enumerate}
\item[$(h)$]  if $\eta,\nu \in \cT$ have the same length {\it then} there is
an isomorphism from $N_\eta$ onto $N_\nu$, call it $f_{\eta,\nu}$, which
maps $x$ to itself, so 
\[
\eta,\nu \in \lim(\cT) \Rightarrow \bigcup\limits_{n < \omega} N_{\eta
\rest n} := N_\eta \cong N_\nu := \bigcup\limits_{n < \omega} N_{\nu \rest n}.
\]
\end{enumerate}
\mn
3)  If $\cS \subseteq [\lambda]^{\le \kappa}$ is stationary of
cardinality $\lambda$ \then \, we can in (1) demand 
\mn
\begin{enumerate}
\item[$(i)$]  $N_\eta \cap \lambda\in {\cS}$.
\end{enumerate}
\mn
4)  We can further demand (in parts (1),(2)) that: 
\mn
\begin{enumerate}
\item[$(j)$]  $N_\eta$ is the Skolem hull of $\{x,\eta,\kappa,\lambda\}\cup
\kappa \cup N_{\langle\rangle}$ in $({\cH}(\chi),\in,<^*_\chi)$
\sn
\item[$(k)$]  if $\kappa=\kappa^{<\partial}$ we can add 
$[N_\eta]^{< \partial} \subseteq N_\eta$.
\end{enumerate}
\end{claim}

\begin{remark}
1) Used in \cite[1.11=L7.6(2),(3),(4)]{Sh:331} and 
\cite[3.23=L7.14,Case 5,clause (k)]{Sh:331} and \cite[3.25=L7.151]{Sh:331}.

\noindent
2) See \cite[Ch.IV]{Sh:f} use 1.10 + the functions witnessing successor.
\end{remark}

\begin{PROOF}{\ref{a48}}
Let $\cS_* \subseteq [\lambda]^{\le \kappa}$ be stationary of
cardinality $\lambda$; why
exists? if $\lambda = \lambda^\kappa$ trivially if just $\lambda =
\cf([\lambda]^{\le \kappa},\subseteq)$ by \cite{Sh:420}.

\noindent
1) We apply \ref{a54} below.  

In detail let
\mn
\begin{enumerate}
\item[$(a)$]  $\kappa = \theta^+,\partial = \aleph_0$ and $\lambda^+$
  here stands for $\lambda$ in \ref{a54}
\sn
\item[$(b)$]  $\cT = \{\eta:\eta$ an increasing sequence of ordinals
  $< \lambda^+\} $
\sn
\item[$(c)$]  if $\eta \in \cT$ then $\bold I_\eta$ is the ideal of
  non-stationary subsets of $\lambda^+$ plus the set $\{\delta <
  \lambda^+:\cf(\delta) \le \kappa\}$
\sn
\item[$(d)$]  $\kappa_\eta = \lambda^+$ for $\eta \in \cT$
\sn
\item[$(e)$]  for some $(g^0,g^1)$ witnessing $\lambda^+$, see below.
\end{enumerate}

\begin{equation*}
\begin{array}{clcr}
\cS = \{u \in [\lambda^+]^{\le \kappa}:& u \text{ is closed under }
g^0,g^1,\kappa +1 \subseteq u \text{ and} \\
  &u \cap \theta \in \theta \text{ and } u \cap \lambda \in \cS_*\}.
\end{array}
\end{equation*}

\mn
Now we can check that the assumptions of \ref{a54} holds hence its
conclusion give the desired conclusion.

\noindent
2) The game proof, but using clause $\oplus_2(g)$ of the conclusion of
   \ref{a54}.

\noindent
3) We could choose $\cS_*$ as the given $\cS$ and use the proof 
of (1).  Of course we can combine part (3) with parts (2),(4) if $\cS$
is as in $(h)$ of \ref{a54}.

\noindent
4) Clause (j) is really proved in \ref{a54}.  As for clause (a) we can
   in the proof of part (j) replace
\mn
\begin{enumerate}
\item[$(a)'$]  $\kappa = \theta^+$ and $\partial$ is the one given,
  \wilog \, regular and use $\lambda' = (\lambda^+)^{< \partial}$ in
  \ref{a54}.
\end{enumerate}
\mn
As $\theta = \kappa^+$ clearly $\alpha < \theta \Rightarrow
|\alpha|^{< \partial} \le \kappa^{< \partial} = \kappa < \theta$ by
the present proof
\mn
\begin{enumerate}
\item[$(d)'$]  in the definition of $\cS$ demand $u$ is closed under
  $h$ as there (exists as we are assuming $\kappa = \kappa^{<
  \kappa}$).
\end{enumerate}
\end{PROOF}

\begin{claim}
\label{a54}
Assume that:
\mn
\begin{enumerate}
\item[$\oplus_1$]  $(a) \quad \theta$ is an uncountable regular cardinal,
\sn
\item[${{}}$]  $(b) \quad ({\cT},\bar{\bold I})$ is a tagged tree,
\sn
\item[${{}}$]  $(c) \quad$ for $\eta \in {\cT},\bold I_\eta$ is a
  normal\footnote{pedantically we should use $\bold I^\dagger_\eta$}
  ideal on some regular uncountable 

\hskip25pt cardinal $\kappa_\eta$,
\sn
\item[${{}}$]  $(d) \quad A_\eta$ is a set of cardinality $<\theta$,
for $\eta \in \cT$
\sn
\item[${{}}$]  $(e) \quad \lambda \ge \Sigma\{\kappa_\eta:\eta \in
  \cT\}$ and $\cS \subseteq [\lambda]^{< \theta}$ is stationary
\sn
\item[${{}}$]  $(f) \quad$ if $\eta 
\triangleleft \nu \in{\cT}$ then $\kappa_\eta\le \kappa_\nu$,
\sn
\item[${{}}$]  $(g) \quad ({\cT},\bar{\bold I}),\langle A_\eta:\eta\in {\cT}
\rangle \in {\cH}(\chi)$ and $x\in{\cH}(\chi)$
\sn
\item[${{}}$]  $(h) \quad$ if $\eta \in \cT$ and $\alpha <
  \kappa_\eta$ then $\cS \rest \alpha$ has cardinality $< \kappa_\eta$
  where $\cS \rest \cU :=$

\hskip25pt $\{u \cap \cU:u \in \cS\}$ and so a sufficient
  condition is 

\hskip25pt $(\forall \alpha < \kappa_\eta)(|\alpha|^{< \theta} <
  \kappa_\eta)$
\sn
\item[${{}}$]  $(i)(\alpha) \quad \bold I_\eta$ is a normal ideal on
  $\kappa_\eta$ 
\sn
\item[${{}}$]  \hskip10pt $(\beta) \quad \{\delta <
  \kappa_\eta:\cf(\delta) < \theta\} \in \bold I_\eta$
\sn
\item[${{}}$]  \hskip10pt $(\gamma) \quad$ if $\eta_1 \ne \eta_2 \in
  \cT$ and $\kappa_{\eta_1} = \kappa_{\eta_2}$ and $\eta_1 \char 94
  \langle \alpha_1\rangle,\eta_2 \char 94 \langle \alpha_2 \rangle \in
  \cT$ then

\hskip52pt  $\alpha_1 \ne \alpha_2$ or at least
  $\cP(\kappa_\eta)/\cI_\eta$
\sn
\item[${{}}$]  $(j) \quad \partial < \theta$ and $\alpha < \theta
  \Rightarrow |\alpha|^{< \partial} < \theta$ and $h:{}^{\partial
  >}\lambda \rightarrow \lambda$ is one to one

\hskip25pt  and $u \in \cS \wedge
  \rho \in {}^{\partial >} u \Rightarrow h(\rho) \in u$.
\end{enumerate}
\mn
\Then \, there is a sequence $\langle N_\eta:\eta\in {\cT}^*\rangle$
such that
\mn
\begin{enumerate}
\item[$\oplus_2$]  $(a) \quad ({\cT},\bar{\bold I}) \le 
({\cT}^*,\bar{\bold I} \rest {\cT}^*)$
\sn
\item[${{}}$]  $(b) \quad N_\eta \prec (\cH(\chi),\in)$ and $x \in
  N_\eta$
\sn
\item[${{}}$]  $(c) \quad$ if $\eta \in \cT^*$
  then $N_\eta \cap \kappa_\eta \in \cS \rest \kappa_\eta$
\sn
\item[${{}}$]  $(d) \quad \eta \in {\cT}^* \Rightarrow A_\eta 
\cup \{x\} \subseteq N_\eta$
\sn
\item[${{}}$]  $(e) \quad \eta \in N_\eta$
\sn
\item[${{}}$]  $(f) \quad \langle N_\eta:\eta\in 
{\cT}^*\rangle$ is a $\Delta$--system, i.e.,
$N_\eta \cap N_\nu=N_{\eta\cap \nu}$ 
\sn
\item[${{}}$]  $(g) \quad$ if $\alpha < \theta \Rightarrow
  2^{|\alpha|} < \kappa_{\langle \rangle}$ then
$\eta,\nu \in \cT \ \&\ \ell g(\eta) = \ell g(\nu) \Rightarrow N_\eta \cong
N_\nu$.
\end{enumerate}
\end{claim}

\begin{remark}
\label{a60}
1) What if $\theta$ is singular? 
Let $\theta = \sum\limits_{\zeta< \partial} \theta_\zeta,
\theta_\zeta$ regular uncountable increasing with
$\zeta,\partial = \cf(\theta)<\theta$.  Now let
$f:{\cT} \rightarrow \partial$ be $f(\eta) = \min\{\zeta:|\bigcup\{A_{\eta
\restriction\ell}:\ell \le \ell g(\eta)\}| < \theta_\zeta\}$ 
and use ?

\noindent
2) Used in the proofs of \cite[1.14=L7.6B]{Sh:331}, \cite[2.15=L7.9]{Sh:331}.
\end{remark}

\begin{PROOF}{\ref{a48}}
\Wilog \, $x$ codes $(\cT,\bar{\cT}),\langle A_\eta:\eta \in
\cT\rangle,\theta,\bar\kappa,\cS$.  Let $\gB$ expand
$(\cH(\chi),\in,<^*_\chi)$ by $x$ and the functions $F_i$ (for $i < \partial$)
where $F_i$ is an $i$-place function from $\cH(\chi)$ to 
$\cH(\chi)$ and $F_i(\ldots,a_j,\ldots)_{j<i} = \langle
a_j:j<i\rangle$ and the functions $G_i$ (for $i<\theta$): 
$G_i(a)$ is: $i$ if $a \in \theta
\backslash i, 0$ if otherwise.

So
\mn
\begin{enumerate}
\item[$(*)_0$]  if $u \subseteq \cH(\chi),|u| < \theta$ then $N =
  \Sk(u,\gB) \prec \gB$ satisfies
\sn
\begin{enumerate}
\item[$\bullet$]  $N \cap \theta \in \theta$
\sn
\item[$\bullet$]  $N$ has cardinality $< \theta$
\sn
\item[$\bullet$]  $N^{< \partial} \subseteq N$.
\end{enumerate}
\end{enumerate}
\mn
Let $\bold N$ be the set of pairs $(\eta,\bar N)$ such that:
\mn
\begin{enumerate}
\item[$(*)^1_{\eta,\bar N}$]  $(a) \quad \eta \in \cT$
\sn
\item[${{}}$]  $(b) \quad \bar N = \langle N_\ell:\ell \le \ell
  g(\eta)\rangle$
\sn
\item[${{}}$]  $(c) \quad N_\ell \prec (\cH(\chi),\in,<^*_\chi)$
\sn
\item[${{}}$]  $(d) \quad x \in N_\eta,\eta \rest \ell \in N_\ell$ and
  $\|N_\ell\| < \theta$
\sn
\item[${{}}$]  $(e) \quad N_\ell$ is the Skolem hull of $N_\ell \cap
  \kappa_{\eta \rest \ell}$ in $\gB$
\sn
\item[${{}}$]  $(f) \quad N_\ell \cap \lambda \in \cS$
\sn
\item[${{}}$]  $(g) \quad N_\ell \subseteq N_{\ell +1}$ (equivalently
  $N_\ell \prec N_{\ell +1}$) and moreover, $N_\ell <_{\kappa_{\eta
  \rest \ell}} N_{\ell +1}$

\hskip25pt  which means $N_\ell \subseteq N_{\ell +1}$
  and $N_\ell \cap \kappa_{\eta \rest \ell} \triangleleft N_{\ell +1} \cap
  \kappa_{\eta\rest \ell}$.
\end{enumerate}
\mn
Let $\bold N_n = \{(\eta,\bar N) \in \bold N:\ell g(\eta) = n+1\}$.

We define a two-place relation $\le_{\bold N}$ on $\bold N$:
\mn
\begin{enumerate}
\item[$(*)_2$]  $(\eta_1,\bar N_1) \le_{\bold N} (\eta_2,\bar N_2)$
  \Iff \, both are from $\bold N$ and $\eta_1 \trianglelefteq
  \eta_2,\bar N_1 \trianglelefteq \bar N_2$.
\end{enumerate}
\mn
Obviously
\mn
\begin{enumerate}
\item[$(*)_3$]  $(a) \quad \bold N$ is non-empty
\sn
\item[${{}}$]  $(b) \quad \le_{\bold N}$ is a partial order on $\bold
  N$, in fact $(\bold N,\le_{\bold N})$ is a tree with $\omega$ levels,

\hskip25pt   the $n$-th level being $\bold N_n$
\sn
\item[${{}}$]  $(c) \quad$ if $(\eta,\bar N) \in \bold N_{n_2}$ and
$n_1 \le n_2$ then $(\eta,\bar N) \upharpoonleft n_1 := 
(\eta \rest n_1,\bar N \rest (n_1+1))$

\hskip25pt  belongs to $\bold N_1$ and is $\le_{\bold N} (\eta,\bar N)$.
\end{enumerate}
\mn
Now we define a function $\rk:\bold N \rightarrow \Ord \cup
\{\infty\}$ by defining when $\rk(\eta,\bar N) \ge \alpha$ by
induction on the ordinal $\alpha$:
\mn
\begin{enumerate}
\item[$(*)_4$]  $\rk(\eta,\bar N) \ge \alpha$ \Iff \, for some
  $n,(\eta,\bar N) \in \bold N_n$ and for every $\beta < \alpha$ there
is $\bold x = \langle (\eta_s,\bar N_s):s \in S\rangle$ such that
\sn
\begin{enumerate}
\item[$(a)$]  $(\eta_s,\bar N_s) \in \bold N_{n+1}$
\sn
\item[$(b)$]  $(\eta,\bar N) \le_{\bold N} (\eta_s,\bar N_s)$ and
  $\rk(\eta_s,\bar N_s) \ge \beta$ for every $s \in S$
\sn
\item[$(c)$]  $\{\eta_s:s \in S\} \in \bold I^+_\eta$
\sn
\item[$(d)$]  if $s_1 \ne s_2 \in S$ then $N_{s_1,n+1} \cap
  N_{s_2,n+1} = N_n$ where $\bar N_s = \langle N_{s,\ell}:\ell < |\bar
  N_s|\rangle$. 
\end{enumerate}
\end{enumerate}
\mn
Clearly $\rk$ is indeed a function from $\bold N$ into $\Ord \cup
\{\infty\}$.
\mn
\begin{enumerate}
\item[$(*)_5$]  if $\rk(\eta,\bar N) = \infty$ for some
$(\eta,\bar N) \in \bold N_0$ then the desired conclusion holds.
\end{enumerate}
\mn
Why?  In short, here we use $\eta \triangleleft \nu \Rightarrow
\kappa_\eta \le \kappa_\nu$ and $\bold I_\eta$ fails
$\kappa^+_\eta$-c.c. and $\bold I^+_\eta$ is a normal ideal on
$\kappa_\eta,\cP(\kappa_\eta)/\bold I^+_\eta$ fails the
$\kappa^+_\eta$-c.c. everywhere (see later on normal ideals on
$[\kappa_\eta]^{< \partial(\eta)}$).  Fully, first we can ignore
$\oplus_2(g)$ as we can apply \ref{1.9}.

Let $\bold N'_\eta = \{(\eta,\bar N) \in \bold N_\eta:\rk(\eta,\bar N)
= \infty\}$.  

Now we shall choose $\cT'_n \subseteq \cT_n := \{\eta \in \cT:\ell g(\eta)
=n\}$ and $\bar N_\eta$ for $\eta \in \cT_n$ such that $(\eta,\bar
N_\eta) \in \bold N$ and $\rk(\eta,\bar N_\eta) = \infty$.
\mn
\begin{enumerate}
\item[$(*)_{5.1}$]  if $n=0$ then $\cT'_0 = \{\langle \rangle\},\bar
  N_{\langle \rangle}$ is such that $(\langle \rangle,\bar N_{\langle
  \rangle}) \in \bold N_0$ and $\rk(\langle \rangle,\bar N_{\langle
\rangle}) = \infty$.
\end{enumerate}
\mn
This holds by the assumption of $(*)_5$
\mn
\begin{enumerate}
\item[$(*)_{5.2}$]  if $\eta \in \cT'_n$ so $\bar N_\eta$ is
well defined then for every ordinal $\alpha$ there is $\bold x_\alpha =
\langle(\eta^\alpha_s,\bar N^\alpha_s):s \in S_\alpha\rangle$
witnessing $\rk(\eta,\bar N_\eta) = \alpha$, hence for some
$\beta = \beta(\eta)$ we have that
$\{\alpha:\bold x_\alpha = \bold x_\beta\}$ is a proper class and let
\sn
\begin{enumerate}
\item[${{}}$]  $(a) \quad \cT'_{n+1} \cap \Succ_{\cT}(\eta) =
\{\eta^{\beta(\eta)}_s:s \in S_{\beta(\eta)}\}$
\sn
\item[${{}}$]  $(b) \quad \bar N_{\eta^{\beta(\eta)}_s} =
  N^\alpha_s$.  

So
\sn
\item[${{}}$]  $(c) \quad \cT'_{n+1} = \cup\{\cT'_{n+1} \cap
\Succ_{\cT}(\eta):\eta \in \cT'_n\}$.
\end{enumerate}
\end{enumerate}
\mn
Clearly
\mn
\begin{enumerate}
\item[$(*)_{5.3}$]  $(a) \quad (\cT,\bar{\bold I}) \le^* (\cT',\bold I)$
\sn
\item[${{}}$]  $(b) \quad$ if $\eta \in \cT$ and $\nu_1 \ne \nu_2 \in
  \Succ_{\cT'}(\eta)$ then $N_{\nu_1} \cap N_{\nu_2} = N_n$.
\end{enumerate}
\mn
Our problem is to find $\cT''$ such that $(\cT',\bar{\bold I}) \le
(\cT'',\bar{\bold I})$ and $\langle N_{\eta,\ell g(\eta)}:\eta \in
\cT''\rangle$ is a $\Delta$-system because then by the assumption on the $\bold
I_\eta$'s, i.e. by $\oplus_1(f)(\gamma)$ we are done.  We still have
to prove the assumption of $(*)_5$
\mn
\begin{enumerate}
\item[$(*)_6$]  there is $(\eta,\bar N) \in \bold N_0$ such that
  $\rk(\eta,\bar N) = \infty$.
\end{enumerate}
\mn
Why?  For every $\eta \in \cT$ and $\alpha < \kappa_\eta$ 
\mn
\begin{enumerate}
\item[$(*)_{6.1}$]  let $\bold N_{\eta,\alpha}$ be $\{\bar
  N:(\eta,\bar N) \in \bold N_{\ell g(\eta)}$ and $N_{\ell g(\eta)} 
\cap \kappa_\eta \subseteq \alpha\}$.
\end{enumerate}
\mn
Now
\mn
\begin{enumerate}
\item[$(*)_{6.2}$]  if $\eta \in \cT$ and $\alpha < \kappa_\eta$ then 
$\bold N_{\eta,\alpha}$ has cardinality $< \kappa_\eta$.
\end{enumerate}
\mn
Why?  Because $|\cS \rest \alpha| < \kappa_\eta$.
\mn
\begin{enumerate}
\item[$(*)_{6.3}$]  If $\eta \in \cT,\alpha < \kappa_\eta,\bar N \in  
\bold N_{\eta,\alpha}$ and $\rk(\eta,\bar N) < \infty$ \then \,
  $C_{\eta,\bar N} \in \bold I_\eta$ where $C_{\eta,\bar N} := \{\beta
  < \kappa_\eta$: there is $\bar N'$ such that $(\eta,\bar N)
  \le_{\bold N} (\eta \char 94 \langle \beta \rangle,\bar N') \in
  \bold N_{\ell g(\eta)+1}$ and $\rk(\eta \char 94 \langle \beta
  \rangle,\bar N') \ge \rk(\eta,\bar N)\}$.
\end{enumerate}
\mn
Why?  By the definition of $\rk(-)$.
\mn
\begin{enumerate}
\item[$(*)_{6.4}$]  if $\eta \in \cT'$ then $C_\eta \in \bold I_\eta$
  where $C_\eta$ is the set of $\beta < \kappa_n$ satisfying at least
  one of the following:
\sn
\begin{enumerate}
\item[$(a)$]  $\cf(\beta) < \theta$
\sn
\item[$(b)$]  $\eta \char 94 \langle \beta \rangle \notin \cT'$
\sn
\item[$(c)$]  for some $\alpha < 1 + \beta$ and $\bar N \in \bold
  N_{\eta,\alpha}$ we have $\beta \in C_{\eta,\bar N}$
\sn
\item[$(d)$]  in the Skolem hull of $\beta \cup \{x\}$ there is an
  ordinal from $[\beta,\kappa_\eta)$.
\end{enumerate}
\end{enumerate}
\mn
Why?  Because $\bold I_\eta$ is a normal ideal on $\kappa_\eta$ and
$(*)_{6.3}$.

Now we choose $\eta_n$ by induction on $n$ such that:
\mn
\begin{enumerate}
\item[$(*)_{6.5}$]  $(a) \quad \eta_n \in \cT'$ has length $n$
\sn
\item[${{}}$]  $(b) \quad$ if $n=m+1$ \then \, $\eta_n = \eta_m \char
  94 \langle \delta_m \rangle$ for some $\delta_m \in \kappa_{\eta_m}
  \backslash C_{\eta_m}$.
\end{enumerate}
\mn
Clearly possible as we are assuming ``$\cS \subseteq [\lambda]^{<
  \theta}$ is stationary" there are $M,u$ such that:
\mn
\begin{enumerate}
\item[$(*)_{6.6}$]  $(a) \quad u \in \cS$
\sn
\item[${{}}$]  $(b) \quad M_u$ is the Skolem hull of $u \cup \{x\}$ in
 $\gB$
\sn
\item[${{}}$]  $(c) \quad \delta_n \in u$ for every $n$.
\end{enumerate}
\mn
Let $N_n$ be the Skolem hull in $\gB$ of $(N \cap \delta_n)
\cup\{x\}$.  Let $\bar N_n = \langle N_\ell:\ell \le n\rangle$.

Now
\mn
\begin{enumerate}
\item[$(*)_{6.7}$]  $(a) \quad (\eta_n,\bar N_n) \in \bold N_n$
\sn
\item[${{}}$]  $(b) \quad$ if $\rk(\eta_n,\bar N_n) < \infty$ then
 $\rk(\eta_n,\bar N_n) > \rk(\eta_{n+1}\bar N_{n+1})$.
\end{enumerate}
\mn
Why?  By the choice of the $C_n$'s.

It follows that $\rk(\eta_0,\bar N_0) = \infty,(\eta_0,\bar N_0) \in
\bold N_0$, so we are done.
\end{PROOF}

\noindent
In \ref{a48}(1) we can replace $\lambda^+$, $\kappa^+$ by $\lambda_1$,
$\kappa_1$, that is
\begin{claim}
\label{a64}
1)  If
\mn
\begin{enumerate}
\item[$(i)$]  $\lambda_1 = \cf(\lambda_1) > \kappa_1 = \cf(\kappa_1) >
\aleph_0$, 
\sn
\item[$(ii)$]  $\alpha < \lambda_1 \Rightarrow \cov (|\alpha|,\kappa_1,
\kappa_1,2) < \lambda_1$,
\sn
\item[$(iii)$]  $\chi > \lambda^+$ and $x \in {\cH}(\chi)$,
\end{enumerate}
\mn
\then \, we can find $\langle N_\eta:\eta\in {\cT}\rangle$ such that
\mn
\begin{enumerate}
\item[$(a)_1$]  $\cT$ is a subtree of ${}^{\omega>}(\lambda_1)$, each $\eta\in
\cT$ is strictly increasing, 
\sn
\item[$(b)_1$]  $N_\eta \prec ({\cH}(\chi), \in, <^*_\chi)$,
\sn
\item[$(c)_1$]  $x,\lambda_1,\kappa_1$ belong to $N_\eta$, $N_\eta 
\cap \kappa_1 \in \kappa_1$ and $\|N_\eta\|=|N_\eta\cap \kappa_1|$, 
\sn
\item[$(d)$]  $\nu \triangleleft \eta \in T \Rightarrow N_\nu \prec
  N_\eta$, 
\sn
\item[$(e)$]  $N_\eta\cap N_\nu = N_{\nu\cap \eta}$ for
  $\eta,\nu\in\cT$, 
\sn
\item[$(f)$]  $\eta\in N_\eta$,
\sn
\item[$(g)$]  if $\eta_\ell \char 94 \langle\alpha_\ell\rangle \in
\cT$ for $\ell=1,2$ and $\alpha_1 < \alpha_2$ \then 
\[
\sup(N_{\eta_1 \char 94 \langle \alpha_1\rangle}\cap \alpha^+)< \min(
N_{\eta_2 \char 94 \langle\alpha_2\rangle}\cap \lambda^+, N_{\eta_2}).
\]
\end{enumerate}
\mn
2)  If in addition $\alpha< \lambda_1 \Rightarrow 
|\alpha|^{< \kappa_1} < \lambda_1$ \then\ we can add 
\mn
\begin{enumerate}
\item[$(h)$]   if $\eta,\nu \in \cT$ have the same length then there is an
isomorphism from $N_\eta$ onto $N_\nu$, call it $f_{\eta,\nu}$, and it maps
$x$ to itself, so 

\[
\eta,\nu \in \lim(\cT) \Rightarrow \bigcup\limits_{n < \omega} 
N_{\eta \rest n} := N_\eta \cong N_\nu := \bigcup\limits_{n < \omega} 
N_{\eta \rest n}
\]
\end{enumerate}
\mn
3)  If $\bar\cS=\langle {\cS}_\alpha: \alpha< \lambda_1\rangle$ is
$\subseteq$-increasing with $\alpha,{\cS}_\alpha \subseteq 
[\alpha]^{<\kappa_1},\alpha\in a \in 
S_\beta \Rightarrow a \cap \alpha \in \cS_\alpha$ then we can demand
\mn
\begin{enumerate}
\item[$(i)$]  $N_\eta \cap \lambda_1\in \bigcup {\cS}_\alpha$.
\end{enumerate}
\mn
4)  We can further demand
\mn
\begin{enumerate}
\item[$(j)$]  $N_\eta$ is the Skolem hull of $\{x,\eta,\kappa_1,\lambda_1\}
\cup (N_\eta\cap \kappa_1)\cup N_{\langle\rangle}$. 
\end{enumerate}
\mn
5)  If $({\cT}_0,\bar{\bold I})$ is a tagged tree, 
$\bold I_\eta$ a normal ideal on $\lambda_1$ such that 
$\{\delta:\cf(\delta) < \kappa_1\} \in \bold I_\eta$ then
we can demand $({\cT}_0,\bar{\bold I}) \le^* ({\cT},\bar{\bold I})$.
\end{claim}

\begin{PROOF}{\ref{a64}}
Similarly to \ref{a48} by \ref{a54}.
\end{PROOF}

\begin{remark}
\label{a57}
The isomorphism is unique, hence if the isomorphism is called 
${\bold f}_{\eta,\nu}$ then ${\bold f}_{\eta_0,\eta_1} = 
{\bold f}_{\eta_0,\eta_1} \circ {\bold f}_{\eta_1,\eta_0}$ when 
they are well defined.
\end{remark}

\begin{claim}
\label{1.17}
Suppose that
\mn
\begin{enumerate}
\item[$(a)$]   $\lambda$ is singular, $\kappa = \cf(\lambda) >\aleph_0$,
\sn
\item[$(b)$]  $f$ is a function from ${}^{\omega>}\lambda$ to finite subsets of
${}^{\omega \ge}\lambda$ or even subsets of ${}^{\omega\ge} \lambda$ of
cardinality $< \cf(\lambda)$, 
\sn
\item[$(c)$]  $\lambda = \sum\limits_{i<\kappa} \lambda_i$, where $\lambda_i$
is (strictly) increasing and continuous with $i< \kappa$; 
\sn
\item[$(d)$]  $S\subseteq \{i < \kappa:\cf(i) = \aleph_0\}$ is stationary,
\sn
\item[$(e)$]  for $i \in S$ we have $\lambda_i = \sum\limits_{n<\omega}
\lambda_{i,n}$, where $\kappa < \lambda_{i, 0}$, and 
$(\forall n)(\lambda_{i,n} < \lambda_{i,n+1}\ \&\ \cf(\lambda_{i,n})=
\lambda_{i,n})$ 
\sn
\item[$(f)$]  $\bold I^n_\mu$ is a $\kappa^+$--complete ideal on 
$\mu$ containing the co-bounded subsets of $\mu$ for $\mu$ regular $<\lambda$ 
\sn
\item[$(g)$]   if $i_1,i_2\in S$ and 
$\{j:\lambda_j<\lambda_{i_1,n}\}=\{j:\lambda_j< \lambda_{i_2,n}\}$ 
when $i_1,i_2\in S$ and $n<\omega$ then $\lambda_{i_1,n}
=\lambda_{i_2,n}$ and $\bold I^n_{\lambda_{i_1,n}} =
\bold I^n_{\lambda_{i_2,n}}$.
\end{enumerate}
\mn
\Then \, there is a closed unbounded $C\subseteq \kappa$ such that for
each $i \in C \cap S$ there is a $\cT$ such that:
\mn
\begin{enumerate}
\item[$(*)_1$]  $\cT \subseteq \bigcup\limits_{n<\omega}\prod\limits_{m<n}
\lambda_{i, m}$,  $\langle\rangle \in \cT$, and $\cT$ is closed under
initial segments; 
\sn
\item[$(*)_2$]  if $\eta\in \cT $ and $\lg(\eta)=n$ then $\{\alpha<\lambda_{i,
n}: \eta \char 94 \langle\alpha\rangle\in \cT\} \ne \emptyset  
\mod \bold I^n_{\lambda_{i, n}}$; 
\sn
\item[$(*)_3$] if $\eta\in \cT$ then  $f(\eta) \subseteq {}^{\omega\ge} 
\lambda_i$.  
\end{enumerate}
\end{claim}

\begin{remark}
Claim \ref{1.17} is used in \cite[1.11]{Sh:331}.
\end{remark}

\begin{PROOF}{\ref{1.17}}
This is a variant of \ref{a45}. 
For each $\eta \in {}^{\omega>}\lambda$ choose $g(\eta)<\kappa$ so
that $f(\eta) \subseteq \bigcup \{{}^{\omega\ge}\zeta:\zeta<
\lambda_{g(\eta)}\}$. Then instead of $(*)_3$, it suffice to demand
\mn
\begin{enumerate}
\item[$(*)_3'$]   $\forall \eta\in \cT (g(\eta)< i)$.
\end{enumerate}
\mn
Now we define a game $\Game_i$ for each $i\in S$: the game is of length
$\omega$, and in the $n$-th move, the second player chooses $A_n\in
\bold I^n_{\lambda_{i,n}}$ with $|A_n|< \lambda_{i, n}$, and the first player
chooses $\alpha_n \in \lambda_{i,n}$. The first player wins if [$\ell<n\  
\Rightarrow\ \alpha_\ell<\alpha_n]$, $\alpha_n \notin A_n$, and
$g(\langle \alpha_0, \ldots, \alpha_n\rangle)\leq i$; otherwise the second
player wins. 

Now
\mn
\begin{enumerate}
\item[$(*)_4$]  if $i \in S$, $g(\langle\rangle) \le i$, and the first player
has a winning strategy, then a tree ${\cT} = {\cT}_i$ as desired exists. 
\end{enumerate}
\mn 
Why? Let $i \in S$ be as in $(*)_4$ and let ${\bold f}_i$ be 
a winning strategy for the first player in the game $\Game_i$. 
Thus for $n < \omega$ and $\bar A \in {}^{n+1}{\cP}(\lambda)$ 
such that $\forall m\leq n (A_m \in \bold I^n_{\lambda_i,m})$ 
we have ${\bold f}_i({\bar A})\in \lambda_{i, n},{\bold f}_i({\bar A}
\rest (m+1)) < {\bold f}_i({\bar A})$ for all $m < n,{\bold f}_i
({\bar A})\notin A_n$, and $g(\langle {\bold f}_i({\bar A} \rest 1),
\ldots,{\bold f}_i({\bar A} \restriction n)\rangle)\le i$. 
Then $\cT = \{\langle {\bf f}_i({\bar A}\restriction 1),\ldots,
{\bold f}_i({\bar A}\restriction (n+1))\rangle$: such ${\bar A}\}\cup
\{\langle \rangle\}$ is as desired. Thus we may assume toward contradiction 
\mn
\begin{enumerate}
\item[$(*)_5$]  $S'=\{i \in S:$ the first player does not have a winning
strategy for $\Game_i\;\}$ is a stationary subset of $\cf(\lambda)$. 
\end{enumerate}
\mn
Now, the game $\Game_i$ is open, so by the Gale-Stewart theorem it is
determined.  Hence for each $i \in S'$ we may choose a 
winning strategy ${\bold f}_i$ for the second player. 

Thus 
\mn
\begin{enumerate}
\item[$(*)_6$]  if $n<\omega$ and  $\eta\in\prod\limits_{m<n} \lambda_{i,m}$
then  ${\bold f}_i(\eta) \in \bold I^n_{\lambda_{i,n}}$;
\sn
\item[$(*)_7$]  for any $\eta \in \prod\limits_{m<\omega} \lambda_{i, m}$ one
of the following occurs: 
\sn
\begin{enumerate}
\item[$(a)$]  $\exists \ell<n<\omega$ $(\eta(\ell)\geq \eta(n))$,
\sn
\item[$(b)$]  there is $n<\omega$ such that $\eta(n) \in 
{\bold f}_i(\eta\restriction n)$,
\sn
\item[$(c)$]  there is $n<\omega$ such that $g(\eta\restriction n)>i$.
\end{enumerate}
\end{enumerate}
\mn
Now choose a regular $\chi> \aleph_0$ so that $g$, $\langle {\bold f}_\delta:
\delta\in S'\rangle$, $\langle \lambda_i: i< \cf(\lambda)\rangle$,  $\langle 
\bold I^n_\mu:\mu<\lambda$ regular, $n<\omega\;\rangle$ and $\langle \langle
\lambda_{i,n}: n<\omega\rangle: i\in S'\rangle$ belong to ${\cH}(\chi)$. 
Remember ${\cH}(\chi)$ is the family of sets with the transitive closure
of cardinality $<\chi$, and that $({\cH}(\chi),\in)$ is a model of
$\ZFC^-$. Let $<^*_\chi$ be a well-ordering of ${\cH}(\chi)$. 

For all $\delta<\kappa$ let $A_\delta$ be the closure of $\delta\cup\{x\}$
under Skolem functions within the structure $({\cH}(\chi),\in,<^*_\chi)$.
Then $C=\{\delta< \kappa: A_\delta \cap \kappa=\delta\}$ is a closed unbounded
subset of $\kappa$. Thus there is $\delta\in S'\cap C$ and an elementary
substructure $(N, \in, <)$ of $({\cH}(\chi), \in, <^*_\chi)$ such that
$|N|<\kappa$ and $N\cap \kappa = \delta$, with $x \in N$. 
Clearly $\lambda_{\delta,m},\bold I^m_{\lambda_{\delta,m}}$ 
belong to $N$ for each $m$ (by assumption (e)). 
However $\delta\notin N$, hence $\{\lambda_{\delta,m}:
m<\omega\}\notin N$ though it is a subset of $N$.
 
Now we define $\eta=\langle\alpha_n:n<\omega\rangle\in \prod\limits_{m<
\omega} \lambda_{\delta, m}$ so as to contradict $(*)_7$. Suppose
$\alpha_m\in N$ has been constructed for all $m<n$. Using elementarity and
absoluteness of suitable formulas we see that the set 

\begin{equation*}
\begin{array}{clcr}
A^* = \bigcup\{{\bold f}_j(\langle
\alpha_0,\ldots,\alpha_{n-1}\rangle): &j \in S' \text{ and}\\
  & \lambda_{j,0} = \lambda_{\delta,0},\ldots,\lambda_{j,n-1} = 
\lambda_{\delta,n-1},\lambda_{j,n} = \lambda_{\delta,n}\}.
\end{array}
\end{equation*}

\mn
belongs to $\bold I^n_{\lambda^+_{\delta,n}}$ (being the union of $\le
\kappa$ sets each from $\bold I^n_{\lambda^+_{\delta,n}}$) and 
belongs to $N$.  Since $\exists \alpha(\alpha_{n-1} < \alpha < 
\lambda_{\delta,n}$ and $\alpha \notin A^*)$ holds in $({\cH}(\chi),
\in,<^*_\chi)$, it holds in $(N,\in, <^*_\chi)$ and this gives 
$\alpha_n$.  This completes the construction, and
it is easily seen that $(*)_7$ is contradicted.
\end{PROOF}

\begin{remark}
1) We can interchange the quantifier in \ref{1.17};
one club (C) of $\cf(\lambda)$ is O.K. for every appropriate 
$\langle\langle
\lambda_{\delta,n}:n<u\rangle:\delta<\cf(\lambda)\rangle$.

\noindent
2) If $\lambda_{\delta,n}=\eta_\delta(n)$ and $\langle \Rang(\eta_\delta):
\delta\in S\rangle$ guess clubs of $\cf(\lambda)$ \then \, we can add 
$\eta\in \prod\limits_{\delta,\ell} \Rightarrow g(\eta)
>\lambda_{\delta,n}$.

\noindent
3)  We can get in this direction more results. If $2^{\cf(\lambda)} 
< \lambda$, $\lambda_{i+1}$ regular, then we can find a closed unbounded set
$\{\alpha(i): i< \cf(\lambda)\},\alpha(i+1)$ successor and ${\cT}\subseteq
{}^{\omega>}\lambda$, such that: $\langle \rangle \in {\cT},
\eta\in {\cT},\max[\Rang(\eta)]<\lambda_{i+1}< \lambda_j$ implies 
$\{\alpha< \lambda_j:\eta \char 94 \langle \alpha\rangle \in {\cT}\}$ 
has power $\lambda_j$, and implies
also $g(\eta)<j$.(For each club $C$ of $\cf(\lambda)$ we define a
game, etc.).

\noindent
4) In (3) we can replace ``$2^{\cf(\lambda)}<\lambda$" by ``there is
a family ${\cP}$ of closed unbounded subsets of $\cf(\lambda)$ 
such that $|{\cP}|< \lambda$, and every closed unbounded subset 
of $\kappa$ contains one of them".

\noindent
5) On the other hand, if $\mu= \mu^{<\mu}$ in ${\bold V}$ let us 
add $\lambda>\mu$ generic closed unbounded subsets of $\mu$ (by $\bbP
=\{f:\Dom(f)$ a subset of $\lambda$ of power $<\mu,f(i)$ 
the characteristic function of a closed bounded subset of $\mu\}$, 
and let $\name{C}_i$ the following ${\bbP}$-name of a club of $\mu$: 
the characteristic function of $\name{C}_i$ is $\bigcup
\{f(i): f$ in the generic set$\}$). Let ${\bold G}$ be a subset of ${\bbP}$ 
generic over ${\bold V}$ and in ${\bold V}[{\bold G}]$ let
$\{C_\eta: \eta\in {}^{\omega>}\lambda\}$ be another enumeration
of $\{\name{C}_i [{\bf G}]: i<\lambda\}$, and define $g$:

\[
g(\eta \char 94 \langle\alpha\rangle) = \min\{i<\cf(\lambda): i\in C_\eta,
\lambda_i>\alpha\}.
\]

\mn
Clearly for this $g$ the conclusion of remark (4) fails.
\end{remark}
\newpage

\section{On unique linear orders}

Hausdorff has introduced and investigated the class of scattered
linear orders (see \ref{b36}).  Galvin and Laver \cite{Lav71}
investigate the class $M$ of 
linear orders which are a countable union of scattered
linear orders.  They were interested in linear orders up to
embeddability inside the class $M =
\cup\{M_{\lambda,\mu_1,\mu_2}:\mu_1,\mu_2$ regular uncountable
$\lambda \in \mu_1 + \mu_2\}$ where $M_{\lambda,\mu_1,\mu_2}$ is the class
of linear orders $M$ of cardinality $\lambda$ with no increasing sequences of
length $\mu_1$ and no decreasing sequences of length $\mu_2$.
Galvin defined $M_{\lambda,\mu_1,\mu_2}$ and proved the existence of a
universal member.

Laver, solving a long standing conjecture of Fra\"iss\'e, and using the
theory of better quasi orders of Nash Williams proved that the $M$ is well
and even better quasi ordered.  In \cite[pp.308,309]{Sh:220}, 
we continue this investigation being interested in the 
uniqueness of such orders.  We do more here.
Invariants related to the $g_i$ here are investigated in
\cite[Ch.VIII,\S3]{Sh:a} and better in \cite{Sh:E59}, and also in
Droste-Shelah \cite{DrSh:195}, \cite{DrSh:743}.  
This is continued being interested in uniqueness.
\bigskip

\subsection {Classes of Coloured Linear Orders}\
\bigskip

\begin{discussion}
\label{b4}
1)  We may waive ``union of countably many scattered subsets", 
and essentially allow a family of $\le \lambda$ isomorphism types 
of linear orders as basic orders.  So ignoring trivialites 
they are neither well ordered nor anti-well ordered; we
lose stability but can retain everything else. 

\noindent
2) Below in \ref{b64} we may start with closed enough set $\cS \subseteq
\cP(X),|\cS| \le \lambda$.

\noindent
3) Another way to get many of the properties is to build such $N$ of larger
cardinality, so e.g. saturated dense linear orders exists and then use the
L\"owenheim--Skolem argument.
\end{discussion}

\begin{context}
\label{b6}
If not said otherwise, in this subsection we use a fix context $\bold c =
(\lambda,\mu_1,\mu_2,\alpha^*,g_1,g_2) = (\lambda^{\bold c},\mu^{\bold
  c}_1,\mu^{\bold c}_2,\alpha^{\bold c}_*,g^{\bold c}_1,g^{\bold
  c}_2)$ which means it is as in \ref{b8}.
\end{context}

\begin{definition}
\label{b8}
1) We say $\bold c$ is a context if it consists of
  $\lambda,\mu_1,\mu_1,g_1,g_2$ (and $\cF_1,\cF_2,\cF^+_1,\cF^+_2$
  defined from them), \when \,
\mn
\begin{enumerate}
\item[$(a)$]  $\lambda$, $\mu_1$, $\mu_2$ are (infinite) cardinals with 
$\mu_1,\mu_2$ being uncountable regular such that $\lambda^+ = 
\max\{\mu_1,\mu_2\}$ and $\alpha^*=\alpha(*)< \lambda^+$, $\alpha^*
\ge 1$
\sn
\item[$(b)$]  for $\ell=1, 2$ we have $g_\ell$, a function
from $\alpha^*$ into ${\cF}_\ell := \{h: h$ a function from some uncountable
$\theta \in \Reg \cap \mu_\ell$ into $\alpha^*\}$ such that
$\{\Dom(h):h \in \Rang(g_\ell)\}$  is unbounded in $\Reg\cap \mu_\ell$. 
(Hence $\lambda = \sup\{\Dom(g_\ell(\alpha)):\ell\in \{1,2\}$ 
and $\alpha<\alpha^*\}$)
\end{enumerate}
\mn
2) In addition if $\ell \in \{1,2\},\alpha < \alpha^*,h=g_\ell(\alpha)$
\then \,: $h \in \cF^+_\ell$ where $\cF^+_\ell$ is the set of $h \in
\cF_\ell$ satisfying
\mn
\begin{enumerate}
\item[$\boxdot^\ell_h$]   $h \in \cF_\ell$ and if 
$\delta < \Dom(h)$ is a limit 
ordinal of uncountable cofinality and 
$\beta=h(\delta)$ and $\langle \epsilon_i:i<\cf(\delta)\rangle$ 
is an increasing continuous sequence with limit $\delta$ then the set
$\{i<\cf (\delta):(h(\delta))(\epsilon_i) = (g(\beta))(i)\}$ 
contains a club of $\cf(\delta)$.
\end{enumerate}
\mn
3) For notational simplicity assume $\alpha^* \le \lambda$.

\end{definition}

\begin{notation}
\label{b12}
1) For a linear order $M= (A,<)$ let $M^*$ be $(A,>)$, i.e., its
   inverse.

\noindent
2) Below $K = K_{\bold c} = K(\bold c)$ and similarly for other
   versions of $K$.

\noindent
3) Properties and Notations defined for linear orders, can be applied
   to expansions of linear orders (here mainly $N \in K$ or $N \in K^{\all}$).
\end{notation}

\begin{definition}
\label{b16}
$K= K^{\hom} = K(\lambda,\mu_1,\mu_2,\alpha^*,g_1,g_2)$ is the family of models
$N=(M,P_\alpha)_{\alpha< \alpha(*)}$ such that: 
\mn
\begin{enumerate}
\item[$(i)$]   $M$ is a linear order,
\sn
\item[$(ii)$]  $M$ is the union of $\aleph_0$ scattered suborders, i.e., 
$|M|$, the universe (=set of elements of $M$) is
$\bigcup\limits_{n\in\omega} A_n$, where each $M\restriction A_n$ is
scattered (see Definition \ref{b20} below),
\sn
\item[$(iii)$]  each $P_\alpha$ is a dense subset of $M$,
\sn
\item[$(iv)$]  $\langle P_\alpha:\alpha<\alpha(*)\rangle$ is 
a partition of $M$, 
\sn
\item[$(v)$]  every increasing sequence in $M$ has length 
$<\mu_1$, but for each $\alpha< \alpha^*$ in every open 
interval of $M$ there is an increasing
sequence of length $\Dom(g_1(\alpha))$, (hence any $\alpha<\mu_1$ is
O.K.) 
\sn
\item[$(vi)$]  every decreasing sequence in $M$ has length $<\mu_2$, but for
each $\alpha<\alpha^*$ in every open interval there is an decreasing sequence
of any length $\Dom(g_2(\alpha))$, (hence any $\alpha<\mu_2$ is O.K.) 
\sn
\item[$(vii)$]  if $\langle a_i: i<\kappa\rangle$ is an increasing bounded
sequence in $M$, $\aleph_0<\kappa\in\Reg\cap\mu_1$ \then\ for some club $C$ of
$\kappa$, for every $\delta\in C\cup \{\kappa\}$, $\{a_i: i<\delta\}$ has a
least upper bound in $M$,
\sn
\item[$(viii)$]  if $\langle a_i:i<\kappa\rangle$ is a decreasing sequence in
$M$ bounded from below and $\aleph_0<\kappa\in\Reg\cap\mu_2$ \then \, for 
some club $C$ of $\kappa$ for every $\delta\in C\cup\{\kappa\}$ we have:
$\{a_i:i<\delta\}$ has a greatest lower bound in $M$,
\sn
\item[$(ix)$]  if $x \in P_\alpha$, $g_1(\alpha)= h$ \then \, 
$N_{<x} = N \restriction \{y: y<^{M} x\}$ and $M_{<x} = N_{<x} \rest \{<\}$
 has cofinality 
$\Dom(h)$ and if $\Dom(h)>\aleph_0$ then it has up-type $h$ which
  means that 
\sn
\begin{enumerate}
\item[${{}}$]  $(*)^1_{N_{< x,h}} \quad$ there is an increasing continuous 
sequence $\bar y = \langle y_\epsilon:\epsilon <$

\hskip40pt  $\Dom(h)\rangle$ 
in $M_{<x}$ such that $y_\epsilon\in P_{h(\epsilon)}$ for a club of 

\hskip40pt $\epsilon\in\Dom(x)$ and $\{y_\epsilon:\epsilon<\Dom(h)\}$
is 

\hskip40pt unbounded from above in $M_{<x}$, 
\end{enumerate}
\sn 
\item[$(x)$]  if $x\in P_\alpha$, $g_2(\alpha) = h$, 
\then \, $N_{>x} = N \restriction \{y: x <^M y\}$ satisfies: $(M_{>x})^*$,
the inverse of $M_{>x}$, has cofinality $\Dom(h)$ and if $\Dom(h)>\aleph_0$
then $N_x$ has down-type $g_2(\alpha)$ which means that  
\sn
\begin{enumerate}
\item[${{}}$]  $(*)^2_{N_{>x,h}} \quad$ there is a decreasing continuous 
sequence $\bar y=\langle y_\epsilon:\epsilon <$

\hskip40pt $\Dom(h)\rangle$ in 
$M_{>x}$ such that $y_\epsilon \in P_{h(\epsilon)}$ for a club of 

\hskip40pt $\epsilon \in \Dom(h)$ and
 $\{y_\epsilon:\epsilon<\Dom(h)\}$ is

\hskip40pt  unbounded from below in $M_{>x}$. 
\end{enumerate} 
\end{enumerate} 
\end{definition}

\begin{definition}
\label{b20}
1)  For a linear order $M$ we define when 
$\Dep(M) \ge \alpha$ by induction on $\alpha$. If $\alpha=0,\Dep(M) \ge 
\alpha$ for any linear order $M$, even the empty one. If $\alpha=1,
\Dep(M) \ge \alpha$ \underline{if and only if} 
$M$ is non-empty.  If $\alpha$ is limit then $\Dep(M) \ge \alpha$ if
and only if $\Dep(M) \ge \beta$ for every $\beta<\alpha$. If
$\alpha=\beta+1$ then $\Dep(M) \ge \alpha$ if and only if $M$ can be
represented as $M_1+M_2$ where $\Dep(M_1) \ge \beta$ and $\Dep(M_2) \ge
\beta$. 

\noindent
2) We let $\Dep(M)=\alpha$ if and only if $\Dep(M) \ge \alpha$ and
$\Dep(M) \ngeq \alpha+1$. 

We say that $\Dep(M)=\infty$ if $\Dep(M) \ge \alpha$ for all ordinals
$\alpha$.

\noindent
3) A linear order $M$ is scattered if $\Dep(M)<\infty$ equivalently
   (by Hausdorff), the rational order cannot be embedded into $M$.

\noindent
4)  If $N$ is an expansion of a linear order \then \, $\Dep(N)$
means $\Dep(|N|,{<}^M)$. 
\end{definition}

\begin{definition}
\label{b24}
1)  Let $K^* = K^{\all}$ be the class of $N=(M,P_\alpha)_{\alpha<\alpha(*)}$
satisfying clauses (i), (ii), (iv) of Definition \ref{b16}, and the first
half of (v), the first half of (vi), and (vii), (viii)
and clause (ix) for $x$ such that $x$ is neither the first element of
$M$ nor the immediate 
successor of any $y\in M$ and clause (x) for $x$ which are neither 
last nor the immediate predecessor of some $y\in M$ .

\noindent
2)  For $N=(M,P_\alpha)_{\alpha<\alpha(*)}\in K$ and $x\in M$ let 
\mn
\begin{enumerate}
\item[$(a)$]   $P^N_\alpha=P_\alpha$, $<^N=<^M$,
\sn
\item[$(b)$]  $N_{>x},M_{>x}$ and $N_{<x},M_{<x}$ be as in clauses 
(ix), (x) of Definition \ref{b16}
\sn
\item[$(c)$]  so 
$N_{>x}=(M_{>x},P^N_\alpha\cap M_{>x})_{\alpha<\alpha(*)}$ and 
$N_{<x}=(M_{<x},P^N_\alpha\cap M_{<x})_{\alpha<\alpha(*)}$.
\end{enumerate}
\mn
3) For $h_1 \in \cF^+_1,h_2 \in \cF^+_2$ (that is, $h_1 \in 
{\cF}_1$, $h_2\in {\cF}_2$ satisfying $\boxdot^1_{h_1},
\boxdot^2_{h_2}$ from Definition \ref{b8})
let $K_{h_1, h_2} = K^{\hom}_{h_1,h_2} = 
K(\lambda,\mu_1,\mu_2,\alpha^*,g_1,g_2,
h_1,h_2)=K_{h_1,h_2}(\lambda,\mu_1,\mu_2,\alpha^*,g_1,g_2)$ 
be the family of $N = (M,P_\alpha)_{\alpha< \alpha(*)}\in K$ such that $N$ 
satisfies $(*)^1_{N,h_1}$ of clause (ix) of Definition \ref{b16}, 
and $(*)^2_{N,h_2}$ of clause (x) of Definition \ref{b16}.

\noindent
4) For $h_1 \in \cF^+_1,h_2 \in \cF^+_2$ let $K^*_{h_1,h_2} =
   K^{\all}_{h_1,h_2}$ be the family of $N = (M,P_\alpha)_{\alpha <
   \alpha(*)} \in K^*$ such that $N$ satisfies $(*)^1_{N,h_1}$ of
clause (ix) of \ref{b16} and $(*)^2_{N,h_2}$ of clause (x) of
   \ref{b16}.

\noindent
5) Let $K^\otimes = K^{\vhm} = \{M \in K^{\all}$: if $N$ has no last
element then $(*)^1_{M,h}$ for some $h \in \cF^+_1$ and if $N$ has no
first element then $(*)^2_{M,k}$ for some $h \in \cF^+_2\}$.

\noindent
6) $K^\oplus = \cup\{K^*_{h_1,h_2}:h_2 \in \cF^+_1,h_2 \in \cF^+_2\}$.
\end{definition}

\begin{definition}
\label{b28}
1)  For $N_i\in K^*(i < \alpha)$ then 
$N_0 + N_1$ and $\sum\limits_{i< \alpha} N_i$ are
defined naturally, as well as $N_0 \times \alpha$.

\noindent
2)  Similarly for anti-well ordered sums.
\end{definition}

\begin{claim}
\label{b32}
1)  If $N \in K_{h_1,h_2}$ (so $h_1 \in \cF^+_1,h_2 \in \cF^+_2$) 
and $x \in P^N_\alpha$ \then \,:
\mn
\begin{enumerate}
\item[$(i)$]  $N_{<x} = N \rest \{y \in N:y <^M x\} \in K_{g_1(\alpha),
  h_2}$ and 
\sn
\item[$(ii)$]  $N_{>x}= N \rest \{y \in N:x <^M y\} \in K_{h_1, g_2(\alpha)}$.
\end{enumerate}
\mn
2) If $N \in K$ and $I$ is a convex non-empty subset of 
$M$ with neither last nor
first element and $M \restriction I$ satisfies 
$(*)^1_{M,h^*_1}$ of clause (ix) of Definition
\ref{b16} for $h^*_1\in H_1$, and $(*)^2_{M,h^*_2}$ of 
clause (x) of Definition \ref{b16} for
$h^*_2\in H_2$ \then \, $M \restriction I\in K_{h^*_1, h^*_2}$.

\noindent
3) If $N = (M,P_\alpha)_{\alpha<\alpha^*}\in K_{h_1,h_2}$ then 
$N^* = (M^*,P_\alpha)_{\alpha<\alpha^*}\in K_{h_2,h_1}$.

\noindent
4) If $N \in K^{\all}$ and $I$ is a convex subset of $N$  \then \, $N
\rest I \in K^{\all}$.  Moreover, $N \in K^\otimes \Rightarrow N \rest
I \in K^\otimes$.
\end{claim}

\begin{PROOF}{\ref{b32}}
Straightforward.
\end{PROOF}

\noindent
Recall
\begin{claim}
\label{b36}
1) The family of scattered order types is the closure of the singletons
under well ordered sums and inverse of well ordered sums.

\noindent
2) If $M_1 \subseteq M_2$ \then \, $\Dep(M_1) \le \Dep(M_2)$.

\noindent
3)  If $M$ is a scattered linear order \then \, one of the following holds:
\mn
\begin{enumerate}
\item[$(a)$]  $M$ is a singleton,
\sn
\item[$(b)$]  for some $x \in M$ we have $\Dep(M_{<x}) < \Dep(M)$ and 
$\Dep(M_{>x}) < \Dep(M)$,
\sn
\item[$(c)$]   there is an increasing unbounded sequence $\langle x_i:i<\kappa
\rangle$ in $M$, with $\kappa$ a regular cardinal, such that 
\[
i<\kappa \Rightarrow \Dep(M_{<x_i}) < \Dep(M),
\]
\sn
\item[$(d)$]  there is a decreasing sequence $\langle x_i:i<\kappa\rangle$ in
$M$ unbounded from below and such that 

\[
i<\kappa \Rightarrow \Dep(M_{>x_i}) < \Dep(M).
\]
\end{enumerate}
\end{claim}

\begin{claim}
\label{b40}
For any $h_1\in {\cF}_1,h_2 \in {\cF}_2$ we have $K_{h_1, h_2} \ne \emptyset$.
\end{claim}

\begin{PROOF}{\ref{b40}}

First
\mn 
\begin{enumerate}
\item[$\circledast_1$]  there is a scattered $N \in K^*$ such that: 
$N$ has a first element, a last element and $P^N_\alpha \ne 
\emptyset$ for every $\alpha<\alpha^*$.
\end{enumerate}
\mn
[Why?  Recall that $\lambda^+ = \Max\{\mu_1,\mu_2\}$ so let 
$\ell\in\{1,2\}$ be such that $\lambda^+ = \mu_\ell$.  For $\theta \in \{\Dom
(g_\ell(\alpha)):\alpha<\alpha^*\}$, let $\alpha_\theta<\alpha^*$ be 
minimal such that $\theta=\Dom(g_\ell (\alpha_\theta))$.  Now we define 
$N_\theta=(M_\theta, P_{\theta,\alpha})_{\alpha<\alpha(*)}$, i.e. 
$P^{N_\theta}_\alpha = P_{\theta,\alpha}$, as follows:
\mn
\begin{enumerate}
\item[$(a)$]  $M_\theta$ is $(\theta+1,<)$ if $\ell=1$ and is its
  inverse if $\ell=2$
\sn
\item[$(b)$]   for $\epsilon \in \theta+1,
\epsilon\in P_{\theta,\alpha}$ \underline{if and only if} 
$\epsilon=\theta$ and $\alpha=\alpha_\theta$ or 
$(\epsilon<\theta)$ is a limit ordinal and 
$(g(\alpha_\theta))(\epsilon)=\alpha)$ or ($\epsilon = \alpha+1$
so $\epsilon$ is a successor ordinal).
\end{enumerate}
\mn
If $\lambda$ is regular then by \ref{b8} we can choose $\theta=\lambda$ 
and we are done as $\alpha^* \le \lambda$. 
If $\lambda$ is singular we can find an increasing sequence 
$\langle\theta_i:i<\cf(\lambda)\rangle$, with limit $\lambda$,
$\theta_i=\Dom(g_\ell(\alpha_{\theta_i})),\theta_0 > \cf(\lambda)$, 
and we combine them by inserting $N_{\theta_i}$ in the $i$-th
open interval of $N_{\theta_0}$, i.e. in $(i,i+1)_{N_{\theta_0}}$.]
\mn
\begin{enumerate}
\item[$\circledast_2$]  there is a scattered 
$N \in K^*$ such that: $N$ has a first element, $N$ 
has a last element and for every $\ell \in \{1,2\}$
and $\theta = \cf(\theta)<\mu_\ell$ the model $N$ has an increasing
sequence of length $\theta$ if $\ell=1$ and a decreasing sequence 
of length $\theta$ if $\ell=2$.
\end{enumerate}
\mn
[Why?  Similar to the proof of $\circledast_1$ using it].
\mn
\begin{enumerate}
\item[$\circledast_3$]  for any $h_1 \in \cF^+_1,h_2 \in \cF^+_2$
  (i.e. $h_1\in {\cF}_1, h_2\in {\cF}_2$ satisfying 
$\boxdot^1_{h_1}+\boxdot^2_{h_2}$ from Definition \ref{b8}), there is
  a scattered $N \in K^*$ satisfying 
\sn
\begin{enumerate}
\item[$(a)$]  $(*)$ of clause (ix) of Definition \ref{b16} for $h_1$,
  that is, $(*)^1_{N,h_1}$
\sn
\item[$(b)$]  $(*)$ of clause (x) of Definition \ref{b16} for $h_2$,
  that is, $(*)^2_{N,h_2}$ 
\sn
\item[$(c)$]  $P^N_\alpha\neq \emptyset$ for $\alpha<\alpha^*$
\sn
\item[$(d)$]  if $\theta = \cf(\theta)<\mu_1$ then $N$ has an increasing 
sequence of length $\theta$
\sn
\item[$(e)$]  if $\theta = \cf(\theta)<\mu_2$, then $N$ has a
  decreasing sequence of length $\theta$.
\end{enumerate}
\end{enumerate}
\mn
[Why?  Let $N$ be as in $\circledast_2$. We define $M$ as follows:
$M$ has set of elements $\{(\ell,x):\ell\in\{-1,0,1\}$, and 
$\ell=-1 \Rightarrow x\in \Dom (h_2)$ and $\ell=0\Rightarrow x\in 
|N|$ and $\ell=1 \Rightarrow x \in \Dom(h_1)\}$ and 
$(\ell_1,x_1)<^M (\ell_2,x_2)$ \underline{if and only if} $\ell_1<\ell_2$ or 
$\ell_1 = \ell_2 = -1$ and $x_2<x_1$ (as ordinals) or $\ell_1=\ell_2=0$
and $x_1<^N x_2$ or $\ell_1=\ell_2=1$ and $x_1<x_2$ (as ordinals).

Lastly, $N = (M,P^N_\alpha)_{\alpha < \alpha(*)}$ where for $\alpha <
 \alpha^*$ we let $P^M_\alpha = \{(\ell,x) \in M:\ell = -1 \wedge
 h_2(x)=\alpha$ or $\ell=0 \wedge x \in P^N_\alpha$ or $\ell = 1
 \wedge h_1(x) = \alpha\}$].

At last we define by induction on $n<\omega$, $(M^n,P^n_\alpha)_{\alpha <
\alpha(*)}$ such that: 
\mn
\begin{enumerate}
\item[$(i)$]  $(M^n,P^n_\alpha)_{\alpha< \alpha(*)}\in K^*$ 
is a submodel of $(M^{n+1},P^{n+1}_\alpha)_{\alpha< \alpha(*)}$,
\sn
\item[$(ii)$]  $M^n$ is scattered, and so 
every interval contains a jump, i.e., an
empty open interval
\sn
\item[$(iii)$]  $(M^n, P^n_\alpha)_{\alpha< \alpha(*)}\in
  K^*_{h_1,h_2}$, 
\sn
\item[$(iv)$]  If $x \in P^n_\alpha$ has no 
immediate predecessor in $M^n$ (recalling $M_n$ has no first element), 
\then \, clause (ix) of 
Definition \ref{b16} holds for it, really follows by \ref{b24}(1)
\sn
\item[$(v)$]  If $x \in P^n_\alpha$ is neither last nor has an 
immediate successor (recalling $M_n$ has no last element), 
\then \, clause (x) holds for it, really follows
by \ref{b24}(1)
\sn
\item[$(vi)$]  If $x\in M^{n+1} \setminus M^n$, then for some $y,z \in M^n$: 
\[
y<x<z, \text{ and } \neg(\exists t\in M^n) y<t<z.
\]
\sn
\item[$(vii)$]  For every $y< z$ in $M^n$: in $M^{n+2}$ the 
element $y$ has no immediate successor and the element $z$ has no 
immediate predecessor, $\bigwedge\limits_{\alpha} P^{n+2}_\alpha \cap (y,
z)^{M^{n+2}} \ne \emptyset$, in $(y, z)^{M^{n+2}}$ there are
increasing sequences of any length $\theta=\cf(\theta)
< \mu_1$ in $(y, z)^{M^{n+2}}$ there are
decreasing sequences of any length $\theta=\cf(\theta)< \mu_2$. 
\end{enumerate}
\mn
There is no problem in this and $(\bigcup\limits_{n} M^n,\bigcup\limits_{n}
P^n_{\alpha})_{\alpha< \alpha(*)}$ is as required. 

That is, for $n=0$ use $\circledast_3$ for $(h_1,h_2)$. 
Given $(M^n,P^n_\alpha)_{\alpha<\alpha(*)}$, to get
$(M^{n+1},P^{n+1}_\alpha)_{\alpha<\alpha(*)}$, for each empty open
interval $(x,y)$ of $M^n$, we insert in this interval a copy of 
$N$ as constructed in $\circledast_3$ \underline{but} 
with $(g_2(\alpha_1), g_1(\alpha_2))$ 
here standing for $(h_1,h_2)$ there when $x \in P^n_{\alpha_1},y \in 
P^n_{\alpha_2}$.
\end{PROOF}

\begin{claim}
\label{b44}
If $h_1\in {\cF}^+_1$ and $h_2 \in {\cF}^+_2$, \then \, every two members 
of $K_{h_1, h_2}$ are isomorphic.
\end{claim}

\noindent
We shall prove this below. 
\begin{claim}
\label{b48}
1) If $N \in K^*$ and $(I_0, I_1)$ is a cut of $N$, i.e. of $M =
(|N|,<^N)$ (as a linear order, i.e. $M=I_0 \cup I_1,I_0\cap I_1=\emptyset$ and 
$t_0\in I_0\wedge t_1\in I \Rightarrow t_0 <^N t_1$), \then \, exactly
one of the following occurs: 
\mn
\begin{enumerate}
\item[$(i)$]  $I_0$ has a last element,
\sn
\item[$(ii)$]  $I_0$ is empty,
\sn
\item[$(iii)$]  $I_1$ has a first element,
\sn
\item[$(iv)$]  $I_1$ is empty,
\sn
\item[$(v)$]  $\cf(I_0) = \cf(I^*_1) = \aleph_0$.
\end{enumerate}
\mn
2)  If $N \in K$ then the set of cuts of case (v) above is
dense.

\noindent
3) If $N \in K$ and $I$ is an infinite subset of $N$ \then \,  we
can find $J$ such that:
\mn
\begin{enumerate}
\item[$(i)$]   $I\subseteq J \subseteq N$,
\sn
\item[$(ii)$]  $|J| = |I|$,
\sn
\item[$(iii)$]  $J$ has neither a first nor a last member,
\sn
\item[$(iv)$]  if $x \in N \setminus J$ and $N_{J,x} = N \rest
 A_{J,x},M_{I,x} = (A_{J,x},<^N \rest A_{J,x})$ where $A_{J,x}
 = \{y \in M: x,y$ realize the same cut of $J\}$ \then \,
\sn
\begin{enumerate}
\item[$(\alpha)$]  $N_{J,x}$ has no last element,
\sn
\item[$(\beta)$]  if $N_{J,x}$ is bounded in $M$ and $\cf(N_{J,x})>
\aleph_0$ \then \, it has a least upper bound in $J$,
\sn
\item[$(\gamma)$]   $N_{J,x}$ has no first element,
\sn
\item[$(\delta)$]   if $N_{J,x}$ is bounded from below in $N_{J,x}$ and
$\cf(N^*_{J, x})> \aleph_0$ \then \, it has a maximal lower bound in $J$.
\end{enumerate}
\sn
\item[$(v)$]   the number of members in $\{N_{J,x}:x\in N \setminus J\}$ is 
$\leq |J|+1$.
\end{enumerate}
\end{claim}

\begin{PROOF}{\ref{b48}}
Straightforward.
\end{PROOF}

\begin{definition}
\label{b52}
1) For a linear order $M$, if 
$J \subseteq M$ satisfies clauses (iii) + (iv) of claim 
\ref{b48}(3) \then \, we say that $J$ is quite closed in $M$. 

\noindent
2) Similarly for $J \subseteq N,N$ an expansion of a linear order.
\end{definition}
\bigskip

\begin{PROOF}{\ref{b44}}
\underline{Proof of Claim \ref{b44}}:

Let $h_1 \in {\cF}_1$, $h_2\in {\cF}_2$ and assume 
$N_1,N_2 \in K_{h_1,h_2}$, and we shall prove that $N_1,N_2$ are 
isomorphic.  Let $N_\ell=\bigcup\limits_{n<\omega} A_{\ell,n}$ 
with $M_\ell\restriction A_{\ell,n}$ being scattered, of course,
$M_\ell = N_\ell \rest \{<\}$. 
Let ${\cG}$ be the family of $f$ such that:
\mn
\begin{enumerate}
\item[$(a)$]  $f$ is a one-to-one function,
\sn
\item[$(b)$]  $\Dom(f)$ is a quite closed subset of $M_1$, see
  Definition \ref{b52}
\sn
\item[$(c)$]  $\Rang(f)$ is a quite closed subset of $N_2$, see
  Definition \ref{b52}
\sn
\item[$(d)$]  $f$ is an isomorphism from $N_1 \rest \Dom(f)$ onto
$N_2 \restriction \Rang(f)$
\sn
\item[$(e)$]  $M_1 \rest \Dom(f)$ is a scattered linear order.
\end{enumerate}
\mn
Now
\mn
\begin{enumerate}
\item[$\boxtimes_1$]  there is $f_1 \in {\cG}$ such that $\Dom(f_1)$ is an
unbounded subset of $N_1$, and $\Rang(f_1)$ is an unbounded subset of $N_2$. 
\end{enumerate}
\mn
[Why?   As $N_1,N_2 \in K_{h_1,h_2}$, using $h_1 \in \cF^+_1$, see
  \ref{b8}(3) and Definition \ref{b24}(3).]
\mn
\begin{enumerate}
\item[$\boxtimes_2$]  There is $f_2 \in {\cG}$ such that $\Dom(f_2)$ is
a subset of $N_1$ unbounded from below and $\Rang(f_2)$ is an unbounded
from below subset of $N_2$.
\end{enumerate}
\mn
[Why?   As $N_1,N_2 \in K_{h_1,h_2}$, using $h_2 \in \cF^+_1$, see
  \ref{b8}(3) and Definition \ref{b24}(3).]
\mn
\begin{enumerate}
\item[$\boxtimes_3$]  There is $f_0\in {\cF}$ satisfying the demands in
$\boxtimes_1$ and $\boxtimes_2$.
\end{enumerate}
\mn
[Why?  Let $f_1,f_2$ be from $\boxtimes_1,\boxtimes_2$, respectively. We
choose $x \in \Dom(f_1)$ and $y \in \Dom(f_2)$ such that 
$M_1 \models ``y<x"$ and $M_2\models ``f_2(y)<f_1(x)"$, note that 
this is possible by the choices of $f_1$ and $f_2$. 

Let 
\[
f_0 = (f_1\restriction\{z\in\Dom(f_1):x <^{N_1} z\})\cup (f_2
\rest \{z \in \Dom(f_2):z<^{N_2} y\}).
\]

\mn
Clearly $f_0$ is as required.]

For any $f \in {\cG}$ which extends $f_0$ and $t \in N_1
\setminus \Dom(f)$ we let

\[
N_{1,f,t} = N_1 \restriction\{s\in N_1:(\forall x\in\Dom(f))[x<_{M_1} t\equiv
x<_{M_1} s] \text{ and } s \notin \Dom(f)\}.
\]

\mn
Now
\mn
\begin{enumerate}
\item[$\boxtimes_4$]  if $f \in {\cG}$ extends $f_0$ and $t\in N_1 \setminus
\Dom(f)$ and $n<\omega$ and $A_{1,n} \cap N_{1,f,t} \ne \emptyset$ 
then there is $f'$ such that 
\sn
\begin{enumerate}
\item[$(i)$]  $f \subseteq f' \in \cG$,
\sn
\item[$(ii)$]  $\Dom(f') \setminus \Dom(f) \subseteq N_{1,f,t}$,
\sn
\item[$(iii)$]   if $s\in M_{1,f,t}\setminus\Dom(f')$ and $A_{1,n}
\cap N_{1,f,t}\neq \emptyset$ then 
\[
\Dep(N_{1,f',s} \rest A_{1,n}) < \Dep(N_{1,f,t} \rest A_{1,n}) =
\Dep(N_{1,f,s} \restriction A_{1,n}).
\]
\end{enumerate}
\end{enumerate}
\mn
[Why?  First note that there are 
$t_0<t_1$ in $\Dom(f_1)$ such that $N \models ``t_0 < t < t_1"$, this
holds by the choice of $f_0$ (recalling we are assuming $f \ge f_0$.
Second, we can demand that $N_{1,f,t} = 
N_1 \rest (t_0,t_1)_{N_1}$, just by the definition of ``$\Dom(f)$ is
quite closed" recalling the assumption on $f$. 

By Claim \ref{b36} it is enough to consider the following 
three cases.
\medskip

\noindent
\underline{Case 1}:   There is $s_1 \in N_{1,t}\cap A_{1,n}$ such that
$\Dep((N_{1,f,t})_{>s_1}) < \Dep(N_{1,f,t})$ (so possibly
$(N_{1,f,t})_{>s_1}$ is empty) and $\Dep((N_{1,f,t})_{<s_1}) < 
\Dep(N_{1,f,t})$. 

Let $s_1\in P^{N_1}_\alpha$ (clearly such $\alpha$ exists). 
Now, $(f(t_1),f(t_2))$ is an open interval of $N_2$ hence 
there is in it an $s_2\in P^{N_2}_\alpha$. Let
$f' = f\cup\{\langle s_1,s_2\rangle\}$.
\medskip

\noindent
\underline{Case 2}:   For every $x \in N_{1,f,t}$ we have 
$\Dep((N_{1,f,t})_{<x}) < \Dep(N_{1,f,t})$.

Let $\alpha_1$ be such that 
$t_1\in P_{\alpha_1}^{N_1}$, so also $f(t_1)\in P^{N_2}_{\alpha_1}$, and
imitate the proof of $\boxtimes_1$.
\medskip

\noindent
\underline{Case 3}:   For every $x\in M_{1,f,t}$ we have 
$\Dep((N_{1,f,t})_{>x}) < \Dep(N_{1,f,t})$.

Let $\alpha_0$ be such that 
$t_0\in P^{N_1}_{\alpha_0}$, so also $f(t_0) \in P^{N_2}_{\alpha_0}$ and
immitate the proof of $\boxtimes_2$.

So $\boxtimes_4$ holds indeed.]
\mn
\begin{enumerate}
\item[$\boxtimes_5$]  If $f \in {\cG}$ extends $f_0$ and $n<\omega$ 
\then \, there is $f'$ such that 
\sn
\begin{enumerate}
\item[$(i)$]   $f \subseteq f' \in {\cG}$,
\sn
\item[$(ii)$]   if $t \in N_1 \setminus \Dom(f')$ then 
$\Dep(N_{1,f',t} \rest A_{1,n})< \Dep(N_{1,f,t}\rest A_{1,n})$. 
\end{enumerate}
\end{enumerate}
\mn
[Why?  Let $\{t_\epsilon:\epsilon<\epsilon(*)\}$ be such that
$t_\epsilon\in N_1 \setminus \Dom(f)$ and $\langle N_{1,f,t_\epsilon}:
\epsilon<\epsilon(*)\rangle$ lists $\{N_{1,f,t}:t \in N_1
\setminus \Dom(f)$ and $A_{1,n}\cap N_{1,f,t}\neq \emptyset\}$ 
with no repetitions.  For each $\epsilon$ let
$f_\epsilon'$ be as in $\boxtimes_4$ for $t_\epsilon$, and let
$f'=\bigcup\limits_{\epsilon<\epsilon(*)} f_\epsilon'$. Now check, so
$\boxtimes_5$ holds indeed.]

For any $f \in {\cG}$ which extends $f_0$ and $t\in M_2\setminus\Rang(f)$, let 

\[
N_{2,f,t} = N_2 \restriction \{s\in M_2:(\forall x\in\Rang(f))(x<^{N_2}s
\Leftrightarrow x<^{N_2} t)\}.
\] 

\mn
Just as in $\boxtimes_4$, $\boxtimes_5$ we can show:
\mn
\begin{enumerate}
\item[$\boxtimes_6$]  if $f \in {\cG}$ extends $f_0$ and $n<\omega$ 
\then \, there is $f'$ such that
\sn
\begin{enumerate}
\item[$(i)$]   $f \subseteq f' \in {\cG}$,
\sn
\item[$(ii)$]  if $t \in N_2 \setminus \Rang(f)$ and 
$A_{2,n}\cap M_{2,f,t} \ne \emptyset$ then $\Dep(N_{2,f',t}
\rest A_{1,n}) < \Dep(M_{2,f,t}\rest A_{1,n})$. 
\end{enumerate}
\end{enumerate}
\mn
Lastly, we choose $f_n \in {\cF}$ by induction on $n<\omega$ such that $k<m
\Rightarrow f_k \subseteq f_m$. For $n=0$ we have already chosen $f_0$. If
$n=k^2+2 m<(k+1)^2$, let $f_{n+1}$ relate to $f_n$ as $f'$ relates to $f$ in
$\boxtimes_5$ (for $A_{1,m}$). If $n=k^2+2m+1<(k+1)^2$, 
let $f_{n+1}$ relate to $f_n$ as 
$f'$ relates $f$ in $\boxtimes_6$ (for $A_{2,m}$). 

Let $f = \bigcup\limits_{n\in\omega} f_n$, clearly $f$ is a partial isomorphism
from $N_1$ to $N_2$. Now, $\Dom(f) = N_1$, because if 
$t\in N_1\setminus\Dom(f)$
then for some $n$ we have $t \in A_{1,n}$ and clearly
$\langle \Dep(N_{1,f_m,t}\rest A_{1,n}):n<\omega\rangle$ is a non-increasing
sequence of ordinals (by \ref{b36}(2)), and for every $k>m$ we have
$\Dep(N_{1,f_{k^2+2m},t} \rest A_{1,n})< 
\Dep(N_{1,f_{k^2+2m+1},t}\rest A_{1,n})$ because of 
the use of $\boxtimes_5$.  A contradiction, so really 
$\Dom(f)=M_1$. Similarly $\Rang(f)=M_2$ and we are done. 
\end{PROOF}

\begin{definition}
\label{b56}
We say $N \in K^*$ is almost $\kappa$-homogeneous \when \,:
\mn
\begin{enumerate}
\item[$\bullet$]  if $I \subseteq N$, $|I|<\kappa$ {\it then} we can find $J$, 
$I\subseteq J\subseteq N$, $|J|< \kappa$ such that 
\sn
\begin{enumerate}
\item[$(*)$]  if $s$, $t \in (N \setminus J)$ realize the same
 cut of $J$ and $s\in P^N_\alpha \Leftrightarrow t \in P^N_\alpha$
for every $\alpha<\alpha(*)$, \then \, there is an automorphism of
$N$ over $J$ mapping $s$ to $t$. 
\end{enumerate}
\end{enumerate}
\mn
Similarly to the proof of \ref{b44}.
\end{definition}

\begin{conclusion}
\label{b60}
Assume $h_1 \in \Rang(g_1),h_2 \in \Rang(g_2)$.

\noindent
1)  If $N \in K^{\hom}_{h_1,h_2}$, $n<\omega$ and $x_1<x_2<\ldots <x_n$ in
$N$, and $y_1 < \ldots < y_n$ in $N$, and $x_m\in P^N_\alpha \Leftrightarrow
y_m \in P^N_\alpha$ for $\alpha<\alpha(*)$, $m\in \{1,\ldots,n\}$,
\then \, there is an automorphism of $N$ mapping $x_m$ to $y_m$ for
$m=1,\ldots,n$. 

\noindent
2) If $N \in K_{h_1,h_2}$ and $J \subseteq N$ is quite closed in $M$ \then   
\mn
\begin{enumerate}
\item[$(*)$]  if $s$, $t\in N \setminus J$ realize the same cut of
$J$ and $s\in P^N_\alpha \Leftrightarrow t\in P^N_\alpha$ for $\alpha<
\alpha(*)$, \then \, there is an automorphism of $N$ over $J$ 
mapping $s$ to $t$. 
\end{enumerate}
\mn
3)  Every $N \in K^{\hom}$ is almost $\kappa$-homogeneous (where $\kappa \ge 
\aleph_0$).

\noindent
4)  Assume $N \in K^{\hom}_{h_1,h_2}$ and $J_1$, $J_2 \subseteq N$ are quite
closed and $[J_1$ is unbounded in $N$ iff $J_2$ is unbounded in $N]$ and 
$[J_1$ is unbounded in $N^*$ iff $J_2$ is unbounded in $N^*]$. If $f$ is an
isomorphism from $N \restriction J_1$ onto $N \restriction J_2$
\then \, we can extend $f$ to an automorphism of $M$.
\end{conclusion}

\begin{PROOF}{\ref{b60}}
Should be clear.
\end{PROOF}
\bigskip

\subsection {Examples}\
\bigskip

\noindent
In this subsection we consider some examples.
\begin{content}
\label{b63} 
We do not assume \ref{b8} fully, still $\lambda,\mu_1,\mu_2$ are as in 
\ref{b8}(a) and  $\theta$
will denote a regular cardinal $< \mu_1 \cap \mu_2$, usually uncountable.
\end{content}

\begin{definition}
\label{b64}
Assume $\theta=\cf(\theta) < \Min\{\mu_1,\mu_2\}$ and let $\bar\sigma =
\langle (\sigma^1_\alpha,\sigma^2_\alpha):\alpha < \beta(*)\rangle$ 
list (with no repetitions) the pairs $(\sigma_1,\sigma_2)$ of 
(infinite) regular cardinals such that $\sigma_\ell < \mu_\ell,
\theta \in \{\sigma_1,\sigma_2\}$ and 
$\alpha=0\Rightarrow (\sigma^1_\alpha,\sigma^2_\alpha)= (\theta,\theta)$.

\noindent
1)  We let $g^\theta_\ell = g^{\theta,\bar\sigma}_\ell \in {\cF}_\ell$
be defined by: for $\beta < \beta(*),g^\theta_\ell(\beta)$ is a
function with domain $\sigma^\ell_\beta$ and for $\gamma<\sigma^\ell_\beta$
\mn
\begin{enumerate}
\item[$(*)$]  $g^\theta_\ell(\beta)(\gamma)=\alpha$ iff 
$\alpha<\alpha(*)$ satisfies
\sn
\begin{enumerate}
\item[$(a)$]  if $\gamma<\sigma$ is a limit ordinal
then $\sigma^\ell_\alpha=\cf(\gamma),\sigma^{3-\ell}_\alpha =\theta$
\sn
\item[$(b)$]  if $\gamma<\sigma$ is non-limit then $\alpha=0$.
\end{enumerate}
\end{enumerate}
\mn
1A) For $\ell \in \{1,2\}$ and regular $\theta < \mu_\ell$ let
$h^\ell_\theta$ be the unique $h:\theta \rightarrow \beta(*)$ such
that $\sigma^\ell_\alpha = \theta \Rightarrow g^\theta_\ell(\alpha) =h$.

\noindent
2) Let $\bold c= \bold
c^{\can}_{\lambda,\mu_1,\mu_2,\theta,\bar\sigma}$ be the unique $\bold
c$ such that $(\lambda^{\bold c},\mu^{\bold c}_1,\mu^{\bold
  c}_2,g^{\bold c}_1,g^{\bold c}_2) = (\lambda,\mu_1,\mu_2,
g^\theta_1,g^\theta_2)$, see \ref{b66}(2).

\noindent
3) In (2) we may omit $\bar \sigma$ when $\alpha < \beta(*) \Rightarrow
\alpha = \otp(u_\alpha,<_{\lex})$ where $u_\alpha = \{(\sigma_1,\sigma_2):
\sigma_1 = \cf(\sigma_1) < \mu_1,\sigma_2 = \cf(\sigma_2) 
< \mu_2$ and $(\sigma_1,\sigma_2) <_{\lex} (\sigma^1_\alpha,
\sigma^2_\alpha)\}$; justify by \ref{b66}(1),(2).

\noindent
4) For regular $\theta_1 < \mu_1,\theta_2 < \mu_2$ we let 
$K^{\can}_{\theta_1,\theta_2} =
K^{\hom}_{h^1_{\theta_1},h^2_{\theta_2}}(\bold c)$ where $\bold c$ is
from part (2),(3) and $h^\ell_{\theta_1}$ is from part (1A). 
For $u \subseteq \beta(*)$ non-empty let
$K^{\can}_{\theta_1,\theta_2,u} = \{(|N|,<^N) \rest \bigcup\limits_{\alpha \in
  u} P^N_\alpha:N \in K^{\can}_{\theta_1,\theta_2}\}$.  Can define for
the general case.
\end{definition}

\begin{claim}
\label{b66}
Let $\lambda,\mu_1,\mu_2$ be as in \ref{b8}(a) and $\theta$ be regular
$< \min\{\mu_1,\mu_2\}$.

\noindent
1) There is $\bar\sigma = \langle (\sigma^1_\alpha,
\sigma^2_\alpha):\alpha < \beta(*)\rangle$ as in Definition \ref{b64}
so $|\beta(*)| = |\Reg \cap \mu_1| \times |\Reg \cap \mu_2|$.  
Moreover, there is one and only one $\bar\sigma$ as in \ref{b64}(3).

\noindent
2) For $\bar\sigma$ as in part (1), $\bold
c^{\can}_{\lambda,\mu_1,\mu_2,\theta,\bar\sigma}$ is well defined,
that is, there is a unique context 
$\bold d$, recalling Definition \ref{b8}, such that:
\mn
\begin{enumerate}
\item[$(a)$]  $(\lambda^{\bold d},\mu^{\bold d}_1,\mu^{\bold d}_2) =
  (\lambda,\mu_1,\mu_2)$
\sn
\item[$(b)$]  $\alpha^{\bold d}_* = \beta(*)$
\sn
\item[$(c)$]  $(g^{\bold d}_1,g^{\bold d}_2)$ is as in \ref{b64}(1).
\end{enumerate}
\mn
3) If for $\iota=1,2$ we have $\bar\sigma_\iota = \langle
(\sigma^1_{\iota,\alpha},\sigma^2_{\iota,\alpha}):\alpha <
\beta(\iota)\rangle$ as in part (1), i.e. as in \ref{b64} and 
$\bold d_\iota$ is as in part (2) for $\bar\sigma_\iota$ 
and $h_{\iota,1} \in \cF^{+,\bold d_\iota}_1,h_{\iota,2} \in \cF^{+,\bold d_i}_2$ and $N_\iota =
(M_\iota,P^\iota_\alpha)_{\alpha < \beta(\iota)} \in K^{\bold
  d_\iota}_{h_{\iota,1},h_{\iota,2}}$ \then \,
\mn
\begin{enumerate}
\item[$(a)$]  there is a unique $f:\beta(1) \rightarrow \beta(2)$ such
 that $(\sigma^1_{1,\alpha},\sigma^2_{1,\alpha}) =
 (\sigma^1_{2,f(\alpha)},\sigma^2_{2,f(\alpha)})$ for $\alpha <
 \beta(1)$; moreover $f$ is one to one onto
\sn
\item[$(b)$]  $M_1,M_2$ are isomorphic linear orders
\sn
\item[$(c)$]  moreover, there is an isomorphism $\bold f$ from $N_1$
  onto $N_2$ which maps $P^1_\alpha$ onto $P^2_{f(\alpha)}$ for every
  $\alpha < \beta(1)$.
\end{enumerate}
\mn
4) For $\bold c$ as in \ref{b64}(2) so
$\bar\sigma = \bar\sigma^{\bold c},\theta = \cf(\theta) >
\aleph_0$ as in \ref{b64} letting $(\beta(*) = \ell g(\bar\sigma))$
if $N = (M,P_\alpha)_{\alpha < \beta(*)}
\in K^{\bold c}$, so $M$ a linear order, then $N$ is uniquely
determined by $M$, i.e. $P^N_\alpha = \{d \in M:M_{<d}\}$ has
cofinality $\sigma^1_\alpha$ and $M_{>d}$ has co-initiality $\sigma^2_\alpha$.
\end{claim}

\begin{PROOF}{\ref{b66}}
Should be clear.
\end{PROOF}

\noindent
The following is used in \cite{Sh:1019}.
\begin{claim}
\label{b69}
Assume $\bold c = \bold c^{\can}_{\lambda,\mu_1,\mu_2,\theta}$, see
\ref{b64}(2), \ref{b66}(2).

\noindent
1) If $N_1,N_2 \in K^{\hom}$ \then \,: $N_1,N_2$ are isomorphic \Iff
   \, $N_1,N_2$ has the same cofinality and same co-initiality.

\noindent
2) So $K^{\hom} = \cup\{K^{\hom}_{\theta_1,\theta_2}:\theta_1 <
   \mu_1,\theta_2 < \mu_2$ are regular$\}$.

\noindent
3) Assume $\alpha \in (0,\lambda^+)$ is a successor and $N =
\sum\limits_{\beta \le \alpha} N_\beta$.

A sufficient condition for $N \in K^\oplus_{\theta,\theta}$ is:
\mn
\begin{enumerate}
\item[$(a)$]  each $N_\beta$ is from $K_{\theta,\theta}$ or is a
  singleton
\sn
\item[$(b)$] if $\theta > \aleph_0$ and $\beta < \alpha$ \then \,
  $N_\beta$ is a singleton or $N_{\beta +1}$ is a singleton but not both
\sn
\item[$(c)$]   $N_0 \in K^\oplus_{\theta,\theta}$ and $N_\alpha$ is a singleton
\sn
\item[$(e)$]  if $\delta < \alpha$ is a limit ordinal then $N_\delta$ is a
singleton and $P^{N_\delta}_\alpha = |N_\delta|$ when
$(\sigma^1_\alpha,\sigma^2_\alpha) = (\cf(\delta),\theta)$.
\end{enumerate}
\mn
4) Like (3) for an inverse, well-ordered sum except that in (e) we
deduce $(\sigma^1_\alpha,\sigma^2_\alpha) = (\theta,\cf(\delta))$.
\end{claim}

\begin{PROOF}{\ref{b69}}
Easy.
\end{PROOF}

\begin{definition}
\label{b71}
For $\lambda,\mu_1,\mu_2$ as in \ref{b8}(a) and $\theta = \cf(\theta)
< \Min\{\mu_1,\mu_2\},\bold c = \bold c^{\can}_{\lambda,\mu_1,\mu_2,\theta}$ be
as above.  Let $K^{1-\can}_{\lambda,\mu_1,\mu_2,\theta}$ be
$K^{\can}_{\theta_1,\theta_2,u}$ for $\bold c$ recalling \ref{b64}(4) 
where $u = \{\alpha\}$ with $\alpha$ such that $(\theta_1,\theta_2) =
(\aleph_0,\aleph_0)$. 

Let $K^{2-\can}_{\lambda,\mu_1,\mu_2,\theta}$ be
$K^{\can}_{\lambda,\mu_1,\mu_2,\theta}$ for $\bold i$ where $n=\alpha_*(c)$.
\end{definition}

\begin{claim}
\label{b73}
For $\lambda,\mu_1,\mu_2,\theta$ and $\bold c$ as above.

\noindent
1) There is an $M \in K^{1-\can}_{\lambda,\mu_1,\mu_2,\theta}$ unique up
to isomorphism, it is a dense linear order of cardinality $\lambda$
   with cofinality and co-initiality $\theta$.

\noindent
2) $K^{1-\can}_{\lambda,\mu_1,\mu_2,\theta}$ is closed under well
   ordered sums of length $\alpha +1 < \lambda^+$.

\noindent
3) $K^{1-\can}_{\lambda,\mu_1,\mu_2,\theta}$ is closed under anti-well
   ordered sums of length $\alpha +1 < \lambda^+$.

\noindent
4) If $\mu \ge \theta$ and $\mu \le \mu_1,\mu \le \mu_2$ and $M \in
K^{1-\can}_{\lambda,\mu_1,\mu_2,\theta}$, \then \, for some algebra
   $\gB$ on $|N|$ with $\mu$ functions if $\gB' \subseteq \gB$ has
   cardinality $\mu$ then $N' \rest \gB' \in
 K^{1-\can}_{\mu,\mu^+,\mu^+,\theta}$.

\noindent
5) $M$ also satisfies the conclusion of \ref{b92} and \ref{b60}.
   E.g. Let $N_i \in K^{\can}_{\lambda,\mu_1,\mu_2,\theta}$ for $i \le
   \alpha$ and $N = \sum\limits_{i \le \alpha} N_i$.  So for each $i$
   there is $N^+_i \in K^{\hom}_{\aleph_0,\aleph_0}$.  We can find
   $\langle N^{++}_i:i \le \beta\rangle$ and increasing $g:\alpha +1
   \rightarrow \beta +1$ such that:
\mn
\begin{enumerate}
\item[$(a)$]   if $N^+_i = N^{++}_{g(i)}$
\sn
\item[$(b)$]  $g$ is increasing continuous and $g(i+1) =
  g(i)+2,g(0)=0,g(\delta) = \delta$ and $h(\alpha) = \beta$
\sn
\item[$(c)$]  if $j \in \beta(*) \backslash \Rang(g)$ then  $N^{++}_j$
  is a singleton
\sn
\item[$(d)$] $\langle N^{++}_i:i \le \beta\rangle$ is as in
  \ref{b69}(3).
\end{enumerate}
\end{claim}

\begin{PROOF}{\ref{b73}}
Should be clear (check almost homogeneous), e.g.

\noindent
4) Let $N_i \in K^{\can}_{\lambda,\mu_1,\mu_2,\theta}$ for $i \le
   \alpha$ and $N = \sum\limits_{i \le \alpha} N_i$.  So for each $i$
   there is $N^+_i \in K^{\hom}_{\aleph_0,\alpha_0}$.  We can find
   $\langle N^{++}_i:i \le \beta\rangle$ and increasing $g:\alpha = 1
   \rightarrow \beta +2$ such that
\mn
\begin{enumerate}
\item[$(a)$]  if $N^+_i = N^{\tau^+}_{g(i)}$
\sn
\item[$(b)$]  $g$ is increasing continuous and $g(i+2) =
g(i)+2,g(0)=0,g(\delta)=\delta$
\sn
\item[$(c)$]  if $j \in j(*) \backslash \rang(g)$ then $N^{++}_j$ is a
singleton
\sn
\item[$(d)$]  $\langle N^{++}_i:i \le \beta\rangle$ is as in
\ref{b69}(3).
\end{enumerate}
\mn
Lastly, we apply \ref{b69}(3).
\end{PROOF}

\begin{claim}
\label{b75}
Assume $\lambda > \mu = \cf(\mu) > \theta,\bold c= \bold
c^{\can}_{\lambda,\lambda^+,\mu,\theta}$ so $(\mu_1,\mu_2) =
(\lambda^+,\mu)$ and $N \in K^{\hom}_{\lambda,\lambda^+,\mu,\theta}$
and $M = (N,<^N) \in K^{2-\can}_{\lambda,\lambda^+,\mu,\theta}$ and
let $\sigma = \cf(\sigma) \in [\mu,\lambda)$ and $T =
\Th_{\bbL_{\sigma,\sigma}}(M)$.  For a model $M'$ of $T$ let $N' =
    (M',\dotsc,P^{M'}_\alpha,\ldots)$ be defined as in \ref{b66}(4).

\noindent
1) $N' \in K^{2-\can}_{\lambda_1,\lambda^+_1,\mu,\theta}$ when:
\mn
\begin{enumerate}
\item[$(a)$]  $M'$ is a model of $T$ of cardinality $\lambda_1 \ge
  \sigma$
\sn
\item[$(b)$]  $N'$ is 2-homogeneous (i.e. if $M' \models ``s_1 < t_2
  \wedge s_2 < t_2"$ and $s_1 \in P^{N'}_\alpha \Leftrightarrow s_2
  \in P^{N'}_\alpha,t_1 \in P^{N'}_\alpha \Leftrightarrow t_2 \in
  P^{N'}_\alpha$) for $\alpha < \alpha_*(\bold c)$ then there is an
  automorphism of $N'$ (equivalently $M'$) mapping $(s_1,t_1)$ to
  $(s_2,t_2)$
\sn
\item[$(c)$]  $M'$ is the countable union of scattered sets
\sn
\item[$(d)$]  $(\alpha) \quad$ if 
$\kappa = \cf(\kappa) > \aleph_0$ if $\langle b_i:i <
  \kappa\rangle$ is an increasing bounded sequence in 

\hskip25pt $M$ then for
  some club $E$ or $\kappa$, for every $\delta \in E \cup \{\kappa\}$,

\hskip25pt $\bar b \rest \kappa$ has a $<^N$-lub
\sn
\item[${{}}$]  $(\beta) \quad$ similarly for $M^*$, the inverse of $M$.
\end{enumerate}
\mn
2) There is a first order sentence $\psi \in \bbL(\{<,F\}),F$ a
three-place function symbol such that $\{\{<\}:N$ a model of $T \ cup
\{\psi\}\}$ is equal to
$\cup\{K^{2-\can}_{\lambda^+_1,\lambda_1,\mu}:\lambda_1 \ge \sigma\}$.

\noindent
3) In part (1), if $\sigma$ is a compact cardinal we can omit clauses
   (c),(d).

\noindent
4) If $\sigma$ is a compact cardinal, \then\, the class from part (2)
   is categorical in every $\lambda_1 \ge \sigma$.
\end{claim}

\begin{PROOF}{\ref{b75}}
Should be clear.
\end{PROOF}
\bigskip

\centerline {$* \qquad * \qquad *$}
\bigskip

We now make the connection to \cite[\S3]{Sh:E59}.

\noindent
We may weaken a little the definition of weakly $\kappa$--skeleton like
(Definition \cite[3.1(1)=L3.1(1)]{Sh:E59}).
\begin{claim}
\label{b80}
Assume $\lambda > \kappa = \cf(\kappa)$, and $\bold d_\ell = 
\inv^\alpha_\kappa(I_\ell)$ for $\ell=1,2$ 
(see Definition \cite[3.4=L3.2]{Sh:E59}), $I_\ell$ a linear order of
cardinality $\le \lambda$, $\alpha < \lambda^+$ (for $\ell= 1,2$). 
\Then \, there are $\alpha^*$, $\mu_1,\mu_2$ (hence
${\cF}_1,\cF_2,\cF^+_1,\cF^+_2),g_1,g_2$ as in Context \ref{b8} such that: 
\mn
\begin{enumerate}
\item[$(*)_1$]  if $M \in K_{g_1(0),g_2(0)}$ and $u\subseteq \alpha^*$ is
non-empty then 
\sn
\begin{enumerate}
\item[$(a)$]  $\inv^\alpha_\kappa(\bigcup\limits_{\epsilon\in u}
P^M_\epsilon,<^M) = \bold d_1$
\sn
\item[$(b)$]  $\inv^\alpha_\kappa(\bigcup\limits_{\epsilon\in u}
P^M_{\epsilon},<^{M^*}) = \bold d_2$, recalling $M^*$ is the $M$ inverted
\end{enumerate}
\sn
\item[$(*)_2$]  if $\bold d_1 = \bold d_2$ and $K' = \{(P^M_0,<^M): 
M \in K_{g_1,g_2}\}$ then 
\sn
\begin{enumerate}
\item[$(a)$]   $K'$ is closed under sums of order type $\alpha$ and
$\alpha^*$ for $\alpha < \lambda^+$
\sn
\item[$(b)$]   each member of $K'$ is cardinality $\lambda$,
\sn
\item[$(c)$]  $K'$ is almost $\theta$-homogeneous for every $\theta$.
\end{enumerate}
\end{enumerate}
\end{claim}

\begin{PROOF}{\ref{b80}}
Straightforward.
\end{PROOF}

\noindent
Also in the cases we use 
skeletons from $K^\omega_{\tr}$ we may like to realize
distinct invariants rather than just non-isomorphic models.
\begin{definition}
\label{b82}
1) Let $N$ be a model of cardinality $\lambda$ with $|\tau_N| <
\lambda$ we say $\bar N$ is a $\lambda$-representation (of $N$), or
$\lambda$-filtration (of $M$) when:
\mn
\begin{enumerate}
\item[$(a)$]   $\bar N = \langle N_\alpha:\alpha < \lambda\rangle$
\sn
\item[$(b)$]  $N_\alpha \subseteq N$ has cardinality $< \lambda$
\sn
\item[$(c)$]  $\bar N$ is $\subseteq$-increasing continuous
\sn
\item[$(d)$]  $N = \cup\{N_\alpha:\alpha < \lambda\}$.
\end{enumerate}
\mn
2) For a $\lambda$-representation $\bar N$ let (on splitting see below)

\begin{equation*}
\begin{array}{clcr}
\Sep(\bar N) = \{\delta< \lambda:\, &\delta \text{ is limit, and for some }
\bar a \in \bigcup\limits_{\alpha< \lambda} N_\alpha\\
 &\text{ for every } \beta< \delta,\tp(\bar a,N_\delta,N) 
\text{ splits over } M_\beta\}.
\end{array}
\end{equation*}

\mn
3) $\Sep_{\Delta_1,\Delta_2}(\bar N) = \{\delta< \lambda:
\delta$ limit, and for some $\bar a\in \bigcup\limits_{\alpha<\lambda} 
N_\alpha$ for every $\beta< \delta$ the type $\tp_{\Delta_1}
(\bar a,N_\delta,N)$ does ($\Delta_1, \Delta_2$)-splits over
$M_\beta\}$. 

\noindent
4) Let $\Sep(N)$ be $\Sep(N)/\cD_\lambda$ for every representation
   of $M$.  Similarly $\Sep_{\Delta_1,\Delta_2}(N)$; both are
   justified because
\mn
\begin{enumerate}
\item[$\boxplus$]  $\Sep$ is $\check{\cD}_\lambda$-invariant of $N$, i.e. 
if $\bar N',\bar N''$ are $\lambda$-representations of $N$; 
$\|N\|=\lambda$ then $\Sep(\bar N')/ \check{\cD}_\lambda = 
\Sep(\bar N'')/ \check{\cD}_\lambda$ and
$\Sep_{\Delta_1,\Delta_2}(\bar N')/\check{\cD}_\lambda =
\Sep_{\Delta_1,\Delta_2}(\bar N'')/\check{\cD}_\lambda$ (when $\lambda=
\cf(\lambda)>\aleph_0$).
\end{enumerate}
\mn
5)  We say that $\tp_{\Delta_1}(a,B,N)$ does
$(\Delta_1,\Delta_2)$-split over $A \subseteq N$ (where $\bar a \in M,
B\subseteq N$) if for some $\bar b_1,\bar b_2 \in B,
\tp_{\Delta_2}(\bar b_1,A,N) =
\tp_{\Delta_2}(\bar b_2,A,N)$ but $\tp_{\Delta_1}(\bar a \char 94 \bar
b_1,A,N) \ne \tp_{\Delta_1}(\bar a \char 94 \bar b_2,A,N)$. 

\noindent
6)  If $\Delta_1 = \Delta_2$ is ${\bbL}_{\omega,\omega}(\tau(M))$, we may
omit $(\Delta_1,\Delta_2)$.

\noindent
7) We can replace $\check{\cD}_\lambda$ by appropriate ${\cE}$
giving us an $\omega$--sequence of sets (or an appropriate filters on
the set). 
\end{definition}

\begin{definition}
\label{b84}
1)  $N$ is $(\lambda,\Delta_1,\Delta_2)$-nice if
$\Sep_{\Delta_1,\Delta_2}(N) = \emptyset / D_\lambda$. 

\noindent
2) $N$ is $(<\lambda,\Delta)$-stable if for every $A\subseteq |N|$
of power $< \lambda$ 

\[
\lambda> |\{\tp_{\Delta}(\bar a,A,N):\bar a \in |M|\}|.
\]

\mn
3) $I \in K^\omega_\tr$ is {\it locally} $(\lambda,\bs,\bs)$-nice 
[{\it locally} $(<\lambda,\bs)$-stable] \If \, for every $\eta\in I
\setminus P^I_\omega$ the linear order $(\Suc_I(\eta),<)$ is 
$(\lambda,\bs,\bs)$-nice [$(< \lambda,\bs)$-stable].
\end{definition}

\begin{claim}
\label{b88}
Every $M \in K$ is $(\lambda,\bs,\bs)$-nice and $(<\lambda,bs)$-stable.
\end{claim}

\begin{PROOF}{\ref{b88}}
Easy (and as in \cite[\S6]{Sh:110}, mainly ``crucial fact'' of pg.217 there).
\end{PROOF}

\begin{claim}
\label{b92}
If $(A,<,P_\alpha)_{\alpha< \alpha(*)}\in K,S\subseteq \lambda$, and 

\[
M = (\bigcup\limits_{\alpha\in S} P_\alpha,< \rest 
(\bigcup\limits_{\alpha\in S} P_\alpha), P_\alpha)_{\alpha\in S},
\]

\mn
\then \, $M$ is $(<\lambda,\bs)$-stable and $(\lambda,\bs,\bs)$-nice.
\end{claim}

\begin{PROOF}{\ref{b92}}
Check.
\end{PROOF}
\bigskip

\subsection {Very Homogenous Linear Orders Revisited}\
\bigskip

We here start to indicate how we can generalize \S(2A).  
The case $\kappa = \aleph_0$ is the one in \S(2A).
\begin{definition}
\label{p2}
1) We say $\bold c$ is a context or $(\lambda,\kappa)$-context \when \, it
 consists of (so $\lambda = \lambda_{\bold c}$, etc.)
\mn
\begin{enumerate}
\item[$(a)$]  $\lambda = \lambda^{< \kappa} \ge \kappa = \cf(\kappa)$
\sn
\item[$(b)$]  $\alpha_* < \lambda^+,u_1 \subseteq \alpha_*,u_2 \subseteq
  \alpha_*,u_1 \cup u_2 \ne \alpha_*$ (or just $u_1 \cap u_2 \ne
  \alpha_*$), note here 1,2 stands for right, left
\sn
\item[$(c)$]  vocabulary $\tau = \{<\} \cup \{P_\alpha:\alpha <
  \alpha_*\}$, where $<$ is a binary predicate, $P_\alpha$ is a unary
  predicate
\sn
\item[$(d)$]  $K_{\all}$ is the class of $N$ such that (all stands for
  all)
\sn
\begin{enumerate}
\item[$(\alpha)$]  $N$ is a $\tau$-model
\sn
\item[$(\beta)$]  $<^N$ a linear order
\sn
\item[$(\gamma)$]  $\langle P^N_\alpha:\alpha < \alpha_*\rangle$ a partition
  of $N$
\sn
\item[$(\delta)$]  if $\partial = \cf(\partial) < \kappa$ and $\bar a
= \langle a_i:i < \partial\rangle$ is increasing/decreasing \then \, it has a
  $<_N$-lub/$<_N$-mdb; moreover if $\partial > \aleph_0$ then for a
club of $\delta < \partial$ this holds for $\bar a \rest \delta$, too
\end{enumerate}
\sn
\item[$(e)$]  $g_\ell$ is a $\ell$-nice function from $u_\ell$ into
  $\cF^*_\ell$, see below
\sn
\item[$(f)$]  $K_{\nic} \subseteq K_{\all}$ is defined in part (5) 
below ($\nic$ stands for nice)
\sn
\item[$(g)$]  $K_{\bas} \subseteq K_{\nic}$ has cardinality $\le
  \lambda$ and each $N \in K_{\bas}$ has cardinality $\le \lambda$ and
  some $N \in K_{\bas}$ has cardinality $\lambda$ ($\bas$ stands for
  basic, the generators).
\end{enumerate}
\mn
2) $\cF_\ell = \cF^\ell_{\bold c}$ is the set of function $f$ with domain a
regular $\partial \le \lambda$ into $\alpha_*$ such that any limit
$\delta < \partial,f(\delta) \in u_\ell$.

\noindent
3) $g_\ell:\alpha_* \rightarrow \cF_\ell$ is $\ell$-nice when
\mn
\begin{enumerate}
\item[$(a)$]  for every $\alpha < \alpha_*,h := g(\alpha)$ is
  $(g_\ell,\ell)$-nice, see below
\sn
\item[$(b)$]  if $\partial = \Dom(g_\ell(\alpha_1)) < \kappa,
g_\ell(\alpha) = g_2(\beta)$ then $\alpha = \beta$
\sn
\item[$(c)$]  if $h:\partial \rightarrow \alpha_*$ is
  $(g_\ell,\ell)$-nice and $\partial < \kappa$ then $h \in \Rang(g)$.
\end{enumerate}
\mn
4) $h \in \cF_\ell$ is $(g,\ell)$-nice when: if $\partial = \Dom(h)$
is regular \then \,
\mn
\begin{enumerate}
\item[$\boxdot^\ell_h$]   $h \in \cF_\ell$ and if 
$\delta < \Dom(h)$ is a limit 
ordinal of uncountable cofinality and 
$\beta=h(\delta)$ and $\langle \epsilon_i:i<\cf(\delta)\rangle$ 
is an increasing continuous sequence with limit $\delta$ then 
$\{i<\cf (\delta):(h(\delta))(\epsilon_i) = (g(\beta))(i)\}$ 
contains a club of $\cf(\delta)$.
For notational simplicity assume $\alpha^* \le \lambda$.
\end{enumerate}
\mn
5) $K_{\nic}$ is the class of $N$ such that
\mn
\begin{enumerate}
\item[$(a)$]  $N \in K_{\all}$ has cardinality $\le \lambda$
\sn
\item[$(b)$]  if $a \in P^N_\alpha$ and $\alpha \in u_1$ and $a$ has
  no immediate predecessor in $N$, \then \, there is an increasing
  sequence $\langle b_\alpha:\alpha \in \Dom(g_1(\alpha))\rangle$ in
  $N$ such that
\sn
\begin{enumerate}
\item[$(\alpha)$]  if $\alpha$ is a limit ordinal then $b_\alpha$ is
  the $<^*_N$-lub of $\langle b_\beta:\beta < \alpha\rangle$
\sn
\item[$(\beta)$]  if $\Dom(g_1(\alpha))$ is uncountable then $\{\alpha
  < \Dom(g_1(\alpha)):b_\alpha \in P^N_{g_1(\alpha)}\}$ contains a
  club of $\Dom g_1(\alpha)$
\sn
\item[$(\gamma)$]  if $\alpha$ is non-limit \then \, $g_1(\alpha)
  \notin u_1 \cup u_2$
\end{enumerate}
\sn
\item[$(c)$]  like (b), replacing $u_1,g_1$ increasing, predecessor, lub
  by $u_2,g_2$ decreasing, successor, $\glb$.
\end{enumerate}
\end{definition}

\begin{convention}
\label{p6}
In this sub-section, $\bold c$ will be a fixed context, if not said otherwise.
\end{convention}

\begin{definition}
\label{p9}
We define a two-place relation $\le_{\nic}$ on $K_{\nic},N_1
\le_{\nic} N_2$ \Iff \,
\sn
\begin{enumerate}
\item[$(a)$]  $N_1,N_2 \in K_{\nic}$
\sn
\item[$(b)$]  $N_1 \subseteq N_2$
\sn
\item[$(c)$]  if $a \in P^{N_1}_\alpha$ and $\alpha \in u_{\bold c,1}$
and $a$ has no immediate predecessor in $N_1$, then $(N_1)_{<a}$ is
unbounded in $(N_2)_{<a}$ from above
\sn
\item[$(d)$]  if $a \in P^{N_1}_\alpha,\alpha \in u_{\bold c,2}$
and $a$ has no immediate successor in $N_1$, \then \, $(N_1)_{<a}$ is
unbounded in $(N_2)_{>a}$ from below
\sn
\item[$(e)$]  if $a \in N_2 \backslash N_1$, \then \, $(N_2)_{<a} \cap
  N_1$ has a last element or $(N_2)_{>a} \cap N_1$ has a first element.
\end{enumerate}
\end{definition}

\begin{claim}
\label{p12}
$(K_{\nic},\le_{\nic})$ is a partial order preserved under isomorphisms.
\end{claim}
\newpage

\section {On pcf and other uncountable combinatorics}

In this section we define and quote but do not prove.
\begin{definition}
\label{wd.1}
1)  For $\lambda$ regular uncountable we define the weak diamond ideal, 
$\check I^{\wed}_\lambda = \check I^{\wed}[\lambda]$ as the 
family of small subset of $\lambda$, where:

\noindent
2) We say $S\subseteq \lambda$ is small \If \, it is 
${\bold  F}$-small for some colouring function $\bold F$ from 
${}^{\lambda>}\lambda$ to $\theta$ where

\noindent
3)  We say $S \subseteq \lambda$ is ${\bold F}$ small \If \, 
(${\bold F}$ is as above and) 
\mn
\begin{enumerate}
\item[$(*)_S$]   for every $\bar c \in {}^{S} 2$ for some $\eta\in 
{}^{\lambda}\lambda$ the set $\{\delta\in S:\bold F 
(\eta \rest \delta) = c_\delta\}$ is not stationary.
\end{enumerate}
\end{definition}

\begin{claim}
\label{wd.2} 
If $\lambda=\mu^+,2^\lambda>2^\mu$ or at least 
$2^\mu=2^{<\lambda}<2^\lambda$ 
($\lambda$ is regular uncountable) then $\lambda \notin 
\{\check I^{\wed}\}_\lambda$.
\end{claim}

\begin{PROOF}{\ref{wd.2}}
By \cite{DvSh:65}, see more in \cite[AP,\S1,pgs.942-961]{Sh:f}.
\end{PROOF}

\begin{remark}
\label{c7}
1) Used in \cite[6.4=constr6.4=f12,stage C]{Sh:482}.

\noindent
2) On $\check I^{\ged}[\lambda]$ see \cite[3.8,3.9]{Sh:309}.
\end{remark}

\begin{definition}
\label{cdl.1}
For $\lambda$ regular uncountable let $\check I[\lambda] = 
\check I^{\ged}[\lambda]$ be the family of sets 
$S \subseteq \lambda$ which have a witness 
$(E,\bar{\cP})$ for $S \in \check I^{\ged}[\lambda]$, which means 
\mn
\begin{enumerate}
\item[$(*)$]  $E$ is a club of $\lambda,{\cP}=
\langle {\cP}_\alpha:\alpha<\lambda\rangle,{\cP}_\alpha 
\subseteq {\cP}(\alpha),|{\cP}_\alpha|<\lambda$, 
and for every $\delta \in E \cap S$ there is an unbounded subset 
$C$ of $\delta$ of order $<\delta$ such that $\alpha\in C \Rightarrow 
C \cap \alpha \in {\cP}_\alpha$. 
\end{enumerate}
\end{definition}

\begin{claim}
\label{c13}
Let $\lambda$ be regular uncountable.

\noindent
1) For $S \subseteq \lambda,S \in \check I^{\ged}[\lambda]$ \Iff \,
equivalently there is a pair $(E,\bar{a})$, $E$ 
is a club of $\lambda,\bar a = \langle a_\alpha:\alpha <
\lambda\rangle,a_\alpha \subseteq \alpha,\beta \in a_\alpha 
\Rightarrow a_\beta=a_\alpha \cap \beta$ 
and $\delta \in E \cap S \Rightarrow \delta=\sup 
(a_\delta) > \otp(a_\delta)$ (or even $\delta=\sup(a_\delta),
\otp(a_\delta)= \cf(\delta)<\delta$.

\noindent
2) If $\kappa < \lambda$ are regular, \then \, there is a stationary
   $S \subseteq S^\lambda_\kappa$ in $\check I^{\ged}[\lambda]$.
\end{claim}

\begin{remark}
\label{c15}
Used in \cite[\S4]{Sh:482}.
\end{remark}

\begin{definition}
\label{c18}
1) For an ideal $\bold I$ on $I$ 
and $f \in {}^\theta(\Ord \setminus \{0\})$ let 
$T_{\bold I}(f) = \sup\{|{\cF}|:{\cF}\subseteq 
\prod\limits_{t\in I} f(t)$ and $h \ne g \in {\cF}$
implies $\{t\in I:h(t) \ne g(t)\}\in \bold I\}$, generally 
$T_{\bold I}(f) = \sup\{|{\cF}|:\cF \in \Xi\}$ where $\Xi$ is the set
of $\cF$ such that:
\mn
\begin{enumerate}
\item[$(a)$]  $\cF \subseteq {}^I{\Ord}$
\sn
\item[$(b)$]   $g \in {\cF}$ implies $g <_{\bold I} h$
\sn
\item[$(c)$]   $h \ne g \in {\cF}$ implies $\{t\in I:h(t) \ne g(t)\} \in 
\bold I$,
\end{enumerate}
\mn
(if $(\forall t) f(t) \ge 2^\kappa$ 
the supremum is obtained and only $f/ \bold I$ matters).

\noindent
1A) We may replace $\bold I$ by the dual ideal. 

\noindent
2) For a partial order $\bbP$ let $\tcf(\bbP)$, the true cofinality of
   $\bbP$ be equal to $\lambda$ \when \, $\lambda$ is a regular
   cardinal and some sequence $\langle p_\alpha:\alpha <
   \lambda\rangle$ witness this which means:
\mn
\begin{enumerate}
\item[$\bullet$]  $\alpha < \beta < \lambda \Rightarrow p_\alpha
  <_{\bbP} p_3$
\sn
\item[$\bullet$]  if $q \in \bbP$ then for some $\alpha < \lambda$ we
  have $q <_{\bbP} p_\alpha$.
\end{enumerate}
\end{definition}

\begin{claim}
\label{c21} 
Assume that $\langle 2^{\lambda_i}:i<\delta\rangle$ is strictly increasing
and $\mu=\sum\langle \lambda_i:i<\delta\rangle< 2^{\lambda_0}$. 
\Then \, for arbitrarily large regular cardinals $\lambda < \mu$ there is 
tree with $<\mu$ nodes and $\ge 2^{\lambda_0},\kappa$-branches 
(hence a linear order of cardinality $<\mu$ with 
$\ge 2^{\lambda_0}>\mu$ Dedekind cuts with both cofinality exactly $\lambda$).
\end{claim}

\begin{remark}
\label{c22}
This is used in \cite[3.28=L3c.16]{Sh:E58} and will 
be used in proving the properties from \cite{Sh:511}.
\end{remark}

\begin{PROOF}{\ref{c21}}
By \cite[3.4]{Sh:430}.
\end{PROOF}

\begin{claim}
\label{prd.f}
Assume:
\mn 
\begin{enumerate}
\item[$(A)$]  $\lambda = \cf(\lambda) \ge \mu>2^\kappa$,
\sn
\item[$(B)$]  $\dot D$ is a $\mu$-complete\footnote{in parts 
(0),(1),$\mu=(2^\kappa)$ is O.K.} filter on $\lambda$,
\sn
\item[$(C)$]  $f_\alpha:\kappa\rightarrow \Ord$ for $\alpha<\lambda$,
\sn
\item[$(D)$]  $\dot D$ contains the co-bounded subsets of $\lambda$. 
\end{enumerate}
\mn
\Then

\noindent
0) We can find $w\subseteq \kappa$ and $\bar \beta^*=
\langle \beta^*_i:i<\kappa\rangle$ such that: $i\in \kappa\setminus w
\Rightarrow \cf(\beta^*_i)> 2^\kappa$ and for every $\bar \beta 
\in \prod\limits_{i\in \kappa\setminus w} \beta^*_i$ for $\lambda$ 
ordinals $\alpha<\lambda$ (even a set in $\check{\cD}^+$) we have 
$\bar \beta< f_\alpha\rest (\kappa\setminus w)< \bar \beta^* \rest 
(\kappa\setminus w), f_\alpha\rest w=\bar \beta^*\rest w$, and 
$\sup\{\beta^*_j:\beta^*_j<\beta^*_i\} < f_\alpha(i)<\beta^*_i$.

\noindent
1) We can find a partition $\langle w^*_\ell:\ell<2\rangle$ 
of $\kappa,X \in \check{\cD}^+$ and $\langle A_i:i<\kappa\rangle, 
\langle\bar \lambda_i:i<\kappa\rangle, \langle h_i:i<\kappa\rangle,
\langle n_i:i<\kappa\rangle$ such that:
\mn
\begin{enumerate}
\item[$(a)$]   $A_i \subseteq \Ord$,
\sn
\item[$(b)$]  $\bar \lambda_i=\langle\lambda_{i,\ell}:\ell <
  n_i\rangle$ and $2^\kappa<\lambda_{i,\ell} \le \lambda_{i,\ell+1} 
\le \lambda$, 
\sn
\item[$(c)$]  $h_i$ is an order preserving function from 
$\prod\limits_{\ell<n_i} \lambda_{i,\ell}$ onto $A_i$ so 
$n_i=0 \Leftrightarrow |A_i|=1$. 
(The order on $\prod\limits_{\ell<n_i} \lambda_{\ell,i}$ being 
lexicographic, $<_{\ell x}$),
\sn
\item[$(d)$]  $i<\kappa$ and $\alpha \in X \Rightarrow f_\alpha(i) 
\in A_i$, and we let $f^*_\alpha(i,\ell)=[h^{-1}_i (f_\alpha
(i))] (\ell)$, so $f^*_\alpha \in 
\prod\limits_{\stackrel{i<\kappa}{\ell<n_i}} \lambda_{i,\ell}$,
\sn
\item[$(e)$]  $i\in w^*_0 \Leftrightarrow n_i=0$ (so $|A_i|=1$),
\sn
\item[$(f)$]  if $i\in w^*_1$ then $|A_i| \le \lambda$, hence 
$|\bigcup\limits_{i\in w^*_1} A_i| \le \lambda$, 
\sn
\item[$(g)$]  if $g\in \prod\limits_{\stackrel{i<\kappa}{\ell<n_i}} 
\lambda_{i,\ell}$ then $\{\alpha\in X:g<f^*_\alpha\}
\in \check{\cD}^+$ and letting $\beta^*_j = \sup\Rang(h_i)$, 
clearly the condition of part $(\gamma)(0)$ holds 
\sn
\item[$(h)$]  if $\dot D$ is $(|\alpha|^{\kappa})^+$-complete for 
any $\alpha<\mu_1$ then $\mu_1 \le \sup\{\lambda_{i,\ell}:i \in
w^*_1$; and $\ell < n_i\} \le \lambda$ when $w^*_1 \ne \emptyset$ 
(so, e.g., if $\mu=\lambda$ and assuming $\GCH$
\[
\sup\{\cf(\lambda_{i,\ell}):i \in w^*_1 \text{ and } \ell < n_i\}=\lambda).
\]
\end{enumerate}
\mn
2) In part (1) we can add $(*)_1$ to the conclusion if (E) below
holds,
\mn
\begin{enumerate}
\item[$(*)_1$] if $\lambda_{i,\ell}\in [\mu,\lambda)$ 
then $\lambda_{i,\ell}$ is regular.
\sn
\item[$(E)$]  For any set ${\ga}$ of $\le \kappa$ singular 
cardinals from the interval $(\mu,\lambda)$, we have 
$\max \pcf\{\cf(\chi):\chi\in {\ga}\}<\lambda$.
\end{enumerate}
\mn

3)  Assume in part (1) that (F) below holds. \Then \, we can demand $(*)_2$.
\mn
\begin{enumerate}
\item[$(*)_2$]  $\lambda^i_\ell \ge \mu_1$ for $i\in w_2,\ell<n_i$.
\sn
\item[$(F)$]   $\cf(\mu_1)>\kappa$ and $\alpha<\mu_1 \Rightarrow \dot D$ 
is $[|\alpha|^{\le\kappa}]^+$-complete.
\end{enumerate}
\mn
4) If in part (1) in addition (G) below holds, \then \, we can add:
\mn
\begin{enumerate}
\item[$(*)_3$]  $\lambda \in \pcf_{\partial-\complete} 
\{\lambda^i_\ell:i \in w^*_1$; and $\ell<n_i\}$ if $w^*_1 \ne 
\emptyset$, moreover
\sn
\item[$(*)_4$]  if $\ell_i<n_i$ for $i \in w^*_1$ then $\lambda\in 
\pcf_{\partial-\complete}\{\cf(\lambda^i_{\ell_i}):i\in w^*_1\}$.
\sn
\item[$(G)$] 
\sn
\item[${{}}$] $(i) \quad (\forall \alpha<\lambda) (|\alpha|^{<\partial}
<\lambda)$ and $\partial=\cf(\partial)>\aleph_0$,
\sn
\item[${{}}$]  $(ii) \quad \dot D$ is $\lambda$-complete
\sn
\item[${{}}$]  $(iii) \quad f_\alpha \ne f_\beta$ for $\alpha \ne \beta$ 
(or just $\alpha \ne \beta \in X$ for some $X \in \dot D^+$)
\end{enumerate}
\mn
5)  If in part (1) in addition (H) below holds \then \, we can add :
\mn
\begin{enumerate}
\item[$(*)_5$]  if $m<m^*,A \in \bold J_m$ and $\ell_i<n_i$ 
for $i\in \kappa \setminus A$ (so $w^*_0\subseteq A$) \then \, 
$\lambda \in \pcf\{\lambda^i_{\ell_i}:i \in \kappa\setminus A\}$. 
\sn
\item[$(H)$] 
\sn
\item[${{}}$]  $(i) \quad m^*<\omega$ and $\bold J_m$ is an 
$\aleph_1$-complete ideal on $\kappa$ for $m<m^*$,
\sn
\item[${{}}$]  $(ii) \quad \dot D$ is $\lambda$-complete. 
\end{enumerate}
\end{claim}

\begin{PROOF}{\ref{prd.f}}
By \cite[7.1=L7.0]{Sh:620}.
\end{PROOF}

\begin{claim}
\label{pcf.1}
Assume that $\bar \lambda=\langle \lambda_i:i<\kappa\rangle$ is a sequence 
of regular cardinals $>\mu$ and $J$ is an ideal of $\kappa$ 
and $\lambda$ is a regular cardinal.

\noindent
1) If $\prod\limits_{i<\kappa}\lambda_i/ \bold J$ is $\lambda^+$-directed 
then we can find $\lambda'_i = \cf(\lambda'_i)\in (\mu,\lambda_i)$
such that:
\mn
\begin{enumerate}
\item[$(a)$]  $\prod\limits_{i<\kappa}\lambda'_i/ \bold J$ has 
true cofinality $\lambda$
\sn
\item[$(b)$]  if $\lambda > \lim_{\bold J} \langle \lambda_i:
i < \kappa\rangle = \mu_* > \cf(\mu_*)$ then $\lim_{\bold J}
\langle \lambda'_i:i<\kappa\rangle=\mu^*$
\sn
\item[$(c)$]  there is an $<_{\bold J}$-increasing sequence 
$\langle f_\alpha:\alpha<\lambda\rangle$
of members of $(\prod\limits_{i<\kappa}\lambda_i,<_{\bold J})$ 
and is $\mu^+_*$-free,
i.e. if $A \subseteq \lambda,|A| \le \mu_*$ then there is a sequence 
$\langle u_\alpha:\alpha\in A\rangle$ such that 
$u_\alpha \in J$ and $\alpha\in A$ and $ \beta \in A$ and 
$\alpha < \beta$ and $\epsilon \in \kappa \setminus u_\alpha
\setminus u_\beta\Rightarrow f_\alpha (\epsilon)<f_\beta(\epsilon)$.
\end{enumerate}
\end{claim}
 
\begin{remark}
Used in \cite[1.16=L7.7]{Sh:331}.
\end{remark}

\begin{PROOF}{\ref{pcf.1}}
By \cite[\S6]{Sh:430}.
\end{PROOF}

\begin{theorem}
\label{4.EK} 
1)  Assume that $\mu=\mu^{<\kappa}<\lambda \le 2^\mu$ \then \, there is a 
sequence 
$\langle f_i:i<\mu\rangle$ of functions from $\lambda$ to $\mu$ such that for 
every $u \subseteq \mu$ of cardinality $<\kappa$ of function $g$ from $u$ to 
$\mu$, for some $i<\mu$ we have $g\subseteq f_i$.
\end{theorem}

\begin{remark}
Used in \cite[1.11=L7.6]{Sh:331}.
\end{remark}

\begin{PROOF}{\ref{4.EK}}
This is Engelking-Karlowicz \cite{EK65}.
\end{PROOF}

\begin{theorem}
\label{4.Ha}
(Hajnal free subset theorem). 
If $f:\lambda\rightarrow [\lambda]^{<\kappa}$ and 
$\lambda>\kappa \ge \aleph_0$ then some $A \in [\lambda]^\lambda$ 
is $f$-free which means that $\alpha \ne \beta \in A
\Rightarrow \alpha\notin f(\beta)$.
\end{theorem}

\begin{PROOF}{\ref{4.Ha}}
This is \cite{Haj62}.
\end{PROOF}

\begin{definition}
\label{prf.2}
1)  For $\mu$ singular let $\pp(\mu)=\sup\{\lambda$: 
for some filter $\bold J$ on $\cf(\mu)$ and sequence 
$\langle\lambda_i:i < \cf(\mu)\rangle$ of regular cardinals 
$<\mu$ such that $\mu'<\mu\Rightarrow \{i:\lambda_i> \mu'\} \in \bold J$, the 
product $\prod\limits_{i<\cf \mu}\lambda_i/ \bold J$ 
has true cofinality $\lambda\}$.

For $\mu$ singular $\pp^+(\mu) = \Min\{\lambda:\lambda$ regular and 
there are no $\bold J$ and $\lambda_i$ as above$\}$.

\noindent
2)  For a set ${\ga}$ of regular cardinals $\ge |{\ga}|$ let 
$\pcf({\ga}) = \{\cf(\prod\limits_{\theta\in{\ga}} (\theta,<)/\dot D):
\dot D$ an ultrafilter on ${\ga}\}$.

\noindent
3) If ${\ga}$ is as above, $\bold J$ is an ideal on ${\ga}$ then we let 
$\pcf_{\bold J}({\ga}) = \{\cf(\Pi{\ga}/\dot D):\dot D$ is an ultrafilter on 
${\ga}$ disjoint to $I\}$.
\end{definition}

\begin{remark}
Used in \cite[1.16=L7.7]{Sh:331}.
\end{remark}

\begin{remark}
Used in \cite[2.15=L7.9]{Sh:331}.
\end{remark}

\begin{claim}
\label{pcf.6}
If $\mu \ge \kappa = \cf(\kappa) > \aleph_0$. \Then \, there is a stationary  
${\cS} \subseteq [\mu]^{<\kappa}$ of cardinality 
$\cf([\mu]^{<\kappa},\subseteq)$.
\end{claim}

\begin{remark}
Used in \cite[5.3]{Sh:482}, stage E statement of $\otimes_3$.
\end{remark}

\begin{PROOF}{\ref{pcf.6}}
By \cite[\S1]{Sh:420}.
\end{PROOF}

\begin{claim}
\label{pcf.6a}
Assume that $\lambda$ is singular of uncountable cofinality $\kappa,
\langle \lambda_i:i<\kappa\rangle$ is increasing continuous 
with limit $\lambda$ and 
$S=\{\delta<\kappa:\pp(\lambda_\delta)=\lambda_\delta^+\}$ 
is a stationary subset of $\kappa$ \then \, $\pp(\lambda)=\lambda^+$.
\end{claim}

\begin{PROOF}{\ref{pcf.6a}}
By \cite[Ch.II,\S2]{Sh:g}.
\end{PROOF}

\noindent
We repeat \cite[Ch.IX,3.7,pg.384,5]{Sh:e}
\begin{claim}
\label{pcf.7}
Suppose $\lambda = \aleph_{\alpha(*)+\delta},\delta$ a limit ordinal
$< \aleph_{\alpha(*)}$.

\noindent
1) $\pp(\lambda) =^+ \cov(\lambda,\lambda,\cf(\lambda)^+,2)$.

\noindent
2) If $\cf(\delta) \le \kappa \le \delta$ then $\pp_\kappa(\lambda)
=^+ \cov(\lambda,\lambda,\kappa^+,2)$.

\noindent
3) If $\cf(\delta) = \kappa,(\aleph_{\alpha(*)+i})^\kappa <
\aleph_{\alpha + \delta}$ for $i < \delta$ then $\\lambda^\kappa =
 \pp(\lambda)$.

\noindent
4) If $\cf(\delta) = \kappa,(\aleph_{\alpha(*)})^\kappa < \lambda$ then

\[
\lambda^\kappa = \sum\{\pp(\aleph_{\alpha(*)+i}):i \le \delta \text{
  limit}, \cf(i) \le \kappa\}.
\]

\mn
5) $\cS_{< \aleph_{\alpha(*)+1}}(\lambda)$ has a stationary subset of
cardinality

\[
\sum\{\pp(\aleph_{\alpha(*)+i}):i \le \delta \text{ limit}\}.
\]
\end{claim}

\begin{claim}
\label{pcf.8}
Assume $\mu>\kappa = \cf(\mu)$. There is an increasing sequence 
$\langle \lambda_i:i<\kappa\rangle$ of regular cardinals 
$<\mu$ and $\lambda = \tcf(\prod\limits_{i<\kappa} 
\lambda_i,<_{J^{\bd}_\kappa})$ \when
\mn 
\begin{enumerate}
\item[$\circledast$]  $(a) \quad \lambda=\cf(\lambda)\in 
(\mu,\pp^+_\kappa(\mu))$
\sn
\item[${{}}$]  $(b)_1 \quad \mu<\mu^{+\kappa}$ or 
\sn
\item[${{}}$]  $(b)_2 \quad \kappa>\aleph_0$ and 
for some $\mu_0<\mu$ for every $\mu'\in (\mu_0,\mu)$ 
of cofinality $<\kappa$

\hskip35pt  we have $\pp(\mu')<\mu$.
\end{enumerate}
\end{claim}

\begin{remark}
1) Used in \cite[1.16=L7.7]{Sh:331}, \cite[2.20=7.11]{Sh:331},
\cite[3.23=L7.14]{Sh:331}, \cite[3.24=L7.7]{Sh:331}.

\noindent
2) It is helpful in applying \cite[2.13=L7.8I]{Sh:331}.
\end{remark}

\begin{PROOF}{\ref{pcf.8}}
By \cite[Ch.VIII,\S1]{Sh:g}.
\end{PROOF}
\newpage

\section {On normal ideals}

The results here are from \cite{Sh:247}.

\begin{theorem}
\label{d4}
If $\check{\cD}$ is a fine normal filter on ${\cU} = \{u\subseteq \lambda:\cf 
(\sup(u)) \ne \cf(|u|)\}$, and $\lambda$ is regular \then \, 
there are functions $f^*_i$ for $i < \lambda^+$ such that: $\Dom(f^*_i)=
{\cU}, f^*_i(u) \in u$ and for $i \ne j, \{u\in I:f^*_i(u)=
f^*_j(u)\}= \emptyset \mod \check{\cD}$.
\end{theorem}

\begin{remark}
\label{d5}
1) Used in \cite[2.15=L7.9]{Sh:331}.

\noindent
2) So ${\cU} = [\lambda]^{\aleph_1}$ is an interesting case.

\noindent
3) This is a strong form of ``not $\lambda^+$-saturated".
\end{remark}

\begin{PROOF}{\ref{d4}}
We can find $A_i(i<\lambda^+)$ such that:
\mn 
\begin{enumerate}
\item[$(*)_1$]   $A_i$ is a subset of $\lambda$, unbounded in $\lambda$ 
and for $j<i, A_i\cap A_j$ is bounded in $\lambda$
\end{enumerate}
\mn
[e.g. let $A_i(i<\lambda)$ be pairwise disjoint subsets 
of $\lambda$ of power $\lambda$, and then choose
 $A_i(\lambda \le i<\lambda^+)$ by induction on $i$ on such that the 
relevant demands hold. Assuming to $i \in [\lambda,\lambda^+)$ let 
$\{j:j<i\}$ be listed as $\{j^i_\alpha:\alpha<\lambda\}$, 
and let $A_i=\{\gamma^i_\beta:\beta<\lambda\}$
where $\gamma^i_\beta = \Min(A_{j_\beta}
\setminus \bigcup\limits_{\alpha<\beta} A_{j_\alpha})$, 
listed without repetitions it exists as 
$|A_{j_\beta}\cap A_{j_\alpha}|<\lambda=\cf (\lambda)$ for $\alpha<\beta$].

For $i<\lambda^+$ let $g_i:i\rightarrow \lambda$ be such that 
$\{A_j\setminus g_i(j):j<i\}$ are 
pairwise disjoint. Let $f_i$ be the strictly increasing function 
from $\lambda$ onto $A_i$ 
(for $i<\lambda^+)$ hence $\alpha < \lambda \Rightarrow 
f_i(\alpha)\ge \alpha$. 
So $C_i = \{u \in \cU:u$ is closed under $f_i$ and 
$\alpha\in u\Rightarrow \alpha+1\in u\}$ belongs 
to $\check{\cD}$. 
For each $u\in {\cU}$ let $u=\{x^u_\alpha:\alpha<|u|\}$. 

Now for each $u\in C_i$ the set $u\cap A_i$ is unbounded in $u$, 
(by the choice of $C_i$ and $f_i$) so for some 
$\alpha_i(u)<|u|$, the set $A_i\cap \{x^u_\alpha:\alpha<\alpha_i(u)\}$ 
is unbounded in $u$.
(Why?  Recall that $\cf(\sup u) \ne \cf(|u|)$ because $u\in {\cU}$). 

Next for $i<\lambda^+$ let $h_i$ be a one-to-one function from $\lambda$ onto 
$\lambda\cup\{j:j<i\}$ and define by induction on $i$:

\begin{equation}
\begin{array}{clcr}
C^1_i = \{u\subseteq i\cup\lambda: \, &u \text{ is closed under } 
h_i,h_i^{-1} \text{ and } u \cap \lambda\in {\cU} \\
\, & u\cap \lambda \text{ is closed under } f_i, f_i^{-1}, \\
\, & u \text{ is closed under } g_j, (\text{when } j \in u \text{ or } j=i)\\
\, & \text{and for every } j \in u \text{ we have } 
u \cap (j \cup \lambda)\in C^1_j\}.
\end{array}
\end{equation}

\mn
Clearly $C^1_i \rest\lambda = \{u\cap \lambda:a\in C^1_i\}$ belongs to 
$\check{\cD}$, and by the choice of $h_i$ for each $u\in {\cU}$
there is at most one $u'\in C^1_i$ satisfying 
$u'\cap\lambda=a$, namely $h_i``(u)$.

Now we define for $i<\lambda^+$ a functions $\xi_i$ and 
$d_i$ with domain ${\cU}$.

\[
\xi_i(u) = \otp(\{j \in h_i(u):\alpha_j(u)=\alpha_i(u)\}, 
\]

\mn
$d_i(u)=(\alpha_i(u),\xi_i(u))$ if $h_i(u)\cap \lambda=u$ and
$h_i(u)\in C^1_i$ and $d_i(u) = \Min(u)$ otherwise.

Now we shall finish by showing:
\mn
\begin{enumerate}
\item[$(A)$]   for $i_1 \ne i_2<\lambda^+$ we have 
$\{u\in \cU:d_{i_1}(u)=d_{i_2}(u)\}=\emptyset \mod \check{\cD}$
\sn
\item[$(B)$]   for $a\in {\cU},\{d_i(u):i<\lambda^+\}$ has 
cardinality $\le |u|$.
\end{enumerate}
\mn
Why does this suffice?  As for each $u\in {\cU}$ by clause (B) 
we can find a one-to-one function ${\bold f}_u$ from 
$\{d_i(u):i<\lambda^+\}$ into $u$ and now use the $\lambda^+$ functions 
$\langle {\bold f}_u(d_i(u)):i < \lambda^+\rangle$, that is for 
$i<\lambda^+$ we define the function $f_i^*$ with 
domain ${\cU}$ such that $f^*_i (u)\in u$ by 
$f^*_i(u)=: {\bf f}_u (d_i(u))$, now by clause (A) we have 
$i<j<\lambda^+\Rightarrow f^*_i \ne f^*_j \mod \check{\cD}$.
\medskip

\noindent 
\underline{Proof of Clause (A)}: 

\Wilog \, $i_1<i_2$ 
and we assume that $\lambda \le i_1$ for notational simplicity. 
Clearly ${\cU}':=\{u\in {\cU}:h_{i_2}(u)\in C^1_{i_2}$ and $i_1 \in 
h_{i_2}(u)$ (hence $h_{i_1}(u) = h_{i_2}(u) 
\cap i_1\in C^1_{i_1})\}$ belongs to $\check{\cD}$.  Let $u$ be in it, 
and assume that $d_{i_1}(u)=d_{i_2}(u)$. For $\ell=1,2$ in 
the definition of $d_{i_\ell} (u)$ the first 
case applies so $d_{i_\ell} (u) = (\alpha_{i_\ell}(u),\xi_{i_\ell}(u))$
hence by the first coordinate $\alpha_{i_1}(u)=\alpha_{i_2}(u)$. Now 
$\{\xi\in h_{i_1}(u):\alpha_\xi(u)=\alpha_{i_1}(u)\}$ 
is an initial segment of $\{\xi\in h_{i_2}(u):\alpha_\xi(u)=\alpha_{i_2}\}$ 
(as $a \in {\cU}'$) and a proper one (as $i_1$ 
belongs to the latter but not the former).
As the ordinals are well ordered, the order types 
$\xi_{i_1}(u),\xi_{i_2}(u)$ are not equal. That means that 
the second coordinates in the $d_{i_1}(u), d_{i_2}(u)$ 
are distinct.  So $d_{i_1}(u) \ne d_{i_2}(u)$ is true when 
$i_1 \ne i_2, a \in {\cU}'$ as required. 
\medskip

\noindent
\underline{Proof of Clause (B)}:
 
The number of possible $\alpha_i(u)$ is $\leq|u|$, and the number of 
order types of well orderings of power $<|u|$ is $|u|$ 
hence by $(*)$ below, the number of pairs $(\alpha_i(u),\xi_i(u))$ is 
$\le |a| \times |u|=|u|+\aleph_1$ and recalling the additional 
 value $\Min(u)$ we are done. So it suffices to prove: 
\mn
\begin{enumerate}
\item[$(*)$]   for $i<\lambda^+, u\in C^1_i$, the set $w=\{j\in u:
\alpha_j(u\cap \lambda)=\alpha_i(u\cap\lambda)\}$ has power $<|u|$.
\end{enumerate}
\mn
Why $(*)$ holds?   Clearly for $j\in w$ the set

\[
A_j\cap \{x^u_\alpha:\alpha<\alpha_i(u\cap\lambda)\}
\]

\mn 
is unbounded in $u\cap\lambda$ \underline{but}  
$A_j\cap g_i(j)$ is bounded in $u\cap\lambda$ (as $u$ is closed under $g_i$)
hence

\[
B_j := (A_j\setminus g_i(j))\cap 
\{x^u_\alpha:\alpha<\alpha_i(u\cap\lambda)\}
\]

\mn
is an unbounded subset of $u\cap\lambda$, hence non-empty. 

But $\langle B_j:j\in w\rangle=
\langle B_j:j\in
u,\alpha_j(u\cap\lambda)=\alpha_i(u\cap\lambda)\rangle$ is a 
sequence of pairwise disjoint subsets of 
$\{x^u_\alpha:\alpha<\alpha_i(u\cap \lambda)\}$
(by the choice of $g_i$). As they are non-empty their number 
is $\le |\{x^u_\alpha:\alpha <\alpha_i(u\cap\lambda)\}| <|u|$.  So
have proved $(*)$, which suffice.
\end{PROOF}

\begin{claim}
\label{d7} 
Let $\check{\cD}$ be a fine normal filter on ${\cU}
=[\lambda]^{<\kappa},\lambda$ singular of cofinality $\partial>\aleph_0$ 
and $(\forall u\in {\cU}) (|u|\ge \partial$ and 
$\cf(|u|)\ne \partial$ and $\partial=\sup(\partial\cap u))$ and 
$\Rk(|u|,\check{\cD}^{\cb}_\partial)\le |u|^+$.

\Then \, there are functions $f_i$ for $i<\lambda^+,\Dom(f_i)={\cU}, 
(\forall u\in {\cU})[f_i(u)\in u]$ and for $i\neq j$ we have 
$\{u\in I:f_i(u) = f_j(u)\}=\emptyset \mod \check{\cD}$.
\end{claim}

\begin{PROOF}{\ref{d7}}
Let $\partial=\cf(\lambda),\lambda=
\sum\limits_{\zeta<\partial}\lambda_\zeta$, each $\lambda_\zeta$ regular, 
$\sum\limits_{\xi<\zeta}\lambda_\xi<\lambda_\zeta<\lambda$ for 
$\zeta<\partial$. We can find for 
$i<\lambda^+$ functions ${\bold f}_i$ from $\partial$ 
to $\lambda,\sum\limits_{\xi<\zeta}\lambda_\xi < \bold f_i(\zeta)
<\lambda_\zeta$ such that for $i<j<\lambda^+$ there is $\xi<\partial$
such that

\[
\xi \le \zeta<\partial\Rightarrow {\bold f}_i(\zeta)<{\bold f}_j(\zeta).
\]

\mn
Let again $u=\{x^u_\alpha:\alpha<|u|\}$, so for each 
$i<\lambda^+$ and $u \in {\cU}$, if $\Range({\bold f}_i\rest u)$ is 
unbounded in $u$ then let $\alpha_i(u)<|u|$ be minimal such that  
$(\Range({\bold f}_i\rest u))\cap \{x^u_\alpha:\alpha< 
\alpha_i(u)\}$ is unbounded in $u$ 
(and $\alpha_i(u) = \Min(u)$ otherwise). 

Now for $i<\lambda^+$ we define functions $\xi_i,d_i$ with domain ${\cU}$ 
($h_i$ is a one-to-one function from $\lambda$ onto $i\cup\lambda$): 

\[
\xi_i := \otp\{j \in h_i(u):\alpha_j(u)=\alpha_i(u)\}
\]

\mn
$d_i(u)$ is $(\alpha_i(u),\xi_i(u))$  when $u = h_i(u)\cap\lambda$ and 
$(\forall\zeta \in(u\cap \cf \lambda)){\bold f}_i(\zeta)\in u$
and $(\forall j \in u)(u=h_j(u)\cap \lambda)$ and $d_i(x) = \Min(u)$ otherwise.

We finish as in \ref{d4}. 
\end{PROOF}

\begin{remark}
\label{d8}
1)  $\check{\cD}^{\cb}_\partial$ is the filter of co-bounded 
subsets of $\partial$.

\noindent
2) Really we use $\Rk(|u|,\check{\cD}^{\cb}_{\partial}) \le |u|^+$ 
just to get, that for every $\zeta<|u|$ for some $\xi_\zeta<|u|^+$
we have 
\mn
\begin{enumerate}
\item[$(*)$]   there are no $f_i:\partial\rightarrow \zeta$ 
for $i<\xi_\zeta, [i<j \Rightarrow f_i <_{\check{\cD}^{\cb}_\partial} f_j]$.
\end{enumerate}
\mn
We should observe that for $u\in {\cU}, u\cap\partial$ 
has order type $\partial$. 

Note that if for each $\zeta<|u|$ there is such $\xi_\zeta$ then $\xi(*)=
\bigcup\limits_{\zeta<|u|} \xi_\zeta$ is $<|u|^+$ and work for all $\zeta$'s.
\end{remark}

\begin{claim}
\label{d10}
Suppose $\kappa \le \partial= \cf(\lambda)<\lambda,{\cU} 
\subseteq \{u\in [\lambda]^{<\kappa}:\, \cf(|u|)\ne \quad 
\cf(\sup (u\cap\partial))$  and $\Rk(|u|,
\check{\cD}^{\cb}_{\cf(\sup(u\cap\partial))})\le |u|^+$  when 
$\cf(\sup u)>\aleph_0$ and $|u|^{\aleph_0}=|u|$ 
or just when  $(\forall \mu<|u|)(\mu^{\aleph_0}\le |u|)$ and 
$\cf(\sup(u))=\aleph_0\}$, and $\check{\cD}$ 
a normal fine filter on ${\cU}$. 

\Then \, there are for $i<\lambda^+$ functions 
$f_i:{\cU} \rightarrow \lambda,f_i(u)\in u$ such that  
for $i \ne j$ we have $\{u\in I:f_i(u)=f_j(u)\}=\emptyset \mod \check{\cD}$.
\end{claim}

\begin{PROOF}{\ref{d10}}
Let ${\bold f}_i,\lambda_\zeta$ be as in the proof of 
\ref{d7}, $u=\{x^u_\alpha:\alpha<|u|\}$. 
Let $h_i$ be a one-to-one function from $\lambda$ onto $\lambda\cup\{j:j<i\}$. 
For each $i$ the set $C^1_i := \{u\in {\cU}:u$ is closed under 
${\bold f}_i$, and $(\Range({\bold f}_i))\cap u$ 
is unbounded in $u,h_i(u)\cap\lambda=u$ and $u\in C^1_j$ for 
$j\in h_i(u)$ and $\cf(\sup u)=\cf(\sup(u\cap\partial))\}$ belongs 
to $\check{\cD}$, and for $u\in C^1_i$ let $\alpha_i(u)<|u|$ be
minimal such that $(\Range({\bold f}_i))\cap 
\{x^u_\alpha:\alpha<\alpha_i(u)\}$ is unbounded in $u$. 

We then let 

\[
\xi_i (u) = \otp\{j:j\in h_i(u),\alpha_j(u)=\alpha_i(u)\}
\]

\[
d_i(a) = \alpha_i(u), \text{ if } u\in C^1_i,
\]

\[
\Min(u) \quad \text{ otherwise}.
\]

\mn
and we proceed as in the proof of \ref{d4}, \ref{d7} 
(and see \ref{d8}).
\end{PROOF}

\begin{definition}
\label{d12}
1) For a filter $D$ on $[\kappa]^{<\theta}$ let $\diamondsuit_D$ mean:
fixing any countable vocabulary $\tau$ there are $S\in D$ and $N=\langle
N_a: a\in S\rangle$, each $N_a$ a $\tau$-model with universe $a$, such that
for every $\tau$--model $M$ with universe $\lambda$ we have

\[
\{a\in S: N_a\subseteq M\}\neq \emptyset \mod D.
\] 

\mn
2) Similarly, let $\diamondsuit^*_D$ (or $\diamondsuit^*(D))$ mean:
there is a $\bar P = \langle P_u:u \in [\lambda]^{< u}\rangle$ such
that:
\mn
\begin{enumerate}
\item[$(a)$]  $P_u \subseteq P(u)$ has cardinality $\le |u|$
\sn
\item[$(b)$]  $\{u \in [X]^{<u}:X \cap u \in P_u\} \in D$ for every $X
  \subseteq \lambda$.
\end{enumerate}
\mn
Recall that for two filters $\cD$ and $U$ on $[\lambda]^{<u}$ the set
$\cD + U$ is defined to be the smallest filter on $[\lambda]^{<u}$
which extends both $\cD$ and $U$.
\end{definition}

\begin{fact}
\label{d14}
1)  For ${\cU}_1\subseteq {\cU}_2 \subseteq 
[\lambda]^{<\kappa}$ and $\check{\cD}_1 \subseteq \check{\cD}_2$ 
normal fine filter we have on $[\lambda]^{<\kappa}$,
\mn
\begin{enumerate}
\item[$(i)$]   $\diamondsuit^*(\check{\cD}_1 + {\cU}_2)
\Rightarrow \diamondsuit^* (\check{\cD}_2+ {\cU}_1)$
\sn
\item[$(ii)$]   $\diamondsuit^*(\check{\cD}_1+{\cU}_2)
\Rightarrow \diamondsuit(\check{\cD}_1+{\cU}_2)$
\sn
\item[$(iii)$]  $\diamondsuit(\check{\cD}_2+{\cU}_1)
\Rightarrow \diamondsuit(\check{\cD}_1+{\cU}_2)$
\sn
\item[$(iv)$]  $\diamondsuit^*(\check{\cD}_1+{\cU}_2)
\Rightarrow \diamondsuit (\check{\cD}_2+{\cU}_2)$
\end{enumerate}
\mn
(remember $\check{\cD}_{<\kappa}(\lambda)+{\cU}_1
\subseteq \check{\cD}$ for any fine normal filter $\check{\cD}$ on ${\cU}_1$).

\noindent
2)  Suppose $\kappa<\lambda=\lambda^{<\kappa}$, and we let

\begin{equation*}
\begin{array}{clcr}
{\cU} = \{a:\, &\text{ for some } \theta,a \in T_{\kappa,\lambda} 
(N^0_\theta), |u|^\theta=|u| \\
  &\text{ or } u \in T_{\kappa,\lambda}(N^1_\theta), \text{ and } 
\cf(|u|) \ne \theta \wedge (\forall \partial<|u|)\partial^\theta\leq|u| \\
  &\text{ or } (\exists \chi,\partial,\alpha)(2^\chi \le 
\lambda\cap\lambda=\chi^{+\alpha} \wedge |u|^{<\partial} =
|u| \wedge (\forall \gamma <\alpha) \\
  % [\cf(u\cap \chi^{+(\gamma+1)}<\partial] \frown \alpha<\partial)\}
\end{array}
\end{equation*}

Suppose further ${\cU} \ne \emptyset \mod \check{\cD}_{<\kappa}(\lambda)$. 
\Then \, $\diamondsuit {}^*(\check{\cD}_\kappa(\lambda) + {\cU})$.
\end{fact}

\begin{remark}
\label{d16}
Used in the proof of \cite[2.13=L7.8I]{Sh:331}.
\end{remark}

\begin{PROOF}{\ref{d14}}
By straightforward generalization of the proof for the 
case $\lambda=\kappa$, due to Kunen for (1), (i.e., 1(ii), 
the rest being trivial) Gregory and Shelah for (2) 
(see e.g. \cite{Sh:108}). 
I.e. for 1)(ii), suppose $\langle{\cP}_u:u\in{\cP}_{<\kappa}
(\lambda)\rangle$ exemplifies $\diamondsuit^*(\check{\cD}_1+J)$. 
Let ${\cP}_u= \{A^u_i:i\in u\}$. 
Let $\pr$, i.e. $\pr(-,-)$ be a pairing function on $\lambda$, and for each 
$i<\lambda, u\in {\cP}_{<\kappa}(\lambda)$ let 

\[
B^i_u = \{\alpha:\alpha\in u, <\alpha,i> \in A^u_i\}.
\]

\mn
So $B^i_u\subseteq u_i$ is $\langle B^i_u:u\in [\lambda]^{<\kappa}\rangle$ 
a $\diamondsuit(\check{\cD}_1)$-sequence 
for some $i$?   If yes we finish, if not let $B^i \subseteq 
\lambda$ exemplify this i.e., 

\[
C^i = \{u\in [\lambda]^{<\kappa}:B^i\cap u \ne B^i_u\} \in \check{\cD}_1.
\]

\mn
Hence

\[
C = \{u\in [\lambda]^{<\kappa}:(\forall i\in u) u \in C^i, \text{ and } u 
\text{ is closed under } \pr(-,-)\} \in \check{\cD}
\]

\mn
and let

\[
A = \{\pr(\alpha,i):\alpha \in B^i \text{ and } i\}.
\]

\mn
So for some $u\in C$, $A \cap u\in {\cP}_u$ hence for some 
$i\in A, A\cap u=A^u_i$ hence $B^i\cap u=B^i_u$ contradiction. 
\end{PROOF}
\newpage

%%\bibliographystyle{abbrv}

%%\bibliographystyle{alphacolon}
%%\bibliography{lista,listb,listx,listf,liste,listz}

\begin{thebibliography}{DvSh65}

%%\bibitem{A:barwise1977}
%%J.~Barwise.
%%\newblock On moschovakis closure ordinals.
%%\newblock {\em The Journal of Symbolic Logic}, 42(2):292--296, June 1977.

%%\bibitem{IP:bosse1993}
%%U.~Bosse.
%%\newblock An ``{E}hrenfeucht-{F}ra\"{i}ss\'{e} game'' for fixpoint logic and
%%  stratified fixpoint logic.
%%\newblock In E.~B\"{o}rger, G.~J\"{a}ger, H.~K. B\"{u}ning, S.~Martini, and
%%  M.~M. Richter, editors, {\em Selected Papers from the 6th Workshop on
%%  Computer Science Logic}, CSL '92, pages 100--114, London, UK, Sept. 1993.
%%  Springer-Verlag.

\bibitem[EK65]{EK65} 
Ryszard Engelking and Monika Kar{\l}owicz,
\newblock {\em Some theorems of set theory and their topological consequences},
\newblock Fundamenta Math. {\bf 57} (1965), 275--285.

\bibitem[Haj62]{Haj62}
Andras Hajnal,
\newblock {\em Proof of a conjecture of S. Ruziewicz},
\newblock Fundamenta Math. {\bf 50} (1961/1962), 123--128.

\bibitem[Lav71]{Lav71} 
Richard Laver,
\newblock {\em On Fraiss\'e's order type conjecture},
\newblock Annals of Mathematics {\bf 93} (1971), 89--111.


\bibitem[Mar75]{Mar75} 
Donald A. Martin,
\newblock {\em Borel Determinacy},
\newblock Annals of Mathematics {\bf 102} (1975), 363--371.

\bibitem[Sh:a]{Sh:a} 
Saharon Shelah,
\newblock {\em Classification theory and the number of nonisomorphic models},
\newblock Studies in Logic and the Foundations of Mathematics, vol. 92, North-Holland Publishing Co., Amsterdam-New York, xvi+544 pp, 1978.


\bibitem[Sh:b]{Sh:b}
-----------,
\newblock {\em Proper forcing},
\newblock Lecture Notes in Mathematics, vol. 940, Springer-Verlag, Berlin-New York, xxix+496 pp, 1982.

\bibitem[Sh:e]{Sh:e}
-----------,
\newblock {\em Non-structure theorem},
\newblock Oxford University Press, accepted.


\bibitem[Sh:f]{Sh:f}
-----------,
\newblock {\em Proper and improper forcing},
\newblock Perspectives in Mathematical Logic, Springer, 1998.

\bibitem[Sh:g]{Sh:g}
-----------,
\newblock {\em Cardinal Arithmetic},
\newblock Oxford Logic Guides, vol. 29, Oxford University Press, 1994.


\bibitem[Sh:E58]{Sh:E58}
-----------,
\newblock {\em Existence of endo-rigid Boolean Algebras}.


\bibitem[Sh:E59]{Sh:E59}
-----------,
\newblock {\em General non-structure theory and constructing from linear 
orders},
\newblock  arxiv:1011.3576.


\bibitem[DvSh:65]{DvSh:65}
Keith Devlin and  Saharon Shelah,
\newblock {\em A weak version of $\diamondsuit $ which follows from $2^{\aleph _{0}}<2^{\aleph _{1}}$},
\newblock Israel Journal of Mathematics {\bf 29} (1978), 239--247.

\bibitem[Sh:108]{Sh:108}
Saharon Shelah,
\newblock {\em On successors of singular cardinals}, Logic Colloquium '78 (Mons, 1978,
\newblock Stud. Logic Foundations Math, vol. 97, North-Holland, Amsterdam-New York, 1979, ps. 367--380.


\bibitem[Sh:110]{Sh:110}
------------, 
\newblock Better quasi-orders for uncountable cardinals,
\newblock Israel Journal of Mathematics {\bf 42} (1982), 177--226.


\bibitem[RuSh:117]{RuSh:117}
Matatyahu Rubin and  Matatyahu and Saharon Shelah,
\newblock Combinatorial problems on trees: partitions, $\Delta$-systems and large free subtree,
\newblock Annals of Pure and Applied Logic {\bf 33} (1987), 43--81.



\bibitem[Sh:136]{Sh:136}
Saharon Shelah,
\newblock Constructions of many complicated uncountable structures and Boolean algebras,
\newblock Israel Journal of Mathematics {\bf 45} (1983), 100--146.


\bibitem[DrSh:195]{DrSh:195}
Manfred Droste, Manfred and Saharon  Shelah,
\newblock A construction of all normal subgroup lattices of $2$-transitive automorphism groups of 
linearly ordered sets,
\newblock Israel Journal of Mathematics {\bf 51} (1985), 223--261.


\bibitem[Sh:220]{Sh:220}
Saharon Shelah,
\newblock Existence of many $L_ {\infty,\lambda}$-equivalent, nonisomorphic models of $T$ of power $\lambda$,
\newblock  Annals of Pure and Applied Logic {\bf 34} (1987) 291--310,
\newblock Proceeding of the Model Theory Conference, Trento, June 1986.



\bibitem[Sh:247]{Sh:247}
-----------,
\newblock {\em More on stationary coding}, Around classification theory of models,
\newblock Lecture Notes in Mathematics 1182,  Springer, Berlin, 1986, 224--246.

\bibitem[Sh:309]{Sh:309}
-----------, 
\newblock {\em Black Boxes}.
\newblock arxiv:0812.0656


\bibitem[Sh:331]{Sh:331}
-----------,
\newblock {\em A complicated family of members of tress with $ \omega +1 $ levels}.
\newblock arxiv:math.LO/1404.2414


\bibitem[Sh:420]{Sh:420}
-----------,
\newblock {\em Advances in Cardinal Arithmetic}, 
\newblock in: Finite and Infinite Combinatorics in Sets and Logic,
\newblock Kluwer Academic Publishers, 1993, N.W. Sauer et al (eds.),  
pp. 355--383.  
\newblock arxiv:0708.1979.


\bibitem[Sh:430]{Sh:430}
-----------,
\newblock {\em  Further cardinal arithmetic},
\newblock Israel Journal of Mathematics {\bf 95} (1996), 61--114.


\bibitem[Sh:482]{Sh:482}
-----------,
\newblock {\em Compactness in ZFC of the Quantifier on ``Complete embedding
	of BA's''}.

\bibitem[Sh:511]{Sh:511}
-----------, 
\newblock {\em Building complicated index models and Boolean algebras}.


\bibitem[Sh:620]{Sh:620}
-----------,
\newblock {\em Special Subsets of ${}^{{\rm cf}(\mu)}\mu$, Boolean Algebras
          and Maharam measure Algebras},
\newblock Topology and its Applications {\bf 99} (1999), 135--235.


\bibitem[DrSh:743]{DrSh:743}
Manfred Droste and Saharon Shelah,
\newblock {\em Outer automorphism groups of ordered permutation groups},
\newblock Forum Mathematicum {\bf 14} (2002), 605--621.


\bibitem[Sh:1019]{Sh:1019}
Saharon Shelah,
\newblock {\em Model theory for a compact cardinal}
\newblock Israel Journal of Mathematics, accepted,
arxiv:math.LO/1303.5247

\end{thebibliography}

\end{document}